\documentclass[11pt]{article}

\usepackage[T1]{fontenc} 
\usepackage[utf8]{inputenc} 
\usepackage{authblk} 
\usepackage{amssymb} 
\usepackage{amsmath} 
\usepackage{amsthm} 
\usepackage{booktabs,array} 
\usepackage[margin=1in]{geometry}  
\usepackage[dvipsnames]{xcolor} 
\usepackage{graphicx}
\usepackage{hyperref}
\usepackage{enumitem}
\allowdisplaybreaks[3] 

\usepackage{combelow} 
\usepackage{tikz}
\usepackage{tikz-cd}
\tikzstyle{every picture}=[level distance = 8mm, baseline=-0.5ex]
\tikzstyle{prop}=[shape=circle,minimum size=6mm, draw=black!80, fill=green!30]
\usetikzlibrary{arrows,matrix}
\usetikzlibrary{decorations.pathmorphing}
\newtheorem{thm}{Theorem}[section]
\newtheorem{cor}[thm]{Corollary}
\newtheorem{lem}[thm]{Lemma}
\newtheorem{prop}[thm]{Proposition}

\theoremstyle{definition}
\newtheorem{defn}[thm]{Definition}
\newtheorem{rk}[thm]{Remark}
\newenvironment{remark}{\begin{rk}}{\hfill$\Box$\par\vspace{.1cm}\end{rk}}
\newtheorem{rks}[thm]{Remarks}

\newtheorem{ex}[thm]{Example}
\newenvironment{example}{\begin{ex}}{\hfill$\Box$\par\vspace{.1cm}\end{ex}}
\newtheorem{exs}[thm]{Examples}
\newenvironment{examples}{\begin{exs}}{\hfill$\Box$\par\vspace{.1cm}\end{exs}}

\newcommand{\R}{\mathbb{R}}

\newcommand{\id}{\mathrm{id}}
\newcommand{\Id}{\mathrm{Id}}
\newcommand{\Sym}{\mathrm{Sym}}
\newcommand{\Aut}{\mathrm{Aut}}
\newcommand{\End}{\mathrm{End}}
\newcommand{\Hom}{\mathrm{Hom}}

\newcommand{\calG}{\mathcal{G}} 
\renewcommand{\ss}{\mathbf{s}} 
\renewcommand{\tt}{\mathbf{t}} 
\newcommand{\mm}{\mathbf{m}} 
\newcommand{\uu}{\mathbf{u}} 
\newcommand{\ii}{\mathbf{i}} 
\renewcommand{\aa}{\mathbf{a}} 
\newcommand{\toto}{\rightrightarrows} 
\newcommand{\Pair}{\mathrm{Pair}} 
\newcommand{\GL}{\mathrm{GL}} 
\newcommand{\Diag}{\mathrm{Diag}} 
\newcommand{\calU}{\mathcal{U}} 
\newcommand{\calV}{\mathcal{V}} 
\newcommand{\Par}{\mathcal{P}} 
\newcommand{\Jac}{\mathrm{Jac}} 
\newcommand{\Exp}{\mathrm{Exp}} 
\newcommand{\dircon}{\Gamma}

\newcommand{\loc}{\mathrm{loc}}
\newcommand{\locmor}{\ \circ\hspace{-5pt}\rightarrow}
\newcommand{\locmap}{\ \ast \hspace{-6pt}\rightarrow}
\newcommand{\longlocmor}{\ \circ\hspace{-6pt}\longrightarrow}
\newcommand{\longlocmap}{\ \ast \hspace{-7pt}\longrightarrow}
\newcommand{\simdown}{\mathrel{\vcenter{\hbox{$\sim$}\vskip-.8ex}}}
\newcommand{\locdif}{\overset{\simdown}{\ast\hspace{-5pt}\rightarrow}}

\newcommand{\Bis}{\mathcal{B}}
\newcommand{\Diff}{\mathrm{Diff}}

\newcommand{\X}{\mathfrak{X}}
\newcommand{\Lie}{A}
\newcommand{\Tlin}{{T^\mathrm{lin}}}
\newcommand{\Der}{\mathrm{Der}}
\newcommand{\Ad}{\mathrm{Ad}}

\newcommand{\GLnd}{GL^n_d} 
\newcommand{\GLzerod}{GL^0_d} 
\newcommand{\GLund}{GL^1_d} 

\newcommand{\Tay}{\mathrm{Tay}}

\newcommand{\Rem}{\mathrm{Rem}}
\newcommand{\TDef}{a}
\newcommand{\PDef}{b}
\newcommand{\Coef}{c}
\newcommand{\calO}{\mathcal{O}}
\newcommand{\calD}{\mathcal{D}}

\DeclareFontFamily{U}{mathb}{\hyphenchar\font45}
\DeclareFontShape{U}{mathb}{m}{n}{<-> mathb10}{} 
\DeclareSymbolFont{mathb}{U}{mathb}{m}{n}
\DeclareMathSymbol{\lefttorightarrow}{\mathrel}{mathb}{"FC}
\DeclareMathSymbol{\righttoleftarrow}{\mathrel}{mathb}{"FD}

\newcommand{\JetSpace}{
\parbox[c]{6cm}{  
\begin{picture}(5,2.4)(-1.1,-2)
\put(0,0){$M$}
\put(3,0){$M'$}
\put(0.7,0.1){\vector(1,0){2}}
\put(0.55,0){$\ast$}
\put(1.5,0.3){$f$}
\put(0,-1.7){$\R^d$}
\put(0.2,-1.2){\vector(0,1){1}}
\put(0.1,-1.35){$\ast$}
\put(-0.3,-0.8){$\varphi$}
\put(0.6,-1.2){\vector(2,1){2.2}}
\put(0.45,-1.35){$\ast$}
\put(1.5,-0.5){$\bar{f}$}
\put(-.7,-1.7){\Large $\lefttorightarrow$}
\put(-.2,-1.6){$\ast$}
\put(-1.1,-1.7){$h$}
\put(3,-1.7){$\R^{d'}$}
\put(3.2,-1.2){\vector(0,1){1}}
\put(3.1,-1.35){$\ast$}
\put(3.3,-0.8){$\varphi'$}
\put(0.7,-1.5){\vector(1,0){2}}
\put(0.6,-1.6){$\ast$}
\put(1.5,-1.3){$\tilde{f}$}
\put(3.6,-1.7){\Large $\righttoleftarrow$}
\put(3.5,-1.7){$\ast$}
\put(4.3,-1.7){$h'$}
\end{picture}}
}

\begin{document}

\title{Direct connections on jet groupoids}

\author{Sara Azzali, Youness Bouta\"ib, Alessandra Frabetti, Sylvie Paycha}


\maketitle

{\it This paper is dedicated to the memory of the late Kirill Mackenzie (1951-2020)}

\begin{abstract}
Direct connections on Lie groupoids are smooth local right inverses of the anchor map which preserve the units. They were introduced by N.~Teleman and J.~Kubarski as generalised parallel transports on a frame groupoid. 
A direct connection on a Lie groupoid (over a connected base manifold) forces it to be a gauge groupoid and induces an infinitesimal connection on its Lie algebroid. We show that the curvature of a direct connection, which measures the multiplicative defect, contain symmetric data that are not detected by the curvature of its infinitesimal connection.
Then we study jet prolongation of direct connections, by introducing the notion of jet primitive bisection. 
We exhibit three natural examples for the pair groupoid, built from Jacobi fields of geodesic variations, one of which recovers the exponential jet connection of I.~Kol\'a\v{r} and W.~M.~Mikulski.
Finally, we use direct connections to define a covariant Taylor expansion of sections of a vector bundle and a higher-order direct connection on the frame groupoid of its jet bundle which is not a jet prolongation, and show that this pair gives a local model for jet regularity structures, as defined by M. Hairer and collaborators. 
\end{abstract}

\noindent \textbf{2020 Mathematics Subject Classification:} Primary 58H05; Secondary 53C05, 58A20, 53C22. 
\bigskip

\noindent \textbf{Keywords:} Lie groupoids, Lie algebroids, Connections, Jet prolongations, Jacobi fields, Taylor expansion, Transport maps. 

\tableofcontents  


\section{Introduction} 
\label{sec:intro}

Groupoids $\calG\toto M$ are a generalisation of groups consisting of arrows attached to a source and a target in the base manifold $M$, with two key properties: 
the arrows with equal source and target
form a group above each base point, the others can be composed whenever the target of the first
coincides with the source of the second, and they are invertible.
Groupoids were introduced by H.~Brandt in 1926 \cite{Brandt:1926} and endowed with topological and differentiable structures by C.~Ehresmann \cite{Ehresmann:1952}. They are nowadays ubiquitous in mathematics, physics and computer science, see for instance the applications mentioned in \cite{Weinstein:1996, Winkelkemper:1983, Altenkirch:2021}.

Lie groupoids, and in particular gauge groupoids, naturally act on fibre bundles and are useful to encode the \emph{global symmetries} of a bundle: the arrows above single points encode the \emph{internal} or \emph{gauge} symmetries acting along the fibres, those above distinct points also encode the \emph{base manifold} symmetries \cite{Costa-Forger-Pegas:2018,Costa-Forger-Pegas:2021}. 
A connection on a fibre bundle turns such symmetries into an actual comparison of nearby fibres, by prescribing a \emph{parallel transport} between them, that is, a specific arrow of the symmetry groupoid with given source and target points. 
In the classical setting of a principal bundle, a \emph{principal connection} selects, among the gauge symmetries, those that realize the parallel transport along the paths of the base and is described infinitesimally by a Lie algebra valued differential form. 
\emph{Cartan connections} extend principal connections by allowing the parallel transports to be governed by the base space symmetries \cite{Cartan:1923}, thus giving a geometric model suitable to describe General Relativity \cite{Wise:2010,Krasnov-Percacci:2018}. 

These theories of connections are \emph{infinitesimal}, in the sense that the parallel transport between fibres is obtained by integrating the connection $1$-forms (which are in fact linked to infinitesimal symmetries) along chosen paths on the base manifold.  
In the groupoid setting, the infinitesimal symmetries form a Lie algebroid, a fibred generalisation of a Lie algebra whose sections generalise the vector fields on the base.  The corresponding infinitesimal notion of connection, due to C.~Ehresmann in 1952 \cite{Ehresmann:1952} and enhanced by S. Kobayashi in \cite{Kobayashi:1957}, is a section of its anchor map.
Integrating these infinitesimal data to a parallel transport living in the groupoid led to the search for a suitable notion of connection on groupoids. 

The first notion of connection on a Lie groupoid is that of a \emph{path connection} introduced by J.~Virsik in 1971 \cite[Appendix]{Virsik:1971} and thoroughly studied by K.~Mackenzie \cite[Chapter 6]{Mackenzie:2005}. It is based on a path lifting and obtained by integration along a path on the base manifold.
Another infinitesimal type of connection is the \emph{multiplicative Ehresmann connection} recently introduced by R.~Loja Fernandes and I.~M\u{a}rcu\cb{t} in \cite{LojaFernandes-Marcut:2023} (and in previous works by other authors cited therein). 
Unlike path connections, by construction a multiplicative Ehresmann connection is compatible with the composition and it is encoded by an intrinsic geometric datum on the groupoid. 

Path-independent notions of connections on a Lie groupoid, where the transport between fibres does not require integration, were first introduced by A.~Kock in the framework of Synthetic Differential Geometry \cite{Kock:1989,Kock:2007,Kock:2017}. They were further developped by  N.~Teleman and J.~Kubarski for the frame groupoid of a vector bundle, to capture the geometric content of its Chern classes, who called them \emph{linear direct connections} \cite{Teleman:2004,Teleman:2007,Kubarski-Teleman:2006}. The generalisation to transitive Lie groupoids was announced by Kubarski in \cite{Kubarski:2008}, we call it \emph{direct connection}.

A direct connection on a Lie groupoid $\calG\toto M$ is a smooth local right inverse $\dircon:\Pair(M)\longlocmap \calG$ of the anchor map that preserves the units (Definition \ref{def:directconnection}). 
It is \emph{flat} if it preserves the composition, and \emph{transitive} if it preserves the composition along the geodesics. 
The basic example is the parallel transport of a linear connection on a vector bundle $E$ along the small geodesic from a source to a nearby target, where the base manifold $M$ is equipped with a torsion-free linear connection. This is a transitive direct connection which is not necessarily flat.

Unlike path connections, a general direct connection does not require the choice of a path. Unlike multiplicative Ehresmann connections, it does not require the compatibility with the groupoid composition. As a result, direct connections carry geometric information that remains invisible at the infinitesimal level.
In other words, while the category of gauge groupoids is equivalent to that of principal bundles, the category of \emph{gauge groupoids with direct connections} is much larger than that of \emph{principal bundles with principal connections}. 
Moreover, as Teleman observed and in contrast with path connections, a direct connection need not be smooth, so that it can also be defined  on singular or non-smooth base manifolds \cite{Kitada:2013} (see also \cite{Emerson-Meyer:2010}).

In this paper we study smooth direct connections and their jet prolongation to jet groupoids. Our original motivation stems from M. Hairer’s theory of regularity structures \cite{Hairer:2014a}, where transport maps $\Gamma(x,y)$ appear naturally. Observing that these maps compose as arrows of a groupoid, namely $\Gamma (x,y)\Gamma(y,z)= \Gamma(x,z)$, lead us to investigate direct connections for frame groupoids of jet bundles.

Let us state our main results. A preliminary result (Proposition \ref{prop:direct}) used all along the paper is that a Lie groupoid which admits a direct connection is a gauge groupoid.
\begin{itemize}
\item Our first set of results concerns geometric by-products of direct connections on gauge groupoids. We show that the latter give rise to both infinitesimal connections on their associated Lie algebroids and to usual linear connections on vector bundles on which the groupoid acts. Among direct connections, we  characterise the usual parallel transports as the transitive ones, that is, those that are groupoid morphisms along geodesics (Theorem \ref{thm:transitive = parallel transport}). 
We further show (Theorem \ref{thm:curvature}) that the curvature of a direct connection encodes symmetric data which are not detected by the curvature of the underlying infinitesimal  connection. 

\item The second set of results concerns jet prolongation of direct connections on jet groupoids. 
We characterise the jet prolongation of both direct connections $\dircon^n = J^n\dircon\circ \Delta^n$ and their infinitesimal connections $\delta^{\dircon^n} = J^n \delta^\dircon \circ \delta^{\Delta^n}$ in terms of direct connections $\Delta^n$ on the jet pair groupoid (Theorem \ref{thm:jet prolongation dircon}). 
Then we construct examples of such direct connections, generalising the usual parallel transport of a linear connection on the base manifold to new transport maps that are solutions of Jacobi second-order geodesic variation equations (Theorems \ref{thm:parallel direct connection}, \ref{thm:Jacobi direct connection} and \ref{thm:exponential direct connection}). 
With the same tools, one can read solutions of higher order geodesic variations as transport maps along the fibres of tangent bundles. 

\item The third set of results concerns a covariant Taylor expansion of sections of a vector bundle built from a direct connection. We show that, for general direct connections $\dircon$, both the principal part and the remainder of the Taylor expansion can be expressed in terms of the usual higher covariant derivatives of the infinitesimal connection $\nabla$ associated to $\dircon$ (Theorem \ref{thm:Taylor with nabla}), but the coefficients contain new terms recording non trivial derivatives of the curvature of $\dircon$ at the base point and along geodesics (the \emph{transitivity defects}). 
Finally, we use the covariant Taylor expansion to construct a direct connection on the gauge groupoid of a jet bundle which is not a jet prolongation (Theorem \ref{thm:Taylor direct connection}) and gives rise to a local model for $E$-valued jet regularity structures (Proposition \ref{prop:regularity structure}). 
\end{itemize}

The paper is organised as follows.
Section \ref{sec:groupoids} recalls the necessary background on Lie groupoids, bisections, groupoid actions and Lie algebroids, with an emphasis on gauge and frame groupoids and a systematic treatment of local maps and local bisections.

In Section \ref{sec:direct-connections} we develop the general theory of direct connections on Lie groupoids.
We show that a Lie groupoid over a connected base which admits a direct connection is necessarily a gauge groupoid (Proposition \ref{prop:direct}),  that every direct connection induces an infinitesimal connection on the induced Lie algebroid (Proposition \ref{prop:dircttoinf}) and hence a linear connection on associated vector bundles (Corollary \ref{cor:deltaIsoE}), and characterise the usual parallel transports as the transitive direct connections (Proposition \ref{prop:tau=dir conn} and Theorem \ref{thm:transitive = parallel transport}). 
We then introduce the curvature of a direct connection, which measures its defect from being a groupoid morphism, and show that its antisymmetric part reproduces the usual curvature tensor of the induced infinitesimal connection (Theorem \ref{thm:curvature}), while its full expression carries finer, non-infinitesimal invariants. 
We close the section by characterizing the Lie groupoids admitting a flat direct connection as precisely the flat ones (Proposition \ref{prop:flat groupoid}), largely studied by K.~Mackenzie \cite{Mackenzie:1989}.

In Section \ref{sec:jetconnections} we study the jet prolongations of direct connections, building on the jet prolongation of groupoids and Lie algebroids after C.~Ehresmann \cite{Ehresmann:1955} and I.~Kol\'a\v{r} \cite{Kolar:2007}.
Using the notion of \emph{jet primitive bisection} (Definition \ref{def:jet primitive bisection}), we prove that the jet prolongation of a direct connection is determined by that of the trivial direct connection $\Id$ on the pair groupoid, and compute its associated infinitesimal connection (Theorem \ref{thm:jet prolongation dircon}).
We then exhibit three explicit jet prolongations of $\Id$ on $\Pair(M)$, arising from three Jacobi fields of geodesic variations that differ by the initial value of the field and of its covariant derivative (Theorems \ref{thm:parallel direct connection}, \ref{thm:Jacobi direct connection} and \ref{thm:exponential direct connection}). We also compute their infinitesimal connections and show that the first two are, in a sense, orthogonal to each other, while the third reproduces the exponential jet connection on $J^nTM$ of I.~Kol\'a\v{r} \cite{Kolar:2007} and W.~M.~Mikulski \cite{Mikulski:2007}.

Section \ref{sec:direct connections on jet bundles} is motivated by the theory of regularity structures for stochastic PDEs, introduced by M.~Hairer \cite{Hairer:2014a,Hairer:2014b} and further developed on manifolds by A.~Dahlqvist, J.~Diehl and B.~Driver \cite{DahlqvistDiehl-Driver:2019} and on vector bundles by M.~Hairer and H.~Singh \cite{Hairer-Singh:2023}.
Starting from a direct connection $\dircon$ on the frame groupoid of a vector bundle $E$, we introduce a covariant Taylor expansion of local sections of $E$ whose displacement along the fibres is ruled by $\dircon$ (Definition \ref{def:Taylor}).
In Theorem \ref{thm:Taylor with nabla} we show that both the principal part of this Taylor expansion and the remainder can be given in terms of higher covariant derivatives of the linear connection induced by $\dircon$, with coefficients containing new terms, called \emph{transitivity defects}, related to the derivatives of the curvature of $\dircon$ at the base point (cf. Definition \ref{def:transitivity defects}). The transitivity defects vanish precisely when $\dircon$ is the parallel transport of a linear connection, so they are a genuinely new feature of general direct connections.
Using this expansion we build a higher-order direct connection $\Sigma^n$ on the frame groupoid of a $n$-jet bundle which is not a jet prolongation of a direct connection on $E$ (Theorem \ref{thm:Taylor direct connection}), and we prove that the pair of maps $(\Tay^n,\Sigma^n)$ provides a covariant local model for $E$-valued $n$-jet regularity structures of transport precision $n+1$ (Proposition \ref{prop:regularity structure}). 

Finally, an appendix (Section \ref{sec:jets}) collects the background on jets and on the jet prolongation of principal and vector bundles.
\bigskip 

\noindent 
\textbf{Acknowledgements.}
The authors are very grateful to the late Kirill Mackenzie and thank Alexander Schmeding for their very helpful comments on a preliminary version of the paper. S. Azzali gratefully acknowledges the support of the Istituto Nazionale di Alta Matematica (INdAM), the SPP 2026 "Geometry at Infinity" (first funding period), and the COST Action CA21109.
Y. Bouta\"ib thanks the DFG for the financial support within the research unit FOR 2402.
A. Frabetti warmly thanks R. Loja Fernandes for correcting an earlier misinterpretation of the gauge group in terms of the gauge groupoid. 
\bigskip 

\noindent 
\textbf{AI disclosure statement.}
The project of this paper started in 2016 and a preliminary version of it circulated from 2022. AI assistance (Claude and ChatGPT6-Astra) was used from August 2026 at the very last step of the writing of the new completely revisited version, to double-check references and for proofreading. 
Claude detected a mistake in Section \ref{subsec:Higher-order dircon}, that was entirely corrected by the authors. 
ChatGPT suggested Example \ref{ex:dircon not jet prolongation} and the direct connections of Remark \ref{rem:interpolation}. He also suggested shorter and more readable proofs for the theorems in Section \ref{sec:direct connections on jet bundles}. 
Except for the AI assistance described above, the authors developed the arguments and wrote the manuscript. The authors take full responsibility for the content of the paper.


\section{Lie groupoids and Lie algebroids} 
\label{sec:groupoids}

In this section we recall basic facts and fix our notations on Lie groupoids, their bisections, their natural actions on fibre bundles, and on their infinitesimal counterpart (Lie algebroids). 
For the purpose of the present work, we only need to focus on gauge groupoids arising from principal bundles and on the frame groupoids of the associated vector bundles.
The only original part of this section is the systematic treatment of local maps between groupoids, local bisections, and local actions on vector bundles. 
Our main references are the standard textbooks by K.~Mackenzie \cite{Mackenzie:2005}, I.~Moerdijk and J.~Mr\v{c}un \cite{Moerdijk-Mrcun:2003} and E.~Meinrencken \cite{Meinrenken:2017}, and the pedagogical introduction by A.~Kumpera \cite{Kumpera:2015}.
 

\subsection{Lie and gauge groupoids}\label{subsec:groupoids}


Let $M$ be a smooth manifold. 
A smooth \textbf{groupoid} on $M$ is a manifold $\calG$ together with the following smooth structure maps:  
\begin{enumerate}[label={(\roman*)}]
\item a \textbf{source} $\ss:\calG\longrightarrow M$ and a \textbf{target} $\tt:\calG\longrightarrow M$,                                                                   
\item an associative \textbf{multiplication} (or \textbf{composition}) $\mm:\calG\times_M\calG\longrightarrow \calG,\ (g_2,g_1)\mapsto g_2 g_1$, defined on the set of composable elements $\calG\times_M\calG= \{(g_2,g_1),\ \ss(g_2)=\tt(g_1)\}$, 

\item a \textbf{unit} $\uu: M\to \calG,\ x\mapsto \uu(x)=:1_x$ with $\tt(1_x)=\ss(1_x)=x$ such that $1_{\tt(g)}\, g=g=g\, 1_{\ss(g)}$ for any $g\in \calG$,
  
\item an \textbf{inversion} $\ii:\calG\longrightarrow \calG,\ g\mapsto g^{-1}$ with $\ss(g^{-1})=\tt(g)$ and $\tt(g^{-1})=\ss(g)$, such that $g\, g^{-1}= 1_{\tt(g)}$ and $g^{-1}\, g= 1_{\ss(g)}$.  
\end{enumerate}
The induced map $\aa=(\tt,\ss):\calG\longrightarrow M\times M$ is called the \textbf{anchor}. From the axioms it follows that the source and the target are surjective maps, the unit is injective and the inversion is bijective. 
The manifold $M$ is called the \textbf{base} of the groupoid and can be identified with the \textbf{set of units} $\uu(M)\subset \calG$. 
A groupoid is compactly denoted by $\calG\toto M$ and the structure maps $(\ss,\tt,\mm,\uu,\ii)$ are tacitely understood. 

A group $G$ can be seen as a groupoid $G\toto \ast$ over a single point, with trivial source and target maps. Hence groupoids generalise groups. 

Given a  groupoid $\calG\toto M$ and points $x,y$ in $M$, we  use the following notations:
\begin{itemize}
\item $\calG^x:=\ss^{-1}(x)$ for the \textbf{$\ss$-fibre} of $x$, with restricted target map $\tt^x:= \tt_{|\calG^x}:\calG^x\to M$, 
\item $\calG_y:=\tt^{-1}(y)$ for the \textbf{$\tt$-fibre} of $y$, with restricted source map $\ss_y:= \ss_{|\calG_y}:\calG_y\to M$, 
\item $\calG^x_y:=\calG^x\cap \calG_y$ for the \textbf{(anchor) fibre} of $(y,x)$, whose elements are sometimes denoted by $g_y^x$. 
\end{itemize}
Similarly, for $U,V\subset M$, we set $\calG^U:= \ss^{-1}(U)$, $\calG_V:= \tt^{-1}(V)$ and $\calG^{U }_V:= \calG^U\cap \calG_V$.

For any $x\in M$, the set $\calG_x^x$ is a (non-empty) group with the composition and unit $1_x$, called the \textbf{vertex group} (or the \textbf{isotropy}) at $x$.
The non emptyness of the set $\calG_y^x$ above two distinct points $x,y$ in $M$ defines an equivalence relation ($x\sim y$ if and only if $\calG_y^x$ is non empty), whose equivalence classes are called \textbf{orbits} of $\calG$.
The orbit of a point $x$ in $M$ is the set $\mathcal{O}_x=\tt^x(\calG^x)=\ss_x(\calG_x)\subset M$.
The \textbf{orbit space} of $\calG$, denoted $M/\calG$, is the quotient of $M$ by the relation and gives a foliation of $M$ which is possibly singular.
The groupoid $\calG$ is \textbf{regular} if the orbits have all the same dimension, that is, if the foliation is regular. The groupoid is \textbf{transitive} if it has a single orbit $M$, which happens if and only if the anchor map is surjective.  
\bigskip

A smooth groupoid $\calG\toto M$ is a \textbf{Lie groupoid} if the source and the target are surjective submersions. 
This guarantees the following nice properties \cite[Theorem 5.4]{Moerdijk-Mrcun:2003}, \cite[\S 1.5]{Mackenzie:2005}: 
\begin{enumerate}[label={\arabic*.}]
\item The units $\uu(M)\subset \calG$ form a submanifold and the composable elements $\calG\times_M \calG$ form a manifold. 
\item The fibre $\calG^x_y$ above each pair of points is a closed submanifold of $\calG$ (possibly empty). 
\item The vertex group $\calG_x^x$ above each point is a Lie group\footnote{
However, the vertex groups over distinct points $x$ and $y$ are not necessarily isomorphic, since the fibre $\calG^x_y$ can be empty and the vertex groups can belong to separate connected components of $\calG$, even if $M$ is connected.}.
\item The orbit $\mathcal{O}_x$ of each point is an immersed submanifold of $M$ and the restricted maps $\tt^x$ and $\ss_x$ are both principal $\calG_x^x$-bundles on $\mathcal{O}_x$.
\item If $\calG$ is transitive then the anchor map is a surjective submersion. 
\end{enumerate}
Thanks to these properties, Lie groupoids admit an infinitesimal calculus (via Lie algebroids) analogous to that defined on Lie groups (via Lie algebras), and are suitable to study smooth actions on fibre bundles and the related symmetries. 

A groupoid can be seen as a small category whose objects are the points of the base manifold and all morphisms are invertible.
A morphism between two groupoids is a functor between the (category-theoretic) groupoids. 
We focus on morphisms between two Lie groupoids $\calG$ and $\calG'$ over the same base $M$, called \textbf{morphisms over $M$}. 
These are smooth maps $\phi: \calG\to \calG'$ such that
\begin{enumerate}[label={(\roman*)}]
 \item $\ss'\circ \phi= \ss$ and $\tt'\circ \phi= \tt$, 
\item $\phi\circ \uu= \uu'$, 
\item $\phi(g\,g')= \phi(g)\, \phi(g')$ for any composable $g,g'$ in $\calG$, and therefore $\phi\circ \ii = \ii'\circ \phi$. 
\end{enumerate} 
A morphism $\phi$ is an \textbf{isomorphism} if the map $\phi$ is a diffeomorphism. 
A \textbf{subgroupoid} of a Lie groupoid $\calG\toto M$ is a groupoid $\calG'\toto M$ together with an injective Lie groupoid morphism $\iota:\calG'\hookrightarrow \calG$ over $M$, giving the inclusion. 

\begin{example}[Trivial Groupoid] \label{ex:trivialGgroupoid}
Given a Lie group $G$, the cartesian product $M\times G\times M$ defines a Lie groupoid on $M$, called \textbf{trivial groupoid with vertex group $G$}, with source $\ss(y,g,x)=x$, target $\tt(y,g,x)=y$ and composition $(z,h,y)(y,g,x)=(z,hg,x)$. Therefore the unit is $1_x=(x,1,x)$ and the inverse is $(y,g,x)^{-1}=(x,g^{-1},y)$.
\end{example}

A Lie groupoid $\calG\toto M$ is called \textbf{trivial} if it is isomorphic to a trivial groupoid. Any trivial Lie groupoid is clearly transitive.

\begin{example}[Pair groupoid] \label{ex:pair}
The \textbf{pair groupoid} of $M$ is the trivial groupoid with trivial vertex group $G=\{1\}$, namely the cartesian product $\Pair(M)=M\times M$, whose elements $(y,x)$ have source $x$, target $y$, composition $(z,y)\,(y,x)=(z,x)$, unit $1_x=(x,x)$ and inverse $(y,x)^{-1}=(x,y)$.
The set of units of the pair groupoid coincides with the \textbf{diagonal of $M$}, $\uu(\Pair(M))=\Diag(M)=\{(x,x),\ x\in M\}$. 
\end{example}

For any Lie groupoid $\calG\toto M$, the anchor map $\aa =(\tt,\ss):\calG\longrightarrow \Pair(M)$ is a Lie groupoid morphism over $M$, which maps the units of $\calG$ onto the units $\Diag(M)$ of $\Pair(M)$, but is not globally surjective in general. 

\begin{example}[Gauge groupoid] \label{ex:gaugegroupoid}
Let $\pi:P\to M$ be a principal bundle with structure group $G$ acting on the fibres from the right, transitively and without fixed points.
The \textbf{gauge} or \textbf{Atiyah groupoid of $P$} is the quotient manifold 
\begin{align*}
\calG(P) = P\times_GP := (P\times P)/\sim
\end{align*}  
under the equivalence relation $(p,q)\sim (p\,g,q\,g)$ for all $p,q\in P$ and $g\in G$, whose elements are equivalence class of pairs $(p,q)$ in $ P\times P$ denoted $[p,q]$.
The source and target maps are given by the bundle projection, namely $\ss([p,q])=\pi(q)$ and $\tt([p,q])=\pi(p)$ for all $p,q\in P$.
The partial composition is defined by $[p,q]\,[p',q']=[p,q'g]$ only if $\pi(q)=\pi(p')$, where $g$ is the unique element in $G$ such that $q=p'\,g$.
The units are $1_x = [p,p]$ for any $p\in P_x=\pi^{-1}(x)$, and the inverse of $[p,q]$ is $[q,p]$. 

For instance, the trivial groupoid $M\times G\times M$ is the gauge groupoid of the trivial principal bundle $P= M\times G$. This holds in particular for the pair groupoid $\Pair(M)$ whose structure group is trivial. 

Gauge groupoids are characterised by a very simple local structure. In fact, if the base manifold $M$ is connected, they are equivalently determined among Lie groupoids $\calG$ by any of the following conditions \cite[Lemma 2]{Kumpera:2015} \cite[Propositions 1.3.3 and 1.3.5]{Mackenzie:2005} \cite[Theorem 3.10]{Meinrenken:2017}: 
\begin{enumerate}[label={\arabic*.}]
\item $\calG$ is \textbf{transitive}, i.e. there is a single orbit $M$ (hence the anchor is a surjective submersion). 

\item $\calG$ is \textbf{locally trivial}, i.e. for any $x$ in $ M$ there exists an open neighbourhood $U$ of $x$ in $M$ such that $\calG_U^U$ is isomorphic to the trivial groupoid $U\times \calG_x^x\times U$.  

\item $\calG$ admits a \textbf{section atlas}, i.e. there exists a point $x_0 \in M$ and a collection of local sections of $\tt^{x_0}:\calG^{x_0}\to M$, that is, an open covering $\{U_\alpha\}$ of $M$ and smooth maps $\tau_\alpha: U_\alpha \to \calG^{x_0}_{U_\alpha}$ such that $\tt^{x_0}\circ\tau_\alpha= \mathrm{Id}_{U_\alpha}$, called \textbf{local decomposing maps}. 
This implies that the restriction $\tt^{x_0}$ of $\tt$ is a surjective submersion (which is not \emph{a priori} ensured in a Lie groupoid), and one can define its \textbf{transition functions} $g_{\alpha \beta}:U_\alpha\cap U_\beta \to \calG_{x_0}^{x_0}$ by $g_{\alpha \beta}(x)=\tau_ \alpha(x)^{-1}\tau_\beta(x)$. 

\item The anchor $\aa:\calG\to M\times M$ is a locally trivial fibre bundle with fibre $\calG_{x_0}^{x_0}$ for any choice of $x_0\in M$, and structure group $\calG_{x_0}^{x_0}\times\calG_{x_0}^{x_0}$ acting on the left on the fibres by $(g',g)\cdot h=g'hg^{-1}$, for any $g',g,h\in \calG_{x_0}^{x_0}$.
  
\item The source $\ss:\calG=\bigcup_{x\in M} \calG^x \to M$ is a locally trivial fibration in principal $G$-bundles on $M$, with $G=\calG_{x_0}^{x_0}$ for any choice of $x_0\in M$. The principal $\calG_{x}^{x}$-bundle $\tt^x:\calG^{x}\to M$  is called the \textbf{vertex bundle} of $\calG$ at the point $x$.
\end{enumerate}
\end{example}

\begin{example}[Frame groupoid] \label{ex:framegroupoid}
The \textbf{frame groupoid} of a vector bundle $E$ is the set of linear isomorphisms $\varphi^x_y: E_x\to E_y$ between its fibres
\begin{align*}
\GL(E)=\bigcup_{x,y \in M} \GL\left(E_x,E_y\right) \toto M
\end{align*}
with source $\ss(\varphi^x_y)=x$, target $\tt(\varphi^x_y)=y$ and the usual composition $\eta^y_z\circ \varphi^x_y: E_x\to E_y\to E_z$ of linear maps.
The frame groupoid is isomorphic to the gauge groupoid of the \textbf{frame bundle} $FE\to M$ of $E$. If $E$ has fibre $\R^r$, this is the principal $GL_r(\R)$-bundle with fibre $F_xE= \GL\left(\R^r,E_x\right)$ given by the set of linear isomorphisms $\varphi_x:\R^r\to E_x$ (the \textbf{frames} of $E_x$) and projection $\pi(\varphi_x)=x$. 
The isomorphism
\begin{align*}
\calG(FE)\cong \GL(E) \toto M
\end{align*}
is given by $[\varphi_y,\varphi_x] \mapsto \varphi^x_y = \varphi_y \circ (\varphi_x)^{-1}:E_x\to \R^r\to E_y$. 

If $E=P\times_G \R^r$ is the vector bundle associated to a principal $G$-bundle $P$, through a linear representation of $G$ on $\R^r$, then the structure group $GL_r(\R)$ of $E$ can be reduced to $G$ and $P\to FE$ is a reduction of the frame bundle. This means that each element $p\in P_x$ can be seen as a frame $p:\R^r\to E_x$. Then, there is an associated groupoid map $\calG(P)\to \GL(E)$, which is injective if and only if the representation $G\to GL_r(\R)$ is faithful. 
\end{example}


\subsection{Local maps and local morphisms}
\label{subsec:local maps} 


Being interested in jet prolongations of maps defined on groupoids, we will need to consider maps possibly defined only locally. 
In a fibre bundle, an object is \emph{local} if it is defined in an open neighbourhood of a base point. For a groupoid, this notion must be adapted to the fact that its very essence is to relate distinct base points: \emph{locality} then means that the points to be related are sufficiently close to each other, wherever they are in the base manifold. This leads to the following definition, which extends that of a \emph{local morphism} given in \cite[Definition 6.1.6]{Mackenzie:2005}. 

\begin{defn}\label{def:localmap} 
A \textbf{local map} between two Lie groupoids $\calG$ and $\calG'$ over the same manifold $M$ is a smooth map $\phi: \calU\subset \calG \longrightarrow \calG'$ defined on an open neighbourhood $\calU$ of the units $\uu(M)\subset \calG$, which preserves the source, the target and the units, that is, 
\begin{enumerate}[label={(\roman*)}]
\item $\ss'\circ \phi = \ss$ and $\tt'\circ \phi = \tt$,
\item $\phi\circ \uu= \uu'$. 
\end{enumerate}
If $\calU=\calG$, we call $\phi$ a \textbf{global map}.
To distinguish local from global maps at a glance, we denote local maps by $\phi:\calG\longlocmap \calG'$. 
The local map $\phi$ is a \textbf{local morphism} of Lie groupoids if it also preserves compositions, i.e. 
\begin{enumerate}[label={(\roman*)},resume]
\item $\phi(g'g)=\phi(g')\phi(g)$ for any composable $g,g'$ in $\calU$ whose product $g'g$ lies in $\calU$. 
\end{enumerate}
In this case, it also preserve inversions, i.e. $\phi(g^{-1})=\phi(g)^{-1}$ for all $g\in \calU$ such that $g^{-1}$ lies in $\calU$. 
We denote local morphisms by $\phi: \calG\longlocmor \calG'$, as in \cite[Definition 6.1.6]{Mackenzie:2005}. 
A global morphism in this sense is the same as a groupoid morphism of Section \ref{subsec:groupoids}.
\end{defn}

\begin{example}[Local maps on the pair groupoid] \label{ex:diagonal domain}
Let $\calG\toto M$ be a Lie groupoid. 
A local map $\phi:\Pair(M)\longlocmap \calG$ is defined on an open neighbourhood $\calU$ of the diagonal $\Diag(M)$, that we shall call a \textbf{diagonal domain}%
\footnote{This terminology is borrowed from \cite[Def. 83]{DahlqvistDiehl-Driver:2019}. Diagonal domains are called first neighbourhoods of the diagonal in \cite[\S 1]{Kock:2017}.}
in $M\times M$. 
Since $\phi$ preserves source and target, we have $\aa \circ \phi = \id_\calU$. In other words $\phi$ has a right inverse, therefore it is an embedding. 
Since $\phi(\Diag(M)) = \uu(M)$, the embedded image $\phi\big(\calU\big)$ in $\calG$ contains $\uu(M)$. 
\end{example}

A local map between two Lie groupoids cannot always be extended to a global one, even if it is a local groupoid morphism. Mackenzie proved in \cite[Theorem 6.1.10]{Mackenzie:2005} that this is possible under rather restrictive conditions, namely when $\calG \longlocmor \calG'$ is a local groupoid morphism of locally trivial Lie groupoids over the same base manifold $M$, if the source-fibres of $\calG$ are connected and simply connected (or alternatively for the target-fibres).
In particular, if $M$ is connected and simply connected, any local groupoid morphism $\phi:\Pair(M)\locmor \calG$ over $M$ can be extended to the whole pair groupoid (but of course not any local map $\phi:\Pair(M)\locmap \calG$). 


\subsection{Bisections and gauge transformations} \label{subsec:bisections}

Bisections of Lie groupoids play the role of sections of principal bundles, with the essential feature that they reproduce not only groups of vertical transformations acting along the fibres (the gauge groups), but also groups of horizontal transformations of the base manifold. 
For the purpose of this work, local bisections are the key notion to define jet groupoids. 

Let $\calG\toto M$ be a Lie groupoid.
A \textbf{local bisection} of $\calG$ is a smooth local section $\sigma: M\locmap \calG$ of the source $\ss$ defined on an open subset $U\subset M$, that is, $\ss\circ \sigma = \mathrm{Id}_U$, whose composition with the target $\tt$ is a diffeomorphism $\varphi_\sigma=\tt\circ \sigma$ between $U$ and an open subset $\varphi_\sigma(U)\subset M$ (which might not intersect $U$).
If $U=\varphi_\sigma(U)=M$, $\sigma$ is a \textbf{(global) bisection}. 
The space $\Bis(\calG)$ of global bisections of $\calG$ is a group with the following operations:
\begin{enumerate}[label={(\roman*)}]
\item \textbf{multiplication} of $\sigma$ and $\sigma'$ given by $(\sigma'\bullet \sigma)(x)=\sigma'(\varphi_{\sigma}(x))\ \sigma(x)$ for any $x\in M$, which corresponds to a semidirect product law on pairs $(\varphi_\sigma,\sigma)$,

\item \textbf{unit bisection} given by the unit map $\uu:M\to \calG$ of the groupoid $\calG$,

\item \textbf{inverse} $\sigma^{-1}$ given by $\sigma^{-1}(x)=\big(\sigma(\varphi_\sigma^{-1}(x))\big)^{-1}$ for any $x\in M$, where $\varphi_{\sigma}^{-1}$ denotes the inverse diffeomorphism of $M$ while the external $\big(\ \big)^{-1}$ denotes the inverse in the groupoid $\calG$. 
\end{enumerate}
See \cite{Schmeding-Wockel:2015} (and in particular Theorem 2.8) for the topological properties of the group $\Bis(\calG)$.
The subspace $\Bis_M(\calG)$ of bisections $\sigma:M\to \calG$ such that $\varphi_\sigma=\Id_M$ is a subgroup of $\Bis(\calG)$.
The map
\begin{align}
\label{eq:BisectionsToDiffeos}
\varphi:\Bis(\calG)\to \Diff(M),\ \sigma\mapsto \varphi_\sigma = \tt\circ\sigma 
\end{align}
is a group homomorphism to the group of diffeomorphisms of $M$ (not necessarily surjective), which sends $\Bis_M(\calG)$ to $\{\Id_M\}$. 
In fact, $\sigma'\bullet \sigma=(\sigma'\circ\varphi_\sigma) \sigma$ (seen as a pointwise groupoid multiplication), and therefore
\begin{align*}
\varphi_{\sigma'\bullet \sigma} &= \tt\circ \big((\sigma'\circ\varphi_\sigma) \sigma\big)
= \tt\circ (\sigma'\circ\varphi_\sigma) = (\tt\circ \sigma')\circ\varphi_\sigma
=\varphi_{\sigma'}\circ \varphi_\sigma. 
\end{align*}

The space of local bisections $\Bis_{\loc}(\calG)$ is a pseudo-group with the above operations, because the product $\sigma'\bullet\sigma$ of two local sections defined respectively on two open sets $U'$ and $U$ is defined only if the target space $\varphi_\sigma(U)$ of $\sigma$ has a non-empty intersection with the source space $U'$ of $\sigma'$.
The subspace $\Bis_{M,\loc}(\calG)$ of local bisections $\sigma:U\to \calG^U_U$ such that $\varphi_\sigma = \Id_U$ is a sub-pseudo-group of $\Bis_{\loc}(\calG)$.   
The map in \eqref{eq:BisectionsToDiffeos} descends to a pseudo-group homomorphism
\begin{align*}
\varphi:\Bis_{\loc}(\calG)\to \Diff_\loc(M),\ \sigma\mapsto \varphi_\sigma = \tt\circ\sigma 
\end{align*}
from local bisections of $\calG$ to the set of diffeomorphisms between open subsets of $M$, that we improperly call \emph{local diffeomorphisms}, which sends $\Bis_{M,\loc}(\calG)$ to the collection $\Id_{\loc}(M)=\{\Id_U,\ \mbox{$U\subset M$ open}\}$ of local identity maps.

To summerize, we have the following group or pseudo-group maps between (local) bisections of $\calG$ and (local) diffeomorphisms of $M$:
$$
\begin{tikzcd}
\Bis_M(\calG) \arrow[r, hook] \arrow[d] & \Bis(\calG) \arrow[d] \\
\{\Id_M\} \arrow[r, hook] & \Diff(M)
\end{tikzcd}
\qquad\qquad
\begin{tikzcd}
\Bis_{M,\loc}(\calG) \arrow[r, hook] \arrow[d] & \Bis_{\loc}(\calG) \arrow[d] \\
\Id_{\loc}(M) \arrow[r, hook] & \Diff_\loc(M)
\end{tikzcd}
$$

\begin{example}[Bisections of gauge groupoids and gauge transformations] \label{ex:bisections}
Given a principal $G$-bundle $P\to M$, a \textbf{principal automorphism} of $P$ is a bundle isomorphism $\Phi:P\to P$ above a smooth diffeomorphism $\varphi:M\to M$ (that is, if $p\in P_x$ for some $x\in M$, then $\Phi(p)\in P_{\varphi(x)}$) such that $\Phi(pg)=\Phi(p)g$ for any $g\in G$. Principal automorphisms can be composed and their set $\Aut(P)$ forms a group. Principal automorphisms above $\Id_M$ are usually called \textbf{gauge transformations} and form a subgroup $\Aut_M(P)$.
All such automorphisms can be considered in a \emph{local} version, defined on the subbundle $P_U$ above open subsets $U\subset M$. The collections $\Aut_{\loc}(P)$ and $\Aut_{M,\loc}(P)$ of local automorphisms form two pseudo-groups with the composition. The elements of $\Aut_{M,\loc}(P)$ are usually called \textbf{local gauge transformations}. 
A principal automorphism $\Phi \in \Aut(P)$ can be identified with the bisection of the gauge groupoid $\calG(P)$ given by
\begin{align*}
\sigma_\Phi:M\to \calG(P),\ x\mapsto \sigma_\Phi(x)=[\Phi(p),p] \qquad \mbox{for any choice of $p\in P_x$}. 
\end{align*}
This gives a series of group and pseudo-group isomorphisms \cite[Examples 3.5]{Meinrenken:2017}
\begin{align*}
\Aut(P) & \cong \Bis(\calG(P)) &\qquad \Aut_{\loc}(P) & \cong \Bis_{\loc}(\calG(P)) \\
\Aut_M(P) &\cong \Bis_M(\calG(P)) &\qquad \Aut_{M,\loc}(P) &\cong \Bis_{M,\loc}(\calG(P))
\end{align*}
which recover the (local) gauge transformations of $P$ as (local) bisections of $\calG(P)$ fixing $M$.
In particular, we have $\Bis(\Pair(M)) \cong \Diff(M)$ and $\Bis_{loc}(\Pair(M))\cong \Diff_\loc(M)$, via the correspondence $\sigma(x)=(\varphi_\sigma(x),x)$.
\end{example}

A smooth (global) map $\phi: \calG\to \calG'$ between two Lie groupoids over $M$ induces a pushforward on their bisections
\begin{align}
\label{eq:MapsBisections-global}
& \Bis(\phi)=\phi_*: \Bis(\calG)\longrightarrow \Bis(\calG'),\ \sigma \mapsto \phi\circ\sigma, 
\end{align}
above the same diffeomorphism of $M$, i.e. $\varphi_{\phi\circ\sigma} = \varphi_\sigma$. 
Moreover, $\Bis(\phi)$ is a group morphism if $\phi$ is a groupoid morphism, in fact
\begin{align*}
\Bis(\phi)(\sigma'\bullet\sigma)(x) &= 
    \phi( \sigma'(\varphi_\sigma(x)) \sigma(x)) \\
& = \phi(\sigma'(\varphi_\sigma(x)))\phi(\sigma(x)) =
\big(\Bis(\phi)(\sigma')\bullet\Bis(\phi)(\sigma)\big)(x),  
\end{align*}
for any $\sigma,\sigma'$ in $\Bis(\calG)$ and for any $x$ in $M$. 
If $\phi:\calG \locmap \calG'$ is a local map, the composite map $\phi\circ \sigma$ is only defined for a local bisection $\sigma$ of $\calG$ taking values in the domain of $\phi$. Therefore, if $\phi$ is defined on a neighbourhood $\calU$ of the diagonal $\uu(M)\subset\calG$, the pushforward 
\begin{align}\label{eq:MapsBisections-local}
& \Bis_\loc(\phi)=\phi_*: \Bis_\loc(\calG) \longlocmap \Bis_\loc(\calG'),\ \sigma \mapsto \phi \circ \sigma 
\end{align}
is a local map defined on the subspace
\begin{align}\label{eq:Bisloc(calU)}
\Bis_\loc(\calU) = \{\sigma\in \Bis_\loc(\calG)\ |\ \mbox{$\sigma$ is defined on $U\subset M$ such that $\sigma(U)\subset\calU$}\}.
\end{align}
This is an open neighbourhood of the unit bisection $\uu$ in $\Bis_\loc(\calG)$, endowed with the graph topology \cite[\S 4.4]{Michor:1980} or with the compact open topology \cite[Definition A.6]{Schmeding-Wockel:2015}, because the set of diffeomorphisms $\varphi_\sigma$ is open in smooth maps \cite[Corollary 5.7]{Michor:1980} and the map \eqref{eq:BisectionsToDiffeos} sending bisections to diffeomorphisms is continuous. 


\subsection{Groupoid actions}
\label{subsec:actions}


Groupoids naturally act on fibre bundles both with fibrewise transformations (as their structure groups) and between fibres on distinct points (as a parallel transport induced by connections). 

Let $\calG\toto M$ be a Lie groupoid and $\pi:E\to M$ a fibre bundle.
A left \textbf{action} of $\calG$ on $E$ is a map \cite[Definition 1.6.1]{Mackenzie:2005}, \cite[\S 5.3. Semi-direct products]{Moerdijk-Mrcun:2003}  
\begin{align}\label{eq:actionGQ}
\rho:\calG\times_M E \longrightarrow E,\ (g,e) \longmapsto g\cdot e
\end{align} 
defined on the set of \textbf{composable} pairs in $\calG\times E$, namely 
$$
\calG\times_M E:=\{(g,e)\in \calG \times E,\ \ss(g)=\pi(e)\}
=  \bigcup_{e\in E} \calG^{\pi(e)} \times \lbrace e \rbrace, 
$$    
such that
\begin{enumerate}[label={(\roman*)}]
\item 
$\rho$ respects the fibres, i.e. $\pi(g\cdot e)=\tt(g)$ for any $e\in E$ and any $g\in \calG$, 
\item
$\rho$ is unital, i.e. $1_{\pi(e)}\cdot e=e$ for any $e\in E$, 
\item 
$\rho$ is associative, i.e. $g' \cdot(g\cdot e) = (g'\,g)\cdot e$ for any $e\in E$, any $g\in \calG$ composable with $e$ and any $g'\in \calG$ composable with $g$. 
\end{enumerate}
In other words, for any $x,y\in M$, the action map $\calG^x_y\times E_x \to E_y$ relates the fibres of $E$ above $x$ to those above $y$, in a manner analogous to the action of the structure group on each fibre.



If $E$ is a vector bundle, the action $\rho$ is called \textbf{linear} if, for any $g\in \calG^x_y$, the induced map $\rho_g:E_x\to E_y$ is linear. Since $g$ is invertible in $\calG$, the map $\rho_g$ is invertible, and therefore it belongs to $\GL(E_x,E_y)$. 
In other words, linear actions of $\calG$ on $E$ are equivalent to groupoid morphisms $\calG \to \GL(E)$ over $M$ \cite[Definition 1.7.1]{Mackenzie:2005} \cite[\S 5.3]{Meinrenken:2017}. 
In this case, the bundle $E$ is also called a \textbf{linear representation} of the groupoid $\calG$.
As usual, the representation is called \textbf{faithful} if the corresponding groupoid morphism $\calG \to \GL(E)$ is injective. 

\begin{defn} \label{def:localaction}
An action $\rho:\calG\times_M E\to E$ is \textbf{local} if it is defined on a subspace $\calU\times_M E$ of $\calG\times_M E$ where $\calU$ is an open neighbourhood of the units $\uu(M)\subset \calG$. A local action is denoted by $\rho:\calG\times_M E\locmor E$.
If $E$ is a vector bundle, a local linear action $\rho$ is equivalent to a local groupoid morphism $\calG \locmor \GL(E)$.\footnote{ 
Note that local groupoid maps $\calG\locmap \GL(E)$ are equivalent to \emph{non-associative} linear actions $\calG\times_M E\to E$, that is, maps of the form \eqref{eq:actionGQ} satisfying the axioms $(\textrm{i})$ and $(\textrm{ii})$ above, but not $(\textrm{iii})$.}
\end{defn}

\begin{examples}\label{ex:groupoidactions}
Let $E\to M$ be a vector bundle. 
\begin{enumerate}[label={(\alph*)}]
\item 
The frame groupoid $\GL(E)\toto M$ has a natural faithful linear representation on $E$ given by the \textbf{evaluation} 
\begin{eqnarray}\label{eq:actionIso(E)}
\mathrm{ev}:\GL(E) \times_M E \to E,\ (\varphi^x_y,e_x) \mapsto \varphi^x_y(e_x)=:\varphi^x_y\cdot e_x. 
\end{eqnarray}

\item 
If the structure group of $E$ reduces to a Lie group $G$, let $P\to M$ be the principal $G$-bundle obtained by reduction of the frame bundle $FE$.
Then the gauge groupoid $\calG(P)$ acts linearly on $E$, with action $\rho:\calG(P)\times_M E\to E$ given by composing the map $\calG(P)\to \GL(E)$ of Example \ref{ex:framegroupoid} with the evaluation $\mathrm{ev}:\GL(E) \times_M E \to E$. 

\item 
A linear representation of the pair groupoid $\Pair(M)$ on $E$ is equivalent to a global trivialisation of $E$, because all fibres are identified by the linear isomorphisms $\rho^x_y$ of the associated groupoid morphism $\rho:\Pair(M)\to \GL(E)$, cf. \cite[\S 5.3]{Meinrenken:2017}.
If $M$ is connected and simply connected, a local linear representation $\rho:\Pair(M)\locmor \GL(E)$ is enough to trivialize $E$, either by extending it to a global one as at the end of Section \ref{subsec:local maps}, or by identifying far away fibres by composing isomorphisms among nearby fibres covered by the diagonal domain of $\rho$.

\end{enumerate}
\end{examples}


\subsection{Lie algebroids}
\label{subsec:algebroids}
 

Lie algebroids are fibred analogues of Lie algebras over a manifold, and play for Lie groupoids the role that Lie algebras play for Lie groups.  
They were introduced by J.~Pradines in \cite{Pradines:1967} and are nowadays used, for instance, to study foliations \cite{Moerdijk-Mrcun:2003}, Poisson geometry \cite{Marle:2006} and sigma models in string theory \cite{Bojowald-Kotov-Strobl:2005}. 
We are mainly concerned with the Lie algebroid defined by a Lie groupoid, but there are many other interesting examples of Lie algebroids which go beyond our scope, cf. \cite{Mackenzie:2005}, \cite{Meinrenken:2017} and \cite{Moerdijk-Mrcun:2003}. 

A \textbf{Lie algebroid} on a manifold $M$ is a vector bundle $q: A\to M$ endowed with
\begin{enumerate}[label={(\roman*)}]
\item
a vector bundle map $a:A\to TM$ over $M$ called the \textbf{anchor}, which can be seen as a $C^\infty(M)$-map from the sections of $A$ to the vector fields on $M$, 
\item a Lie bracket $[\ ,\ ]_A\colon \Gamma(A)\times \Gamma(A)\to \Gamma(A)$ on its space of sections (that is, an antisymmetric $\R$-bilinear map satisfying the Jacobi identity), which satisfies the Leibniz rule
\begin{eqnarray*}
[\xi,f\,\eta]_A = f\,[\xi,\eta]_A+ a(\xi)(f)\,\eta 
\end{eqnarray*} 
for any $\xi,\eta \in \Gamma(A)$ and $f\in C^\infty(M)$, that is, each   $[\xi,\cdot\, ]_A$ is a derivation.
\end{enumerate}
It then follows that the anchor is a morphism of Lie algebras from $\Gamma(A)$ to the space $\X(M)$ of vector fields on $M$, i.e. $a\left([\xi,\eta]_A\right)= \left[a(\xi),a(\eta)\right]$, cf. \cite[Remark 6.2 (b)]{Meinrenken:2017}. 
The Lie algebroid $A$ is \textbf{regular} if $a$ has constant rank, it is \textbf{transitive} if $a$ is fibrewise surjective and it is \textbf{totally intransitive} if $a=0$.
A Lie algebroid on $M$ is often denoted by its anchor map $A\to TM$, instead of its bundle map to $M$.  

A \textbf{morphism} between two Lie algebroids $A\to TM$ and $A'\to TM$ on the same base manifold is a morphism of vector bundles $\phi: A\to A'$ which commutes with the anchors and preserves the brackets.

\begin{example}[Atiyah algebroid of a principal bundle] \label{ex:Atiyah algebroid}
The Atiyah algebroid of a principal $G$-bundle $\pi:P\to M$ plays the role of the Lie algebra of the Lie group $G$: it contains the right invariant vector fields on $P$, which are then determined by their value on any chosen element $p_x$ of each fibre (which represents the unit in $G$ above the point $x\in M$). 

The key fact to implement this analogy is that the group $G$ acts on the tangent bundle $TP\to P$, with action obtained by differentiating the principal action on $P$. The space of sections of the quotient bundle
\begin{align}\label{eq:Atiyahalg}
A(P):=TP/G\to M
\end{align}
is isomorphic to that of $G$-invariant vector fields on $P$ \cite[\S 3.5]{Mackenzie:2005}, and therefore it is closed under the Lie bracket of vector fields. Hence $A(P)$ is a Lie algebroid, called the \textbf{Atiyah algebroid} of $P$, with anchor $a:A(P)\to TM$ induced by $d\pi:TP\to TM$. The anchor is fibrewise surjective, hence $A(P)$ is transitive. 
The quotient map $\natural:TP\to TP/G=A(P)$ is a fibrewise isomorphism between the two bundles $TP\to P$ and $A(P)\to M$ over different base manifolds, and therefore it is fibrewise invertible with multi-valued inverse $\natural^{-1}:A(P)\to \bigcup_{p\in P} T_pP$.
The Atiyah algebroid $A(P)$ is also called the \textbf{bundle of connections} of $P$ in \cite{Kobayashi:1957}, because sections of $d\pi:TP\to TM$ determine horizontal distributions in $TP$, and the $G$-invariant ones, which are equivalent to principal connections on $P$, are exactely the image under $\natural^{-1}$ of the sections of the anchor $a:A(P)\to TM$. 

Being transitive, the Atiyah algebroid fits into a short exact sequence of Lie algebroids over $M$, called the \textbf{Atiyah sequence},
\begin{align}\label{eq:Atiyahseq}
0\longrightarrow \mathrm{Ker}(a)\cong (P\times \mathfrak{g})/G =:\Ad(P) \longrightarrow A(P)= TP/G \overset{a}{\longrightarrow} TM \longrightarrow 0, 
\end{align}
obtained by taking the quotient under the right $G$-action of the analogue sequence 
$$
0\longrightarrow VP = \ker(d\pi)\cong P\times \mathfrak{g} \longrightarrow TP \overset{d\pi}{\longrightarrow} TM \longrightarrow 0, 
$$
where the action of $G$ on $P\times \mathfrak{g}$ is obtained by differentiating the action of $G$ on $P\times G$ given by $(p,g)\cdot h=(ph,h^{-1}gh)$ for $p$ in $P$ and $g,h$ in $G$, cf. \cite[\S 3.2]{Mackenzie:2005}, \cite[\S 6.4]{Moerdijk-Mrcun:2003}.
The kernel $\Ad(P)$ of the anchor is the so-called \textbf{adjoint bundle} of $P$ and it is a totally intransitive Lie algebroid. It is in fact a \textbf{Lie algebra bundle} on M, that is, a bundle with Lie algebras as fibres (isomorphic to $\mathfrak{g}$), which is fibrewise isomorphic to the vertical subbundle $VP \subset TP$.
\end{example}


Unlike for principal bundles, a Lie groupoid $\calG\toto M$ does not globally act on itself, nor on its tangent bundle $T\calG\to \calG$, because the right multiplication by a given element $g\in \calG$ is only defined on the subspace $\calG^{\tt(g)}$ with source depending on the target of $g$, and defines a diffeomorphism $\calG^{\tt(g)}\to \calG^{\ss(g)}$ between $\ss$-fibres of $\calG$. 
To obtain a subspace of $\calG$-invariant vector fields, determined by their value at the units $1_x\in \calG$ and giving rise to a true bundle over $M$, the suitable analogue of $TP$ is the \textbf{$\ss$-tangent} subbundle of $T\calG$
$$
T^s\calG := \ker(d\ss:T\calG\to TM) = \bigcup_{g\in \calG} T_g \calG^{\ss(g)} \longrightarrow \calG, 
$$
made of vectors tangent to the $\ss$-fibres, which inherits the differential of the target $d\tt|_{T^s\calG}:T^s\calG\to TM$ serving as anchor. 
Indeed, $\calG$ acts on $T^s\calG$ on the right and the space of sections of the quotient bundle $T^s\calG/\calG$ on $M$ is isomorphic to that of right-invariant vector fields on $\calG$ (with values necessarily in $T^s\calG$) and therefore it is closed for the Lie bracket \cite[\S 3.5]{Mackenzie:2005} \cite[\S 4.1]{Moerdijk-Mrcun:2003}\footnote{An equivalent definition of the Lie algebroid $\Lie(\calG)$ is given by \cite[\S 9.2]{Meinrenken:2017} or again \cite[\S 4.1]{Moerdijk-Mrcun:2003} as the \emph{normal bundle} $T\calG_{|u(M)}/Tu(M)$ of $M$ in $\calG$ (cf. \cite[\S 8.2]{Meinrenken:2017}).}.
Finally, the \textbf{Lie algebroid} of $\calG$ is the quotient bundle $T^s\calG/\calG$ on $M$, which can also be described as the pull-back of $T^s\calG\to \calG$ along the embedding $\uu: M\to \calG$, that is, 
\begin{align*} 
\Lie(\calG):=T^s\calG/\calG \cong \uu^*(T^s\calG) = \bigcup_{x\in M} T_{1_x}\calG^x. 
\end{align*}
The bundle projection onto $M$ is given by the source, with constant value $x$ on the fibre $T_{1_x}\calG^x$. 
The anchor is the composite map $a=d\tt|_{T^s\calG}\circ i : \Lie(\calG)\to TM$, where $i:\uu^*(T^s\calG)\to T^s\calG\subset T\calG$ is the natural map on the pull-back. The Lie bracket on $\Gamma(A(\calG))$ is that of right-invariant vector fields on $\calG$.

From the definition it follows that any vector $Y\in A_x(\calG)$ determines a $\tt$-horizontal lift of $a(Y)=X\in T_xM$ to $T^s\calG$ above any $g\in \calG^x$, namely the vector
\begin{align*}
\natural^{-1}_g(Y)\in T_g\calG^x\subset T^s\calG,  
\end{align*}
where $\natural^{-1}:\Lie(\calG)\to T^s\calG$ is the multi-valued inverse of the quotient map of bundles
$$
\begin{tikzcd}
T^s\calG \arrow[r, "\natural" above] \arrow[d] & T^s\calG/\calG = \Lie(\calG) \arrow[d] \\
\calG \arrow[r, "\tt" above] & M
\end{tikzcd}, 
$$
which is a fibrewise isomorphism.

\begin{remark} \label{rem:transitive Lie algebroid}
As recalled in Example \ref{ex:gaugegroupoid}, transitive Lie groupoids $\calG$ on a connected base manifold $M$ are equivalent to gauge groupoids $\calG(P)$.
The analogue property for Lie algebroids does not hold: there are transitive Lie algebroids $A$ which are \emph{not} Atiyah algebroids $A(P)$ of a principal bundle, cf. \cite[Example 6.9]{Meinrenken:2017}. This is related to the \emph{integrability} of a Lie algebroid $A$ to a Lie groupoid $\calG$ such that $A=A(\calG)$: the cohomological obstruction to integrability is described by C. Mackenzie in \cite[Ch. 8]{Mackenzie:2005} for transitive Lie algebroids and by M. Crainic and R. Fernandes in \cite{Crainic-Fernandes:2003} for general ones, see also \cite[\S 10.2]{Meinrenken:2017}. 
\end{remark}

If a Lie algebroid $A(\calG)$ is \emph{assumed} to come from a Lie groupoid, then the transitivity of $A(\calG)$ (on a connected base manifold) implies the transitivity of $\calG$ \cite[Corollary 3.5.18]{Mackenzie:2005}. Then $\calG=\calG(P)$ for a principal bundle $P$ and one can compare $A(\calG(P))$ to $A(P)$. 

\begin{examples}[Lie algebroids of gauge groupoids] \label{exs:Lie algebroids}
The Lie algebroid of a gauge groupoid $\calG(P)\toto M$ is isomorphic to the Atiyah algebroid of the principal bundle $P$ \cite[Example 9.5 (c)]{Meinrenken:2017}, i.e. 
\begin{align}\label{eq:L(G(P))=A(P)}
\Lie(\calG(P))\cong A(P)\to TM.
\end{align}
Therefore gauge groupoids give rise to various classical examples of Lie algebroids, which all play a role in the theory of connections. 
\begin{enumerate}[label={(\alph*)}]
\item 
The Lie algebroid of the pair groupoid $\Pair(M)\toto M$
is isomorphic to the tangent bundle $TM\to M$ equipped with the vector field brackets and the identity $TM\to TM$ as anchor, i.e. $\Lie(\Pair(M))\cong TM$ \cite[Example 3.5.11]{Mackenzie:2005}. 
\item If $P=M\times G$ is the trivial principal $G$-bundle, then $\calG(P)=M\times G\times M$ is the trivial groupoid and $\Lie(M\times G\times M) \cong A(M\times G) \cong TM \oplus (M\times \mathfrak{g})$ is called the \textbf{trivial Lie algebroid} with structure Lie algebra $\mathfrak{g}=\mathrm{Lie}(G)$. The anchor is given by the projection to $TM$, and the bracket is
\begin{align*}
[X+v,Y+w] = [X,Y]_{\mathfrak{X}(M)} + \big(X(w)-Y(v)+[v,w]_{\mathfrak{g}}\big),
\end{align*}
where $X,Y\in \mathfrak{X}(M)$, $v,w\in C^\infty(M,\mathfrak{g})$ and $X(w)$ denotes the Lie derivative of $w$ with respect to $X$ \cite[Example 3.5.13]{Mackenzie:2005}. 
When $G=\{1\}$, this gives back $\Lie(\Pair(M))\cong TM$.

\item Let $G$ and $H$ be two (unrelated) Lie groups, and let $Q$ be a principal $H$-bundle over $M$. Then $P=Q\times G$ is a principal $(H\times G)$-bundle over $M$, and its gauge groupoid is
\begin{align}\label{eq:calG(QxG)}
\calG(P) &= (Q\times G)\times_{H\times G} (Q\times G) \cong \calG(Q)\times G, 
\end{align}
with structure maps inherited by $\calG(Q)$ and multiplication extended to $G$.
Its Lie algebroid is
\begin{align}\label{eq:A(QxG)}
\Lie(P) &\cong \Lie(Q)\oplus (M\times G), 
\end{align}
and \emph{not} $\Lie(Q)\oplus (Q\times G)$, because on the tangent bundle $TP\cong TQ \oplus (Q\times \mathfrak{g})$ the group $H\times G$ acts on both factors of $Q\times \mathfrak{g}$ (and on $TQ$ it acts by $H$). The anchor is given by that of $\Lie(Q)$ and the bracket is 
\begin{align*}
[\xi+v,\eta+w]_{\Lie(P)} = [\xi,\eta]_{\Lie(Q)} + \big(a(\xi)(w)-a(\eta)(v)+[v,w]_{\mathfrak{g}}\big),
\end{align*}
where $\xi,\eta\in \Gamma(\Lie(Q))$, $v,w\in C^\infty(M,\mathfrak{g})$ and $a(\xi)\in \mathfrak{X}(M)$ is the anchor of $A(Q)$. 
Note that the previous example of the trivial bundle $P=M\times G$ corresponds to $H=\{1\}$, not $G=\{1\}$. 


\item If $\pi:E\to M$ is a vector bundle, the Lie algebroid of its frame bundle $FE$ and of the frame groupoid $\GL(E)$ is the \textbf{bundle of linear vector fields of $E$}, $\Tlin E\to M$, with fibre above $x$ given by the linear sections of the canonical projection $TE|_{E_x}\to E_x$, that is, smooth maps $\xi_x:E_x\to TE|_{E_x}$ of the form $\xi_x(e)=\frac{d}{dt}(e+t \xi_x)\big|_{t=0}$ for $e\in E_x$. The projection of $\xi_x$ to $M$ through $\pi$ is then a well defined tangent vector $X_{\xi_x}=d\pi_e(\xi_x(e)) = \frac{d}{dt}(\pi(e+t \xi_x))\big|_{t=0}$, independent of the choice of $e\in E_x$. The assignement $a(\xi_x)=X_{\xi_x}\in T_xM$ defines the anchor $a:\Tlin E\to TM$.

In other words, the fibre of $\Tlin E$ at $x\in M$ is made of pairs of maps $(\xi_x,X_{\xi_x})$ which make the following diagram commute
$$
\begin{tikzcd}
E_x \arrow[r, "\xi_x" above] \arrow[d, "\pi" left] & TE|_{E_x} \arrow[d, "d\pi" right] \\
\{x\} \arrow[r, "X_{\xi_x}" above] & T_xM
\end{tikzcd}.
$$
The sections of the bundle $\Tlin E\to M$ are then sections of the canonical double bundle
$$
\begin{tikzcd}
E \arrow[d, "\pi" left] & TE \arrow[l] \arrow[d, "d\pi" right] \\
M & TM \arrow[l] 
\end{tikzcd}
$$
that is, pairs $(\xi,X_\xi)$ where $\xi$ is a vector field on $E$ and also a bundle map $\xi:E\to TE$ above the vector field $X_\xi:M\to TM$ given by $(X_\xi)_x=X_{\xi_x}$ for $\xi_x=\xi|_{E_x}$. These pairs are closed for the Lie bracket of vector fields on $E$ and on $M$, so $\Tlin E$ is a Lie algebroid. Details can be found in \cite[\S 3.4]{Mackenzie:2005} or in \cite[Example 2.3]{Laurent-Gengoux-Voglaire:2020}. 

Since $d\pi\circ \xi = X_\xi\circ\pi$, the field $\xi$ is a horizontal lift of $X_\xi$ and the anchor $a:\Tlin E\to TM$ is just the projection of a linear vector field $\xi$ on $E$ to its base vector field $X_\xi$ on $M$. 
Then, the kernel of the anchor is the bundle $V^\mathrm{lin}E\to M$ of \textbf{linear vertical vector fields of $E$}, that is, the linear vector fields $\xi:E\to TE$ with values in the vertical subbundle $VE=\ker(d\pi)$ of $TE$. Since the vertical fibres $V_eE$ are canonically isomorphic to $E_{\pi(e)}$ (by translation of the $0$ vector), there is a canonical bundle isomorphims $V^\mathrm{lin}E \cong \End(E)$. 
This bundle fits in the Atiyah exact sequence of algebroids 
\begin{align}\label{eq:exactTlinEsequence}
0\longrightarrow V^\mathrm{lin}E\cong \End(E) \longrightarrow \Tlin E \overset{a}{\longrightarrow} TM\longrightarrow 0.
\end{align} 
\end{enumerate}
\end{examples}

\begin{example}[The algebroid of derivative endomorphisms] \label{ex:Der(E)} 
The Lie algebroid $A(\GL(E))=A(FE)=\Tlin E\to TM$ is suitable to describe \emph{horizontal lifts} in $TE$ of vectors tangent to $M$, i.e. Ehresmann connections on the vector bundle $E$. Linear connections on $E$, that is, $\R$-linear smooth maps $\nabla:\mathfrak{X}(M)\times \Gamma(E) \to \Gamma(E), (X,\psi)\to \nabla_X\psi$ which are $C^\infty(M)$-linear in $X$ and $C^\infty(M)$-derivations in $\psi$, are \emph{equivalent} but \emph{not equal} to Ehresmann connections, since they describe the \emph{vertical} variation of sections once a horizontal distribution in $TE$ is fixed.  
Linear connections on $E$ are also related to a Lie algebroid, which is \emph{not} the Lie algebroid of $FE$ or $\GL(E)$, although it is isomorphic to $\Tlin E$. 

A \textbf{linear derivation} on $E$ at $x\in M$ is an $\R$-linear map $D_x:\Gamma(E)\to E_x$ satisfying the Leibniz rule 
\begin{align}\label{eq:Derx}
D_x(f\,\psi) = X_{D_x}(f)\, \psi(x)+f(x)\,D_x\psi \qquad \forall f\in C^\infty(M),\ \forall \psi\in \Gamma(E), 
\end{align}
for a given tangent vector $X_{D_x}\in T_xM$ associated to $D_x$, which acts as a derivation of functions on $M$ at $x$. The \textbf{bundle of derivative endomorphisms} of $E$ \cite[Example 3.3.4]{Mackenzie:2005}, also called the \textbf{bundle of linear derivations} on $E$ \cite[\S 2.1]{Esposito-Vitigliano-Tortorella:2019}, is the collection $\Der(E)$ of the vector spaces $\Der_x(E)$ of derivations of sections of $E$ at each point of $M$. 
Its sections are the first order differential operators $D:\Gamma(E)\to \Gamma(E)$ satisfying the Leibniz rule \eqref{eq:Derx} for an associated vector field $X_D\in \mathfrak{X}(M)$ defined by $(X_D)_x = X_{D_x}$, with $D_x(\psi)=D(\psi)(x)$. 
The vector bundle $\Der(E)$ inherits a Lie bracket from the commutator of linear operators, and the anchor $a:\Der(E)\to TM$ is given by the symbol map $a(D_x)=X_{D_x}$. Therefore $\Der(E)\to TM$ is a Lie algebroid over $M$, and it is clearly transitive.

The kernel of the anchor at $x\in M$ is given by the $\R$-linear maps $D_x: \Gamma(E)\to E_x$ which are also $C^\infty(M)$-linear. Hence $\ker(a)$ is isomorphic to the Lie algebra bundle $\End(E)$ and gives rise to a short exact sequence similar to the Atiyah sequence \eqref{eq:exactTlinEsequence}, 
\begin{align}\label{eq:exactDerEsequence}
  0\longrightarrow  \End(E)\longrightarrow \Der(E)\overset{a}{\longrightarrow} TM\longrightarrow 0. 
\end{align}

Note that, for a given vector field $X\in \mathfrak{X}(M)$, the covariant derivative $\nabla_X$ of a linear connection on $E$ is precisely a section of $\Der(E)$, with $D=\nabla_X$ and $a(D)=X$. The relationship between the two algebroids $\Tlin E$ and $\Der(E)$ mimics the relationship between a horizontal lift to $E$ of $X\in \mathfrak{X}(M)$ and a covariant derivative $\nabla_X$, which is a vertical-valued operator in $TE$: the isomorphism 
\begin{align}\label{eq:TlinE iso Der(E)} 
\Tlin E\overset{\simeq}{\longrightarrow} \Der(E),\ (\xi_x,X_{\xi_x})\longmapsto D_{\xi_x}
\end{align}
is given by
\begin{align}\label{eq:def TlinE iso Der(E)}
D_{\xi_x}(\psi) = d\psi|_x(X_{\xi_x}) - \xi_x(\psi(x)) \qquad \forall \psi\in \Gamma(E), 
\end{align}
where $d\psi|_x:T_xM\to T_{\psi(x)}E$ is the total differential of $\psi$ at $x$ (its tangent map), with symbol $a(D_{\xi_x})=X_{\xi_x}$. The operator $D_{\xi_x}$ is indeed a linear derivation on $E$ at $x$, because $d\psi|_x$ is subject to the Leibniz rule and therefore, for any $f\in C^\infty(M)$, we have
\begin{align*}
  D_{\xi_x}(f\psi) & = d(f\psi)\big|_x(X_{\xi_x}) - \xi_x(f(x)\psi(x)) \\ 
  &= X_{\xi_x}(f) \psi(x) + f(x) d\psi|_x(X_{\xi_x}) - f(x) \xi_x(\psi(x)) \\ 
  &= X_{\xi_x}(f) \psi(x) + f(x) D_{\xi_x}(\psi).  
\end{align*}
See \cite[Proposition 3.4.4]{Mackenzie:2005} or \cite[Example 2.3]{Laurent-Gengoux-Voglaire:2020} for an alternative definition of this isomorphism using the canonical non-degenerate pairing between $E$ and its dual vector bundle $E^*$, and \cite[Theorem 3.6.6]{Mackenzie:2005} for an alternative isomorphism based on an exponential map from $A(\GL(E))$ to $\GL(E)$. 
\end{example}

\begin{remark}\label{rk:A'(G)}
The Lie algebroid $\Der(E)\to TM$ can be related to the Lie groupoid $\GL(E)\toto M$ through an alternative construction of a Lie algebroid $A'(\calG)\to TM$ from any Lie groupoid $\calG\toto M$, obtained by swapping the source and the target everywhere along the construction of $A(\calG)$. 
It is immediate to show that $A(\calG)$ and $A'(\calG)$ share the same Lie algebra bundle $L(\calG)$ as kernel of the respective anchors, i.e. there are two exact sequences 
\begin{align*}
&  0\longrightarrow L(\calG)\cong \mathrm{Ker}(a) \longrightarrow A(\calG) \overset{a}{\longrightarrow} TM \longrightarrow 0, \\
&  0\longrightarrow L(\calG)\cong \mathrm{Ker}(a') \longrightarrow A'(\calG) \overset{a'}{\longrightarrow} TM \longrightarrow 0. 
\end{align*}

When $\calG=\GL(E)$, swapping source and target essentially corresponds to swapping \emph{horizontal} and \emph{vertical} vectors tangent to $E$, at the instar of \eqref{eq:def TlinE iso Der(E)}, which can be written as 
\begin{align}\label{eq:d=xi+D} 
d\psi|_x(X_x) = \xi_x(\psi) + D_x\psi|_{x}, \qquad a(\xi_x)=a'(D_x)=X_x,
\end{align}
and presents the decomposition of the vector $d\psi|_x(X_{\xi_x})\in T_{\psi(x)}E$ into a horizontal part $\xi_x(\psi)\in T_{\psi(x)}E$ (obtained through $A(\GL(E))=\Tlin E)$) and a vertical part $D_x(\psi)=\nabla_{X_x}\psi|_{x}\in E_x\cong V_{\psi(x)}E$ (obtained through $A'(\GL(E))=\Der(E)$). 
\end{remark}


To conclude this section, we complete the \textbf{Lie functor} from groupoids to algebroids by recalling the construction of the algebroid map induced by a groupoid map. 
Just as the tangent map of a morphism of Lie groups, at the unit, induces a Lie algebra morphism between the corresponding Lie algebras, there exists a tangent map of a local morphism between Lie groupoids, at the units, which gives rise to a Lie algebroid morphism. 
More generally, one can differentiate at the units any local map between Lie groupoids, and obtain its infinitesimal part along the diagonal.

Let $\phi:\calG \locmap \calG'$ be a local map between two Lie groupoids on $M$, defined on an open neighbourhood $\calU\subset \calG$ of the units. 
For any $x\in M$, we denote by $d^\tt\phi\vert_x$ the \textbf{target} or \textbf{$\tt$-differential} of $\phi$ at $1_x\in \calU$ given by
\begin{align*}
d^\tt\phi\vert_x(\dot{\gamma}(0)) &= \frac{d\ }{dt} \phi\big(\gamma(t)\big)\vert_{t=0}, 
\end{align*}
where $\gamma:(-\varepsilon,\varepsilon)\to \calU\cap \calG^x$ is any smooth curve living in the source-fibre of $x$ such that $\gamma(0)=1_x$. 
\begin{prop} \label{prop:t-differential} 
The $\tt$-differential at the units of a local map $\phi: \calG \locmap \calG'$ between two Lie groupoids induces a vector bundle morphism
\begin{align*}
d^{\tt}\phi\vert_M: \Lie(\calG) \longrightarrow \Lie(\calG')
\end{align*}
between their associated Lie algebroids. 
If $\phi$ is a local groupoid morphism, then $d^\tt\phi\vert_M$ is a Lie algebroid morphism. 
\end{prop} 

\begin{proof}
By assumption, the $\tt$-differential $d^\tt\phi$ is defined on vector fields tangent to the source fibres of $\calG$ at the units, which span the fibres of the vector bundle $\uu^*(T^s\calG)=\Lie(\calG)$. Since $\phi$ preserves the source and the units, by Definition \ref{def:localmap}, the map $d^\tt\phi$ maps $\uu^*(T^s\calG)$ to $(\uu')^*(T^s\calG')$ and gives a vector bundle morphism $d^\tt\phi\vert_M: \Lie(\calG)\to \Lie(\calG')$.
The second assertion is proved in \cite[\S 3.5, (40) and (41)]{Mackenzie:2005}.    
\end{proof}

\begin{remark} \label{rem:s-differential} 
Using the alternative Lie algebroid $\Lie'(\calG)$ described in Remark \ref{rk:A'(G)}, one can define the \textbf{source} or \textbf{$\ss$-differential} of a local groupoid map $\phi:\calG\locmap \calG'$ at the units, giving rise to vector bundle morphism
\begin{align*}
d^{\ss}\phi\vert_M: \Lie'(\calG) \longrightarrow \Lie'(\calG'),
\end{align*}
by setting 
\begin{align*}
d^\ss\phi\vert_x(\dot{\gamma}(0)) &= \frac{d\ }{dt} \phi\big(\gamma(t)\big)\vert_{t=0}, 
\end{align*}
where $\gamma:(-\varepsilon,\varepsilon)\to \calU\cap \calG_x$ is any smooth curve living in the target-fibre of $x$ such that $\gamma(0)=1_x$. 
\end{remark}

If $G$ and $G'$ are two Lie groups, a Lie algebra morphism $\mathfrak{g}\to \mathfrak{g}'$ between their Lie algebras always integrates to a Lie group morphism $G\to G'$ if $G$ is connected and simply connected. 
The same does \emph{not} hold true for morphisms of Lie algebroids, the obstruction to integrability is studied in \cite{Crainic-Fernandes:2003}, see also \cite[Ch. 8]{Mackenzie:2005}.

\section{Direct connections on Lie groupoids}
\label{sec:direct-connections} 

In this section we consider direct connections on Lie groupoids as generalisations of parallel transport along small geodesics seen as a local map from the pair groupoid of a base manifold to the frame groupoid of a vector bundle, in the spirit of N. Teleman and J. Kubarski \cite{Teleman:2004,Teleman:2007,Kubarski-Teleman:2006,Kubarski:2008}. 
We show that a Lie groupoid which admits a direct connection is a gauge groupoid (Proposition \ref{prop:direct}), that direct connections generalise principal connections on principal bundles (Proposition \ref{prop:tau=dir conn} and Example \ref{ex:dircon on RxR}), that they induce infinitesimal connections on the associated Lie algebroid (Proposition \ref{prop:dircttoinf}), and consequently that they induce a linear connection on associated vector bundles (Corollary \ref{cor:deltaIsoE}). Finally, we characterise direct connections that are exactly the parallel transport of their linear connection (Theorem \ref{thm:transitive = parallel transport}). 
Then we introduce the curvature of a direct connection, which measures how much the direct connection deviates from being a groupoid morphism, and show that its antisymmetric component gives rise to the usual curvature tensor of the induced infinitesimal connection (Theorem \ref{thm:curvature}), whilst producing finer invariants than just the infinitesimal ones. 
We conclude with a strong characterisation of groupoids admitting a flat direct connection (Proposition \ref{prop:flat groupoid}), largely studied by K. Mackenzie as a groupoid morphism \cite{Mackenzie:2005}. 

\paragraph{Assumption:} 
We suppose that the base manifold $M$ is equipped with a linear connection $\nabla^M$ on $TM$, so that one can consider geodesics on $M$.  
We further assume that $\nabla^M$ is {\it torsion-free}, so that around each point $x\in M$ there exists a convex neighbourhood $U$, that is, any two points in $U$ are linked by a unique geodesic contained in $U$ \cite{Whitehead:1932}, \cite[Theorem 8.7]{Kobayashi-Nomizu:1963}.


\subsection{Direct connections and generalised parallel transports}
\label{sec:direct connections}


Recall from Example \ref{ex:diagonal domain} that a \textbf{diagonal domain} in $\Pair(M)$ is any open neighbourhood $\calU$ of the diagonal $\Diag(M)=\{(x,x),\, x\in M\}$.
We call $\calU$ a \textbf{convex diagonal domain} if, for any $(y,x)\in \calU$, there exists an open convex subset in $M$ containing both $x$ and $y$. 
Without further mention, from now on we assume that $\calU$ is also \textbf{symmetric}, that is, if $(y,x)\in \calU$ then also $(x,y)\in \calU$.

\begin{defn}
\label{def:directconnection}
A \textbf{direct connection} on a Lie groupoid $\calG\toto M$ is a smooth 
local right inverse of the anchor which preserves the units, that is, a smooth local map $\dircon: \Pair(M)\longlocmap \calG$ over $M$, defined on a symmetric diagonal domain $\calU\subset \Pair(M)$, which preserves the source, the target and the units, i.e. 
\begin{enumerate}[label={(\roman*)}] 
\item $\dircon(y,x)\in \calG_y^x$ for all $(y,x)\in \calU$, 
\item $\dircon(x,x)=1_x\in \calG_x^x$ for all $x\in M$. 
\end{enumerate}
We call the direct connection \textbf{antisymmetric} if it also preserves the inversion, {\it i.e.}
\begin{enumerate}[resume,label={(\roman*)}]
\item $\dircon(y,x)^{-1}=\dircon(x,y)$ for all $(y,x)\in \calU$, 
\end{enumerate}
\textbf{transitive} if it is a groupoid morphism along geodesics, i.e. 
\begin{enumerate}[resume,label={(\roman*)}]
\item $\dircon(z,y)\dircon(y,x)=\dircon(z,x)$ for all $(y,x)\in \calU$ and $z$ belonging to the geodesic from $x$ to $y$, or viceversa: $(z,x)\in \calU$ and $y$ belonging to the geodesic from $x$ to $z$, 
\end{enumerate}
and \textbf{flat} if it is a local groupoid morphism, i.e. 
\begin{enumerate}[resume,label={(\roman*)}]
\item $\dircon(z,y)\dircon(y,x)=\dircon(z,x)$ for all $x,y,z\in M$ such that $(y,x),(z,y),(z,x)\in \calU$. 
\end{enumerate}
Note that transitiveness implies antisymmetry, by taking $z=x$ the initial point of the geodesic from $x$ to $y$. 
Therefore we have
$$
\text{flatness}\quad\Longrightarrow\quad 
\text{transitiveness}\quad\Longrightarrow\quad 
\text{antisymmetry}. 
$$
Finally, a direct connection is called \textbf{global} if it is defined on the whole pair groupoid $\Pair(M)$.
\end{defn}

\begin{remark}[Direct connections and bisections] \label{rem:dircon vs bisection}
Let $\dircon:\Pair(M)\locmap \calG$ be a direct connection defined on the diagonal domain $\calU$, and denote
\begin{align*}
    \Diff_\calU(M) &= \{ \varphi\in \Diff(M)\ |\ (\varphi(x),x)\in \calU,\ \forall x\in M\}. 
\end{align*}
Then $\dircon$ induces a map $\widehat{\dircon}:\Diff_\calU(M)\to \Bis(\calG)$ which brings a diffeomorphism $\varphi$ into the global bisection $\sigma_\varphi$ of $\calG$ defined by 
\begin{align*}
    \sigma_\varphi(x) &= \dircon(\varphi(x),x),\qquad \forall x\in M.
\end{align*}
The map $\widehat{\dircon}$ is injective, it sends the identity $\id$ on $M$ to the unit map $\uu:M\to \calG$, seen as a global bisection, and in fact it is a morphism of groups. 
\end{remark}

The only direct connection on the pair groupoid is the identity map $\Id:\Pair(M)\to \Pair(M)$, and it is clearly flat. Its associated map $\widehat{\Id}$ gives the well known isomorphism of groups $\Diff(M)\cong \Bis(\Pair(M))$.

\begin{example}[Trivial direct connection] \label{ex:trivial dircon}
Any smooth map $\tau: M\to G$ with values in a Lie group $G$ gives rise to a global direct connection $\dircon_\tau$ on the trivial gauge groupoid $M\times G\times M\toto M$, by setting
\begin{align*}
\dircon_\tau(y,x) &= (y,\tau(y) \tau(x)^{-1},x)
\end{align*}
for any $(y,x)\in \Pair(M)$. This direct connection is flat, we call it \textbf{trivial}. 
\end{example}

\begin{defn} \label{def:locally trivial dircon}
Let $\calG$ be a locally trivial groupoid with section atlas $\tau_\alpha:U_\alpha\to \calG^{x_0}_{U_\alpha}$ with respect to a fixed point $x_0$ (cf. Example \ref{ex:gaugegroupoid}). 
A direct connection $\dircon$ on $\calG$, defined on a diagonal domain $\calU$, is called \textbf{locally trivial} if its restrictions to the charts of the section atlas are trivial (cf. Example \ref{ex:trivial dircon}), i.e. 
\begin{align}\label{eq:exactGamma}
\dircon\vert_{(U_\alpha\times U_\alpha)\cap \calU}(y,x)= \tau_\alpha(y)\,\tau_\alpha(x)^{-1}.
\end{align}
Such a direct connection is automatically flat. 
\end{defn}

Before giving more examples of direct connections, let us prove a necessary condition for Lie groupoids to admit direct connections. 

\begin{prop} \label{prop:direct}
If the base manifold $M$ is connected, a Lie groupoid $\calG\toto M$ equipped with a direct connection is a gauge groupoid.
\end{prop}
 
\begin{proof}
Since $M$ is connected, by Example \ref{ex:gaugegroupoid} it is enough to show that $\calG$ is transitive.
Let $\dircon:\Pair(M)\longlocmap \calG$ be a direct connection defined on a diagonal domain $\calU$.
Given $x\in M$, the pair $(x,x)$ lies in the open set $\calU$, so there is an open neighbourhood $W_x$ of $x$ in $M$ with $W_x\times W_x\subset \calU$. For every $y\in W_x$ the arrow $\dircon(y,x)\in \calG^x_y$ is then defined, so $y$ lies in the orbit of $x$. Hence $W_x\subset \mathcal{O}_x$, and every orbit is open.
The orbits partition $M$, so the complement of each of them is a union of orbits and is open as well: every orbit is therefore a non-empty open and closed subset of $M$. Since $M$ is connected, there is a single orbit, that is, $\calG$ is transitive.
\end{proof}

A direct connection also yields an explicit section atlas.

\begin{remark}\label{rk:section atlas from dircon}
Let $\dircon$ be a direct connection on $\calG$ defined on a diagonal domain $\calU$, let $(U_\alpha)_{\alpha\in A}$ be an open cover of $M$ with $\bigcup_{\alpha\in A}(U_\alpha\times U_\alpha)\subset \calU$, and fix $x_0\in M$. For every $\alpha\in A$ choose, once and for all, a point $x_\alpha\in U_\alpha$ together with an element $\xi_\alpha\in \calG^{x_0}_{x_\alpha}$, which is non-empty by Proposition \ref{prop:direct}. Then
\begin{align}\label{eq:section atlas from dircon}
\tau_\alpha:U_\alpha\longrightarrow \calG^{x_0}_{U_\alpha}, \qquad \tau_\alpha(x):=\dircon(x,x_\alpha)\,\xi_\alpha
\end{align}
is a section atlas based at $x_0$: since $x$ and $x_\alpha$ both lie in $U_\alpha$ and $U_\alpha\times U_\alpha\subset \calU$, the arrow $\dircon(x,x_\alpha)\in \calG^{x_\alpha}_{x}$ is defined and depends smoothly on $x$, and composing it with $\xi_\alpha$ gives a smooth local section of $\tt^{x_0}:\calG^{x_0}\to M$ with $\tt^{x_0}\circ \tau_\alpha=\Id_{U_\alpha}$.
\end{remark}


If the manifold $M$ is not connected, one can work separately on its connected components. 
Thus, from now on we only consider connected base manifolds $M$ and gauge groupoids $\calG(P)\toto M$ associated to a principal $G$-bundle $P\to M$.

\begin{prop}\label{prop:tau=dir conn}
Any principal connection on $P$ defines a direct connection on $\calG(P)$ through its parallel transport along small geodesics on $M$.
Such a direct connection is transitive, and it is flat if and only if the principal connection on $P$ is. 
\end{prop}

\begin{proof}
For any two points $x$ and $y$ of $M$ joint by a unique geodesic, let $\Par^P(y,x):P_x\to P_y$ denote the parallel transport on $P$ along this geodesic. 
The direct connection $\dircon^\Par$ is defined on the pair $(y,x)$ by 
\begin{align}\label{eq:tau=Gamma}
\dircon^\Par(y,x) = [\Par^P(y,x)p,p] \quad\mbox{for any choice of $p\in P_x$}.
\end{align}
The map $\dircon^\Par$ is well defined on the convex diagonal domain
\begin{align*}
    \calU = \{(y,x)\in M\times M\ |\ \text{there exists a convex open subset of $M$ containing both $x$ and $y$}\}, 
\end{align*}
and has values in $\calG(P)$ because for any other element $p'\in P_x$ there exists a unique element $g$ in $G$ such that $p'=pg$, and the equality $[\Par^P(y,x)p',p'] = [\Par^P(y,x)p,p]$ follows from the $G$-equivariance of $\Par^P(y,x)$.

The direct connection $\dircon^\Par$ is transitive, because parallel transports are. It follows from the Ambrose-Singer theorem \cite{Ambrose-Singer:1953}, \cite[Theorem 8.1]{Kobayashi-Nomizu:1963}, that $\dircon^\Par$ is flat if and only if the principal connection on $P$ inducing $\Par^P$ is flat. 
\end{proof}

Principal connections always exist on principal bundles \cite[Ch.II, Proof of Theorem 2.1]{Kobayashi-Nomizu:1963},
therefore any gauge groupoid on a base manifold $M$ equipped with a linear connection admits a direct connection. 

\begin{example}[Parallel transport on the trivial groupoid induced by a global section]\label{ex:tau on MxG}
Example \ref{ex:trivial dircon} corresponds to the case of a trivial principal bundle $P=M\times G$ with global section $\tau$. Each fibre $P_x$ is diffeomorphic to $G$ through $\tau(x)\mapsto 1_G$, so any $p\in P_x$ can be written as $p=\tau(x)g$ for some $g\in G$. The $\mathfrak{g}$-valued $1$-form on $P$ given by $\omega(\tau(x)g):= d(\tau(x)g)\,(\tau(x)g)^{-1}$ yields the principal connection on $P$ with parallel transport
\begin{align*}
\Par^P(y,x) \tau(x)g &= \tau(y)g.  
\end{align*}
The induced direct connection $\dircon_g$ on $\calG(P)=M\times G\times M$ is precisely that of Example \ref{ex:trivial dircon}, namely $\dircon_\tau(y,x) = (y,\tau(y)\tau(x)^{-1},x)$. 
\end{example}


\subsection{Direct connections on frame groupoids}
\label{sec:direct connections on frame groupoids}

When the principal bundle is the frame bundle $FE$ of a vector bundle $E\to M$, direct connections on the frame groupoid $\GL(E)=\calG(FE)\toto M$ are called \textbf{linear direct connections} on $E$ by N. Teleman and J. Kubarski \cite{Teleman:2004,Teleman:2007,Kubarski-Teleman:2006}.

\begin{example}[Parallel transport on a vector bundle] \label{ex:tau^E}
Let $E\to M$ be a vector bundle with a linear connection $\nabla$. Then the parallel transport $\Par^E$ of $\nabla$ along small geodesics of $M$ is a linear direct connection on the frame groupoid $\GL(E)$ \cite{Teleman:2004}.
The linear connection $\nabla$ is equivalent to a principal connection on the frame bundle $FE$, and $\Par^E$ is precisely the direct connection induced on $\GL(E)$ by the parallel transport $\Par^{FE}$ in $FE$: for any $x\in M$ and any $y\in U_x$, the linear isomorphism 
\begin{align*}
\Par^E(y,x):E_x\to E_y
\end{align*}
is the composite $\varphi_y\circ \varphi_x^{-1}$ where $\varphi_x:\R^r\to E_x$ is any frame at $x$ (assuming that $E$ has rank $r$) and $\varphi_y=\Par^{FE}(y,x)(\varphi_x):\R^r\to E_y$ is its parallel transported frame at $y$.
\end{example}

We borrow from \cite[Example 5]{Kubarski-Teleman:2006} an easy and illustrative example of a direct connection which is not induced by parallel transport. 

\begin{example}[Direct connections which are not parallel transports] \label{ex:dircon on RxR}
On the base manifold $M=\R$, with global vector field $\partial_x=\frac{d\ }{dx}$, let us fix the flat linear connection 
whose geodesics are the segments parametrised by $x$. 
Let $E=M\times \R$ be the trivial bundle, with global section $e_1(x)=(x,1)\in E_x$. A linear connection $\nabla^E:\dircon(TM)\otimes \dircon(E)\to \dircon(E)$ on $E$ (necessarily flat) is given by its Christoffel symbol $\lambda\in C^\infty(M)$ such that $\nabla^E_{\partial_x} e_1(x) = -\lambda(x)\,e_1(x)$. 
The induced parallel transport along geodesics defines a global direct connection $\Par^E:\Pair(M)\to \GL(E)$ by 
\begin{align} \label{eq:tau(y,x)}
\Par^E(y,x):E_x \to E_y,\ e_1(x)\mapsto \Par (y,x)\,e_1(x) = e^{\Lambda(y)-\Lambda(x)}\,e_1(y),
\end{align}
where $\Lambda(x)=\int \lambda(x)\,dx$ is a primitive function of the Christoffel symbol. 
Note that $e^{\Lambda(y)-\Lambda(x)}=\sigma(y)\sigma(x)^{-1}=\Par^{FE}(y,x)$ is also the parallel transport on the frame bundle of $E$, which is the trivial principal $GL_1(\R)$-bundle as in Example \ref{ex:tau on MxG}, with global section $\sigma(x)=\pm e^{\Lambda(x)}$. 


In contrast, the smooth maps of the form 
\begin{align}\label{eq:Gamma not tau}
\dircon(y,x)\,e_1(x) &= e^{\Lambda(y)-\Lambda(x)+(y-x)^2\,h(y,x)}\,e_1(y), 
\end{align}
where $\Lambda\colon\R\to\R$ and $H\colon \R^2\to\R$ are any smooth functions, and $H$ is not identically zero, define all possible global direct connection on $\GL(E)$ which are not a parallel transport.
In fact, if $\dircon(y,x)=e^{F(y,x)}$ is a global smooth direct connection, the condition $\dircon(x,x)=1_x$ and its consequence $(\frac{\partial \dircon}{\partial y}+\frac{\partial \dircon}{\partial x})\vert_{y=x}=0$ \cite[Proposition 6]{Kubarski-Teleman:2006} determine $F$ at the 1st order of expansion in $y-x$, in terms of $\Lambda(x)=\int \frac{\partial \lambda}{\partial y}\vert_{y=x}\, dx = -\int \frac{\partial \lambda}{\partial x}\vert_{y=x}\, dx$.
Indeed, if $\tilde{F}$ is of the form $\tilde{F}(y,x)=\Lambda(y)-\Lambda(x)+(y-x)\tilde{h}(y,x)$, with $\tilde{H}:\R^2\to \R$ a non-zero smooth function, its expansion at $y=x$ gives
\begin{align*}
\frac{\partial\tilde{h}}{\partial y}\vert_{y=x} &= \Lambda'(x)+\tilde{h}(x,x). 
\end{align*}
Then, if $\tilde{h}(x,x)= 0$ for all $x\in M$, one can write $\tilde{h}(y,x)=(y-x)\,h(y,x)$ and the direct connection 
\begin{align*}
\tilde{\dircon}(y,x)\,e_1(x) &= e^{\Lambda(y)-\Lambda(x)+(y-x)\,\tilde{h}(y,x)}\,e_1(y)
\end{align*}
can be put in the form \eqref{eq:Gamma not tau} of $\dircon$. Instead, if $\tilde{h}(x,x)\neq 0$ for some $x\in M$, then $\tilde{\dircon}$ is again of the same form \eqref{eq:Gamma not tau} for $\tilde{\Lambda}(x)=\Lambda(x)+\int \tilde{h}(x,x)\,dx$ related to a modified linear connection 
\begin{align*}
\tilde{\nabla}^E_{\partial_x} e_1(x) = -\big(\lambda(x)+\tilde{h}(x,x)\big)\, e_1(x). 
\end{align*}
\end{example}


Direct connections on frame groupoids are also related to the notion of parallelism used in \cite[Definition 86]{DahlqvistDiehl-Driver:2019}. A \textbf{parallelism} on $E$ is a smooth local section $U:M\times M\locmap \Hom(E)$ of the bundle $\Hom(E)\cong E\boxtimes E^* \to M\times M$ of linear maps between fibres of $E$, defined on a diagonal domain and such that $U(x,x)=\Id_{E_x}$.  
Of course a direct connection on $\GL(E)$ is a parallelism.
The converse also holds. 

\begin{prop}\label{prop:paralleldirconn}
A parallelism on $E$ yields a direct connection on $\GL(E)$. 
\end{prop}

\begin{proof}
A parallelism $U$ on $E$ sends a pair $(y,x)\in M\times M$ to a linear map $U(y,x):E_x\to E_y$. All we need to prove is that this map is invertible if $x$ and $y$ are sufficiently close.

Suppose that $E$ has rank $r$, choose a neighbourhood $U_x$ of $x$ trivialising $E$, so that $E_x\cong \R^r$ and $E_y\cong \R^r$ for any $y\in U_x$. Then compose the map $U(\_,x):U_x\to \End(\R^r)\cong M_r(\R)$ with the determinant map $\det:M_r(\R)\to \R$. The real valued map $\det U(\_,x)$ is continuous and has value $1$ on $x$, because $U(x,x)=\Id_{E_x}$. So there exists an open neighbourhood $V_x$ of $x$ in $U_x$ where $\det U(\_,x)>0$. This means that $U(y,x)$ is invertible for any $y\in V_x$, and therefore $U$ is a direct connection on $\GL(E)$ defined on the diagonal domain $\{(y,x)\ |\ x\in M, y\in V_x\}$.
\end{proof}

With this result, one can prove in an easy direct way that frame groupoids always have direct connections. 

\begin{remark} \label{rk:dirconframe}
On a frame groupoid, local trivialisations induce a direct connection. 

Let $\{(U_\alpha,\theta_\alpha)\}_{\alpha\in A}$ be a trivializing atlas of $E$, with trivialisations $\theta_\alpha:E_{U_\alpha}\to U_\alpha\times \R^r$. For any $x,y\in U_\alpha$, set 
\begin{align*}
\dircon_\alpha(y,x) &:= \left({\theta_\alpha}\vert_{E_y}\right)^{-1} \circ {\theta_\alpha}\vert_{E_x} \in \GL(E_x,E_y). 
\end{align*}
Then choose a family of smooth functions $\chi_\alpha$ such that  $\{\chi^2_\alpha\}_{\alpha\in A}$ is a partition of unity 
subordinated to $\{U_\alpha\}_{\alpha\in A}$ and glue the maps $\dircon_\alpha$ as
\begin{align*}
\dircon(y,x) &:= \sum_{\alpha\in A} \chi_\alpha(y) \dircon_\alpha(y,x) \chi_\alpha(x). 
\end{align*}
This gives a parallelism on $E$ defined on the diagonal domain $\bigcup_{\alpha\in A} U_\alpha\times U_\alpha$, and therefore a direct connection on $\GL(E)$. 
\end{remark}



\subsection{Infinitesimal connections induced by direct connections} 
\label{subsec:infinitesimal connection induced by direct connections}

A \textbf{connection on a Lie algebroid $A\to M$} is a splitting of the anchor $a:A\to TM$, i.e. a vector bundle map $\delta: TM\to A$ such that $a\circ \delta=\Id_{TM}$, cf. \cite[Definition 5.2.5]{Mackenzie:2005}.
It is called \textbf{flat} if it is a morphism of Lie algebroids, i.e.
$[\delta(X),\delta(Y)]_A = \delta([X,Y])$
for any $X,Y\in \Gamma(TM)$.

The existence of a connection on $A$ forces $A$ to be transitive (the anchor is surjective), however this does not imply that $A$ is the Atiyah algebroid $A(P)=TP/G$ of a principal $G$-bundle $P\to M$, cf. Remark \ref{rem:transitive Lie algebroid}.
If we assume that $A=A(\calG)$ comes from a Lie groupoid, then the existence of a connection on $A$ implies that $\calG$ is a gauge groupoid and $A=A(P)$ is an Atiyah algebroid. 
In this case, a connection on $A(P)$ is also called an \textbf{infinitesimal connection on $P$} and it is equivalent to a principal connection on $P$ \cite[Ch. IV]{Kobayashi:1957}.
In fact, $\delta$ gives a splitting of the Atiyah sequence given in \eqref{eq:Atiyahseq},
\begin{equation*}
\begin{tikzcd}
0 \arrow[r] & \arrow[l, bend right=15] \Ad(P) \arrow[r, "\iota" below] & \arrow[l, bend right=15, "\lambda" above] A(P) \arrow[r, "a" below] & \arrow[l, bend right=15, "\delta" above] TM \arrow[r] & \arrow[l, bend right=15] 0
\end{tikzcd}
\end{equation*}
where $\lambda: A(P)\to \Ad(P)$ is a bundle map such that $\iota \circ \lambda + \delta \circ a = \Id_{A(P)}$, that is, $\lambda(\xi)=\xi-\delta(a(\xi))$ for any $\xi\in A(P)$ and in particular $\lambda$ is the identity on vertical vectors $v\in \Ad(P)$ \cite[\S 5.2, 5.3]{Mackenzie:2005}. 
Lifting to $TP$ by $\natural^{-1}$ yields a splitting $TP\cong VP \oplus HP$ into vertical and horizontal subbundles, with 
$$
VP=\ker(d\pi)=\natural^{-1}\Ad(P) \qquad\text{and}\qquad HP=\natural^{-1}\delta(TM).
$$
In other words, given $\delta$ and $\lambda$, we recover the \emph{horizontal lift} of any $X\in TM$ to $TP$ as $\natural^{-1}\delta(X)$, and a principal connection 
This motivates the alternative name of the Atiyah algebroid $A(P)$ in the literature, as the \textbf{bundle of connections} of $P$, cf. \cite[Ch. IV]{Kobayashi:1957}. 

\begin{example}[Connection on the trivial Lie algebroid]
If $P=M\times G$ is trivial, then $A(P)=TM \oplus (M\times \mathfrak{g})$ is the trivial Lie algebroid, with $\Ad(P)=M\times \mathfrak{g}$ (cf. Example \ref{exs:Lie algebroids} (c)). 
Then, an infinitesimal connection $\delta:TM\to A(P)$ splits a vector $X+v\in TM\oplus (M\times \mathfrak{g})$ as $X+v = \delta(X) + \lambda(X+v)$ with 
$$
\lambda(X+v)=A(X)+v \qquad\text{and}\qquad \delta(X)=X-A(X), 
$$
where $-A(X)\in M\times \mathfrak{g}$ represents the \emph{vertical} component of the horizontal lift of $X$ in $TP$ fixed by $\delta$, with respect to the given trivialisation of $P$, and 
$$
A=\lambda|_{TM}:TM\to M\times \mathfrak{g},\ X\mapsto A(X)=X-\delta(X)
$$ 
gives the connection $1$-form on $M$ representing the principal connection $\omega$ on $P$ in the given trivialisation, that is, a usual \emph{gauge boson}. 
\end{example}

\begin{defn} \label{def:conngroup}
Following the classical terminology, we call \textbf{infinitesimal connection on a Lie groupoid $\calG\toto M$} a connection $\delta:TM\to \Lie(\calG)$ on its Lie algebroid.
The existence of such a connection implies that $\Lie(\calG)=A(P)$ is an Atiyah algebroid, and therefore, again, that $\calG=\calG(P)$ is a gauge groupoid.
\end{defn}

\begin{prop} \label{prop:dircttoinf}
Let $\calG(P)\toto M$ be a gauge groupoid endowed with a direct connection $\dircon:\Pair(M)\locmap \calG(P)$. Then the $\tt$-differential of $\dircon$ along the diagonal gives an infinitesimal connection $\delta^\dircon = d^\tt \dircon\vert_M : TM\to A(P)$ defined as 
\begin{align}\label{eq:delta^Gamma}
\delta^\dircon(X) & = \frac{d}{dt}\dircon(c(t),x)\vert_{t=0}, 
\end{align}
where $X\in T_xM$ and $c:[0,1]\to M$ is the integral curve of $X$ at $x$.
Moreover, if $\dircon$ is flat (it is a local groupoid morphism) then $\delta$ is flat (it is an algebroid morphism). 
\end{prop}

\begin{proof} 
Apply Proposition \ref{prop:t-differential} to the local map $\dircon:\Pair(M)\locmap \calG$. Then, the vector $\delta^\dircon(X) = d^\tt\dircon\vert_x(\dot{c}(0))$ certainly belongs to $\Lie(\calG(P))\cong A(P)$.
Moreover, $\delta^\dircon$ yields a splitting of the anchor $a=d\tt\vert_{T_{1_{x}}\calG^{x}}$ because, for $X=\dot{c}(0)\in T_xM$, we have
\begin{align*}
  a(\delta^\dircon(X)) &= a\left(\frac{d\ }{dt}\dircon(c(t),x)\vert_{t=0}\right)
  = \frac{d\ }{dt} \tt\left(\dircon(c(t),x)\right)\vert_{t=0}
  = \frac{d\ }{dt} c(t)\vert_{t=0} = X. 
\end{align*}
The last statement follows from the fact that the tangent map along the diagonal of a local groupoid morphism is a Lie algebroid morphism \cite[\S 3.5]{Mackenzie:2005}. 
\end{proof}
 
\begin{cor} \label{cor:deltaIsoE} 
With the above notations, if $\dircon$ is a direct connection on a frame groupoid $\GL(E)\toto M$ or on a reduction $\calG(P)\rightarrow \GL(E)$, the local expression
\begin{align}\label{eq:nablaGamma}
 \nabla^\dircon_X \psi\vert_x &= d^\ss\dircon(X)\psi\vert_x = \frac{d}{dt}\dircon(c(t),x)^{-1}\,\psi(c(t))\vert_{t=0}, 
\end{align}
where $\psi$ is a local smooth section of $E$ around $x\in M$, defines a linear connection on $E$. Together with the infinitesimal connection of \eqref{eq:delta^Gamma}, it yields a decomposition of the tangent map $d\psi\vert_x:T_xM\to T_{\psi(x)}E$ into its horizontal and vertical components
\begin{align}\label{eq:d=delta+nabla} 
d\psi\vert_x(X) &= \delta^\dircon(X)\psi(x) + \nabla^\dircon_X \psi\vert_x, 
\end{align}
where $\delta^\dircon(X) \in T^\mathrm{lin}_xE$ acts linearly on $\psi(x)\in E_x$, and $\nabla^\dircon_X\in \Der_x(E)$ is a linear derivation of $\psi$ at $x$. 
\end{cor}

\begin{proof}
Apply Remark \ref{rem:s-differential}: the expression \eqref{eq:nablaGamma} can be obtained as the $\ss$-differential of $\dircon$ along the diagonal, together with the evaluation of $\GL(E)$ on $E$ of Example \ref{ex:groupoidactions} (a), in the same way as in Proposition \ref{prop:dircttoinf}, and finally using the equality 
$$
\frac{d}{dt}\dircon(x,c(t))\,\psi(c(t))\vert_{t=0} = \frac{d}{dt}\dircon(c(t),x)^{-1}\,\psi(c(t))\vert_{t=0}
$$
of Lemma \ref{lem:d/dt Gamma}. 
Alternatively, this follows from the general relationship between a horizontal and vertical components of the differential of $\psi$, as in Remark \ref{rk:A'(G)}.
In particular, \eqref{eq:d=delta+nabla} is a reformulation of \eqref{eq:d=xi+D}: 
\begin{align*}
  d\psi\vert_x(X) &= \frac{d}{dt} \psi(c(t))\vert_{t=0} 
  = \frac{d}{dt} \dircon(c(t),x)\,\dircon(c(t),x)^{-1}\, \psi(c(t))\vert_{t=0} \\
  &= \frac{d}{dt} \dircon(c(t),x)\vert_{t=0}\,\dircon(x,x)^{-1}\,\psi(x)
  + \dircon(x,x)\,\frac{d}{dt}\dircon(c(t),x)^{-1}\, \psi(c(t))\vert_{t=0} \\
  &= \delta^\dircon (X) \psi(x) + \nabla^\dircon_X \psi\vert_x.
\end{align*}
\end{proof}

\begin{lem} \cite[Proposition 6]{Kubarski-Teleman:2006} \label{lem:d/dt Gamma}
Any direct connection $\dircon:\Pair(M)\locmap \calG$ satisfies the infinitesimal identities
\begin{align*} 
  \frac{d\ }{dt}\dircon(c(t),x)\vert_{t=0} &= -\frac{d\ }{dt}\dircon(x,c(t))\vert_{t=0} \\
  \frac{d\ }{dt}\dircon(c(t),x)^{-1}\vert_{t=0} &= -\frac{d\ }{dt}\dircon(c(t),x)\vert_{t=0} = \frac{d\ }{dt}\dircon(x,c(t))\vert_{t=0}
\end{align*}
\end{lem}

\begin{proof}
The first identity is a restatement of the fact that the directional derivative of $\dircon(y,x)$ along the diagonal vanishes, i.e. $\big(\frac{\partial}{\partial x^i}+ \frac{\partial}{\partial y^i}\big)\dircon(y,x)\vert_{y=x}=0$, which follows from $\dircon(x,x)=1_x$.
The second identity is the derivative of $\dircon(c(t),x)^{-1}\dircon(c(t),x)=\Id_{x}$. 
\end{proof}

A consequence of Proposition \ref{prop:dircttoinf} and Corollary \ref{cor:deltaIsoE} is that one can recover the usual Christoffel symbols of a linear connection on $E$ in terms of the matrix elements of a direct connection on $\GL(E)$, cf. also \cite[Theorem 3]{Kubarski-Teleman:2006}.
If $\{e_\alpha(x)\}_{\alpha=1}^r$ is a basis of local sections of $E$ at $x\in M$, and $\psi(x)=\psi^\alpha(x)\, e_\alpha(x)$ is a local section (with Einstein summation convention), then a direct connection $\dircon$ on $E$ is locally given by matrix elements $\dircon^\beta_\alpha(y,x) \equiv \dircon(y,x)^\beta_\alpha$ such that
\begin{align*}
\dircon(y,x)\big(\psi^\alpha(x)\,e_\alpha(x)\big) &= \psi^\alpha(x)\,\dircon^\beta_\alpha(y,x)\,e_\beta(y). 
\end{align*}
If $\{\partial_i=\frac{\partial}{\partial x^i}\vert_x\}_{i=1}^d$ is a basis of local sections of $T_xM$, the expression 
\begin{align*}
\nabla^\dircon_{\partial_i} e_\alpha\vert_{x}
&= \partial_i\dircon^\beta_\alpha(x,y)\vert_{y=x}\,e_\beta(x)
:= \frac{\partial\ }{\partial y^i} \dircon^\beta_\alpha(x,y)\vert_{y=x}\, e_\beta(x)
\end{align*}
defines the usual \textbf{Christoffel symbols} of the linear connecion $\nabla^\dircon$ on $E$ as
\begin{align} \label{eq:christoffel}
\dircon^\beta_{i \alpha}(x) &:= \partial_i\dircon^\beta_\alpha(x,y)\vert_{y=x}. 
\end{align}
Therefore, by \eqref{eq:d=delta+nabla}, for $X=X^i\,\partial_i\in T_xM$ we have
\begin{align} \nonumber
  \nabla^\dircon_X \psi\vert_x &= X^i\,\Big(\partial_i \psi^\beta(x) + \psi^\alpha(x)\,\partial_i \dircon^\beta_\alpha(x,y)\vert_{y=x}\Big) e_\beta(x)
  \\ \label{eq:connection in coordinates}
  \delta(X) \psi(x) &= -X^i\,\psi^\alpha(x)\,\partial_i \dircon^\beta_\alpha(x,y)\vert_{y=x}\,e_\beta(x)
  \\ \nonumber
  & = X^i\,\psi^\alpha(x)\,\partial_i \dircon^\beta_\alpha(y,x)\vert_{y=x}\,e_\beta(x).
\end{align}

\begin{prop} \label{prop:nabla^Gamma^omega=omega}
Let $\delta:TM\to A(P)$ be an infinitesimal connection on $P$ and let $\dircon^\delta$ be the direct connection on $\calG(P)$ defined by the parallel transport of $\delta$ along small geodesics, as in Proposition \ref{prop:tau=dir conn}.
Then, the infinitesimal connection defined by $\dircon^\delta$ coincides with $\delta$, that is, 
\begin{align}\label{eq:deltaGammadelta}
  \delta(X) = \frac{d}{dt}\dircon^\delta(c(t),x)\vert_{t=0}, 
\end{align}
where $c:[0,1]\to M$ is the integral curve of $X$ at $x$. 
\end{prop}

\begin{proof}
Let $\Par$ be the parallel transport of $\delta$ along small geodesics. Then $\dircon^\delta:\Pair(M)\locmap \calG(P)$ is given by $\dircon^\delta(y,x)=[\Par(y,x)(p),p]$ for any choice of $p\in P_x$, and its infinitesimal connection is 
$$
\delta^{\dircon^\delta}(X) = \frac{d}{dt}\dircon^\delta(c(t),x)\vert_{t=0}
= \frac{d}{dt} [\Par(c(t),x)(p),p]\vert_{t=0}.
$$
By uniqueness of the horizontal lift, the tangent vector at $t=0$ to the curve $\gamma(t)=\Par(c(t),x)(p)$ gives the value of $\delta(\dot{c}(0))$. 
Taking into account the action of $G$ on $TP$, which makes the choice of $p$ irrelevant, one can identify $ {\delta}^{\dircon^\delta}(\dot{c}(0))$ to the class of $\dot{\gamma}(0)$ in $TP/G$, from which it follows that
$$
\delta^{\dircon^\delta}(\dot{c}(0)) = \natural_p(\dot{\gamma}(0)) = \delta(\dot{c}(0)), 
$$
for any $\dot{c}(0)\in T_xM$.
\end{proof}
 
The analogue statement holds for a linear connection $\nabla$ on $E$ and the direct connection $\dircon^\nabla$ on $\GL(E)$ defined by its parallel transport.

The above result shows that infinitesimal connections can always be integrated to direct connections: it suffices to consider their parallel displacement. 
However, the correspondence between direct connections and infinitesimal connections is not one-to-one, as shown by next example. 

\begin{example}[Direct connections with equal infinitesimal connections] \label{ex:Gamma with equal delta}
Consider the trivial vector bundle $E=M\times \R$ on the flat base manifold $M=\R$, as in Example \ref{ex:dircon on RxR}. The two direct connections on $\GL(E)$ given by \eqref{eq:tau(y,x)} (a parallel transport) and by \eqref{eq:Gamma not tau}, namely 
\begin{align*} 
\Par(y,x)\,e_1(x) &= e^{\Lambda(y)-\Lambda(x)}\,e_1(y) \\
\dircon(y,x)\,e_1(x) &= e^{\Lambda(y)-\Lambda(x)+(y-x)^2\,h(y,x)}\,e_1(y), \qquad H\neq 0
\end{align*}
have the same infinitesimal connection 
$\delta(\partial_x)=\Lambda'(x)=\lambda(x)$. 
\end{example}

In fact, we can fully characterise the direct connections which are parallel transport: they are exactly the transitive ones. 

\begin{thm} \label{thm:transitive = parallel transport}
Let $\dircon$ be a direct connection on a frame groupoid $\GL(E)$ defined on a convex diagonal domain $\calU$. 
If $\dircon$ is transitive, then it is the parallel transport of its associated linear connection $\nabla^\dircon$ on $E$.
\end{thm}

\begin{proof}
For any $(y,x)\in \calU$, let $c:[0,1]\to M$ be the geodesic from $x$ to $y$. The direct connection $\dircon$ restricted to $c$ gives a linear map $\dircon(c(t),c(s)): E_{c(s)}\to E_{c(t)}$ for any $s,t\in [0,1]$. 
Its infinitesimal connection $\delta^\dircon$, restricted to $c$, is given by \eqref{eq:delta^Gamma} evaluated at every point $c(t)$, namely
\begin{align*}
    \delta^\dircon(\dot{c}(t)) = \frac{d}{du} \dircon(c(u),c(t))\big|_{u=t}. 
\end{align*}
Since $\dircon$ is transitive, for any $s,t\in [0,1]$ we have the composition property
\begin{align*} 
    \dircon(c(u),c(t))\, \dircon(c(t),c(s)) = \dircon(c(u),c(s)). 
\end{align*}
Differentiating it in the variable $u$ at $u=t$, we see that $\dircon$
satisfies the linear differential equation 
\begin{align*}
    \frac{d}{du} \dircon(c(u),c(s))\big|_{u=t} &= \frac{d}{du} \dircon(c(u),c(t))\big|_{u=t}\ \dircon(c(t),c(s)) 
    =  \delta^\dircon(\dot{c}(t))\ \dircon(c(t),c(s)), 
\end{align*}
with initial condition $\dircon(c(t),c(t))=1_{c(t)}$. 

Now, the parallel transport induced by $\dircon$ along the geodesic $c$ is also a linear map $\Par^\dircon(c(t),c(s)): E_{c(s)}\to E_{c(t)}$ which satisfies the same differential equation
\begin{align*}
    \frac{d}{du}\Par^\dircon(c(u),c(s))\big|_{u=t} &= \delta^\dircon(\dot{c}(t))\ \Par^\dircon(c(t),c(s)), 
\end{align*}
with the same initial condition $\Par^\dircon(c(t),c(t))=1_{c(t)}$. 
In fact, by definition, the parallel transport $\Par^\dircon(c(t),c(s))$ is the linear map such that, for any $e\in E_{c(s)}$, the section $t\mapsto \psi(c(t))=\Par^\dircon(c(t),c(s))(e)$ along $c$ is parallel, i.e. 
\begin{align*}
    \nabla^\dircon_{\dot{c}(t)} \psi\big|_{c(t)} = d\psi(\dot{c}(t)) - \delta^\dircon(\dot{c}(t))\ \psi(c(t)) = 0. 
\end{align*}
Since $d\psi(\dot{c}(t))=\frac{d}{du}\Par^\dircon(c(u),c(s))\big|_{u=t}(e)$, this gives the above equation. 

By uniqueness of the solution of ODEs, then, a transitive direct connection $\dircon$ coincides with the parallel transport $\Par^\dircon$ of its infinitesimal connection on every points along the geodesic from $x$ to $y$, and therefore they coincide on the convex diagonal domain $\calU$ in $\Pair(M)$. 
\end{proof}

\begin{remark} \label{rem:one-to-one}
The fact that direct connections are \emph{not} in one-to-one correspondence with infinitesimal connections constrasts with \emph{path connections} on groupoids, cf. \cite[Definition 6.3.1]{Mackenzie:2005}. A path connection is a horizontal lift to the groupoid of a curve on the base manifold $M$. This lift is by definition so constrained that its initial direction (which is an infinitesimal connection on its Lie algebroid) determines it all along the path on $M$ \cite[Theorem 6.3.5]{Mackenzie:2005}, and finally identifies it with the exponential map from Lie algebroids to Lie groupoids \cite[Corollary 6.3.6]{Mackenzie:2005}.

The situation changes if we consider \emph{flat} connections: there is a one-to-one correspondence between \emph{flat} direct connections, i.e. local groupoid morphisms, and \emph{flat} infinitesimal connections, i.e. algebroid morphisms. This is proved by K. Mackenzie in \cite[Corollary 6.2.7]{Mackenzie:2005} and can be seen as a consequence of \cite[Proposition 6.1.8]{Mackenzie:2005}, that we recall in Section \ref{subsec:flat groupoids}.  
\end{remark}

We conclude this section by collecting the characterisation of flat connections in classical terms and in terms of groupoids, which follows from Propositions \ref{prop:tau=dir conn}, \ref{prop:dircttoinf} and Corollary \ref{cor:deltaIsoE}. 

\begin{cor} \label{cor:flat-equivalence}
Let $M$ be a smooth manifold endowed with a linear connection. 
If $\omega$ is a principal connection on a principal bundle $P\to M$ and $\dircon$ is the direct connection on $\calG(P)\toto M$ induced by the parallel transport on $P$, then $\omega$ is flat as a principal connection if and only if $\dircon$ is flat as a direct connection. 
\end{cor}


\subsection{Curvature of direct connections}
\label{subsec:curvature}


The obstruction of a direct connection to be flat is measured, as usual, by a \emph{curvature}. Our definition slightly differs from that given by N. Teleman and J. Kubarski in \cite{Teleman:2004,Teleman:2007,Kubarski-Teleman:2006,Kubarski:2008} and by A. Kock in \cite{Kock:2007,Kock:2017} but the notions are ssentially equivalent.


\begin{defn}\label{def:curvature} 
Let $\dircon: \Pair(M)\locmap \calG$ be a direct connection defined on a symmetric diagonal domain $\calU\subset M\times M$. For any $x\in M$, we set
$$
\calU(x) = \{ (z,y)\in \calU\ |\ (y,x),(z,x)\in \calU\}
$$
and call \textbf{curvature of $\dircon$ at $x$} the map $K_x: \calU(x)\longrightarrow \calG_x^x$ given on $(z,y)\in \calU(x)$ by 
\begin{align}\label{eq:curvature}
K_x(z,y) := \dircon(z,x)^{-1}\,\dircon(z,y) \dircon(y,x). 
\end{align}
Then $\dircon$ is flat, according to Definition \ref{def:directconnection}, if $K_x=1_x$ for any $x\in M$, and it is transitive if $K_x=1_x$ along geodesics starting at $x$. 
The curvature is \textbf{symmetric} if $K_x(z,y)=K_x(y,z)$ for any points, that is, $K_x(y,z)^{-1} K_x(z,y)=1_x$, and it is \textbf{antisymmetric} if $K_x(z,y)=K_x(y,z)^{-1}$, that is, $K_x(y,z) K_x(z,y)=1_x$.
Note that, in general, we have
\begin{align*}
K_x(y,z) K_x(z,y) &= \dircon(y,x)^{-1} \dircon(y,z) \dircon(z,y) \dircon(y,x) \\
&= \Ad_{\dircon(y,x)^{-1}}\left(\dircon(y,z) \dircon(z,y)\right), 
\end{align*}
therefore $K$ is antisymmetric if and only if the direct connection $\dircon$ is. 
\end{defn}

\begin{example}[Curvature of a non-infinitesimal connection] \label{ex:curvature on RxR}
For the general direct connection $\dircon$ on $\R\times \R\to \R$ of Example \ref{ex:dircon on RxR}, given by \eqref{eq:Gamma not tau}, the curvature is: 
\begin{align*}
    K_x(z,y)=e^{(z-y)^2h(z,y)+(y-x)^2h(y,x)-(z-x)^2h(z,x)}.
\end{align*}
\end{example}

The invariants of a direct connection $\dircon$ are carried by the curvature operator $K_x$ and its derivatives. 

\begin{lem} \label{lem:K'=0}
The first derivatives of $K_x$ along a geodesic $y(t)$ starting at $x$ vanish: 
\begin{align*} 
    \frac{d}{dt} K_x(x,y(t))\big|_{t=0} = \frac{d}{dt} K_x(y(t),x)\big|_{t=0} = \frac{d}{dt} K_x(y(t),y(t))\big|_{t=0} = 0.
\end{align*}
\end{lem}

\begin{proof}
The first derivative vanishes because of Lemma \ref{lem:d/dt Gamma}. 
The second and the third ones because $K_x(y(t),x)=1_x$ and $K_x(y(t),y(t))=1_x$ are constant.
\end{proof}

Note, however, that the first derivatives of $K_x$ computed not along a  geodesic do not necessarily vanish, for instance
$$
\frac{d}{dt} K_x(z,y(t))\big|_{t=0}\;, \qquad \frac{d}{dt} K_x(y(t),z)\big|_{t=0}
$$
when $z\neq x$. 
Let us show that the second derivative of $K_x$ is related to the curvature tensor of the induced infinitesimal connections. 

Let $P$ be a principal bundle on $M$ with associated vector bundle $E$, so that the gauge groupoid $\calG(P)$ is a reduction of the frame groupoid $\GL(E)$ (for instance $P\subset FE$). 
Let $\dircon:\Pair(M)\locmap \calG(P)$ be a direct connection with infinitesimal connection $\delta:TM\to A(P)\subset T^{lin}E$ as in \eqref{eq:delta^Gamma} and linear connection $\nabla:TM\to \Der(E)$ as in \eqref{eq:nablaGamma}. We recall that the \textbf{curvature} of $\delta$ and of $\nabla$ at a point $x\in M$ is given by
\begin{align*}
  R^\delta_x(Z,Y)  \psi(x) &= [\delta(Z),\delta(Y)]_{A(P)} \psi(x) -\delta([Z,Y]) \psi(x) \\
  R^\nabla_x(Z,Y) \psi(x) &= [\nabla_Z,\nabla_Y]_{\Der(E)}\psi\vert_x-\nabla_{[Z,Y]}\psi\vert_x, 
\end{align*}
for any $Y,Z\in T_xM$ and for any local section $\psi$ of $E$ around $x$. In terms of the Lie algebroids, both curvatures are operators with values in the kernel of the anchor, that is, 
\begin{align*}
  R^\delta, R^\nabla:TM\times TM\to \End(E) \cong L(T^{lin}E)\cong L(\Der(E), 
\end{align*}
and measure how much the algebroid connections deviate from being algebroid morphisms. 
It was proven in \cite[Proposition 5.3.14]{Mackenzie:2005} that they coincide: $R^\delta_x(Z,Y)=R^\nabla_x(Z,Y)$. 

The relationship between the curvature $K$ of the direct connection $\dircon$ and the curvatures $R^\delta$ and $R^\nabla$ of its infinitesimal and linear connections was described by J.~Kubarski, N.~Teleman in \cite[Corollary 4]{Kubarski-Teleman:2006} for frame groupoids. We extend their result to general groupoids and give a comprehensive formula in terms of the operator which measures how much $K$ deviates from being symmetric. 

\begin{thm} \label{thm:curvature}
The curvature $R^\delta$ of the infinitesimal connection associated to a direct connection $\dircon$ is
\begin{align*}
R^\delta_x(Z,Y) &= \frac{\partial^2\,}{\partial s\, \partial t} \Big(K_x(y(t),z(s))^{-1}\, K_x(z(s),y(t))\Big)\vert_{t=s=0},
\end{align*}
where $Y,Z\in T_xM$ are two vectors tangent to $M$ at $x$, $y(t)=\exp_x(tY)$ is the integral curve of $Y$ and $z(s)=\exp_x(sZ)$ is the integral curve of $Z$. 
\end{thm}

As a consequence, a direct connection which is infinitesimally flat is not necessarily flat: it is only symmetric up to order $2$.
\bigskip

\begin{proof}
The proof is a simple computation making use of the Leibniz rule, the infinitesimal properties of $\dircon$ given in Lemma \ref{lem:d/dt Gamma} and the dictionary between direct connections and linear connections in terms of Christoffel symbols given in Eqs. \eqref{eq:christoffel} and \eqref{eq:connection in coordinates}. There is however a subtelty in the choice of curves along which we compute the derivatives of $\dircon$, let us give some details. 
For $\psi=\psi^\alpha\,e_\alpha$, $Y=Y^j\,\partial_j$, $Z=Z^k\,\partial_k$ and $\dircon_{i \alpha}^\beta(x)=\frac{\partial}{\partial y^i} \dircon_\alpha^\beta(x,y)\vert_{y=x}$, we have 
\begin{align*}
\nabla_Y\psi\vert_x &= Y^j\,\partial_j \psi^\beta(x)\,e_\beta(x) + Y^j\,\psi^\alpha(x)\,\frac{\partial\,}{\partial y^j} \dircon_\alpha^\beta(x,y)\vert_{y=x}\,e_\beta(x) \\
&= \frac{d}{dt}\psi(y(t))\vert_{t=0} + \frac{d}{dt}\dircon(x,y(t))\vert_{t=0}\,\psi(x). 
\end{align*}
The derivative of $\psi$ at $x$ disappears in the expression of the curvature $R^\delta_x(Z,Y)= R^\nabla_x(Z,Y)$, which acts linearly on $\psi(x)$, and can be neglected. It remains the linear part in $\psi(x)$, given by the derivative of $\dircon$.
Now, as we just showed, the first-order derivative of $\dircon(x,y)$ in the variable $y$, at $y=x$, gives the Christoffel symbol at $x$. In the curvature $R^\nabla(Z,Y)$, there appears a derivative of the Christoffel symbols in their (unique) variable $x$: this corresponds to a second-order derivative of $\dircon(x,y)$ in the variable $x$. Then, the standard expression of the curvature $R^\nabla(Z,Y)$ in terms of the Christoffel symbols gives
\begin{align*}
R^\delta_x(Z,Y)
&= Z^k\,\frac{\partial}{\partial z^k} \dircon(x,z)\vert_{z=x}\, Y^j\,\frac{\partial}{\partial y^j} \dircon(x,y)\vert_{y=x}
- Y^j\,\frac{\partial}{\partial y^j} \dircon(x,y)\vert_{y=x}\, Z^k\,\frac{\partial}{\partial z^k} \dircon(x,z)\vert_{z=x} \\ 
& \qquad + Z^k\,Y^j\,\frac{\partial^2\,}{\partial z^k\,\partial y^j} \dircon(z,y)\vert_{z=y=x}
- Y^j\,Z^k\,\frac{\partial^2\,}{\partial y^j\,\partial z^k} \dircon(y,z)\vert_{y=z=x}. 
\end{align*}
Note that this expression coincides with that given in \cite[Corollary 4]{Kubarski-Teleman:2006}.
To relate this to the formula given above, call $y=\exp_x(Y)$ and $z=\exp_x(Z)$ the end points of the geodesic $y(t)$ and $z(s)$ and consider the geodesic variation of $Y$ along $z(s)$ for small $s$, denoted $Y_s$, and the geodesic variation of $Z$ along $y(t)$ for small $t$, denoted $Z_t$, as in Figure \ref{fig:curvature} (cf. \cite[Ch. 6]{Lee:1997} for geodesic variations).  
\begin{figure}[h]
\centering
\begin{tikzpicture}[scale=1.5,font=\footnotesize]
\draw (0,0) node{$\bullet$} (-0.1,-0.1) node{$x$} ;
\draw [->] (0,0) -- (3,0) ;
\draw (3,0) node{$\bullet$} (3.1,0) node{$y$} (2.7,-0.2) node{$Y$} ;
\draw [->] (0,0) -- (1,2) ;
\draw (1,2) node{$\bullet$} (1.1,2.1) node{$z$} (0.7,1.8) node{$Z$} ;
\draw [->,dashed] (1,0) -- (1,2) ;
\draw (1,0) node{$\bullet$} (1,-0.2) node{$y(t)$} (1.2,1.8) node{$Z_t$} ;
\draw (1,0.7) node{$\bullet$} (1.3,0.7) node{$z_t(s)$} ;
\draw [->,dashed] (0.3,0.6) -- (3,0) ;
\draw (0.3,0.6) node{$\bullet$} (0.1,0.7) node{$z(s)$} (2.7,0.3) node{$Y_s$} ;
\draw (1.6,0.3) node{$\bullet$} (1.8,0.5) node{$y_s(t)$} ;
\end{tikzpicture}
\caption{Geodesic variations and related Jacobi fields}\label{fig:curvature}
\end{figure}
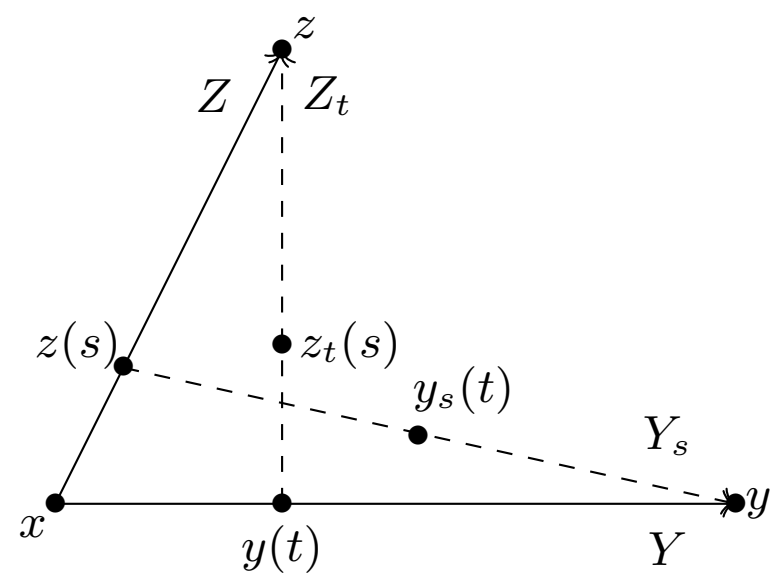 
The curvature can now be given in terms of the integral curves $y_s(t)=\exp_{z(s)}(tY_s)$ and $z_t(s)=\exp_{y(t)}(sZ_t)$ for small $t,s$ namely
\begin{align*}
R^\delta_x(Z,Y)
&= \frac{d}{ds} \dircon(x,z(s))\vert_{s=0}\,\frac{d}{dt} \dircon(x,y(t))\vert_{t=0}
- \frac{d}{dt} \dircon(x,y(t))\vert_{t=0}\, \frac{d}{ds} \dircon(x,z(s))\vert_{s=0} \\ 
& \qquad + \frac{\partial^2\,}{\partial s\,\partial t} \dircon(z(s),y_s(t))\vert_{s=t=0}
- \frac{\partial^2\,}{\partial t\,\partial s} \dircon(y(t),z_t(s))\vert_{t=s=0} \\
&= \frac{\partial^2\,}{\partial s\, \partial t} K_x(z(s),y(t))\vert_{t=s=0} - \frac{\partial^2\,}{\partial s\, \partial t} K_x(y(t),z(s))\vert_{t=s=0} \\
& = \frac{\partial^2\,}{\partial s\, \partial t} \Big(K_x(y(t),z(s))^{-1}\, K_x(z(s),y(t))\Big)\vert_{t=s=0}. 
\end{align*}
The last two equalities are a direct computation from the definition of $K_x(z,y)$ using the Leibniz rule for products and Lemma \ref{lem:d/dt Gamma}. 
\end{proof}


\subsection{Flat groupoids} 
\label{subsec:flat groupoids} 

A direct connection $\dircon:\Pair(M)\locmap \calG$ is flat if and only if it is a local groupoid morphism. Such maps have been extensively studied by K. Mackenzie in \cite[Chapter 6]{Mackenzie:2005}, in particular for what concerns the relationship with their infinitesimal Lie algebroid morphisms. 
We collect here some relevant results, in particular a strong caracterisation of Lie groupoids which admit a flat direct connection.

\begin{defn}
\label{def:flatgroupoid}
A Lie groupoid $\calG\toto M$ is called \textbf{flat} if it is locally trivial and it admits an atlas of local decomposition maps $\tau_\alpha:U_\alpha\to \calG^{x_0}_{U_\alpha}$ (cf. Example \ref{ex:gaugegroupoid}) whose transition functions $g_{\alpha \beta}:U_\alpha\cap U_\beta\to \calG^{x_0}_{x_0}$ given by $g_{\alpha \beta}(x)=\tau_\alpha(x)^{-1}\,\tau_\beta(x)$ are locally constant, i.e. constant on each connected component.
A flat Lie groupoid $\calG$ is called \textbf{locally exact} if furthermore it admits a locally trivial direct connection $\dircon:\Pair(M)\locmap \calG$ (cf. Definition \ref{def:locally trivial dircon}). 
\end{defn}

In essence, the subsequent statement is a reformulation of \cite[Example 6.1.7]{Mackenzie:2005} in terms of direct connections. 

\begin{prop}

\label{prop:flat groupoid} 
A locally trivial Lie groupoid on a connected base admits a flat direct connection if and only if it is flat. In this case, it is locally exact. 
\end{prop}

\begin{proof}
($\Rightarrow$) Assume $\calG\rightrightarrows M$ is a Lie groupoid endowed with a local groupoid morphism $\dircon:\Pair(M) \longlocmor \calG$ defined on a diagonal domain $\calU$, in which case $\calG$ is locally trivial by Proposition \ref{prop:direct}. We choose an open covering $\{U_\alpha\}$ of $M$ such that $\bigcup_\alpha (U_\alpha\times U_\alpha)\subset \calU$ and $x_0\in M$. For every $\alpha$, fix a point $x_\alpha\in U_\alpha$ together with an element $\xi_\alpha\in \mathcal G^{x_0}_{x_\alpha}$ (which exists by the local triviality of $\calG$). Then the collection of maps $\tau_\alpha(x):=\dircon(x,x_\alpha)\,\xi_\alpha$ builds a section atlas with constant transition functions $g_{\alpha \beta}(x)= \xi_\alpha^{-1}\,\dircon(x_\alpha,x_\beta)\,\xi_\beta$.
Moreover, for fixed $\alpha$, the map
$$
\begin{array}{ccl}
\calG_{U_\alpha}^{U_\alpha}& \longrightarrow & {U_\alpha}\times \calG_{x_0}^{x_0}\times {U_\alpha}\\
g &\longmapsto & (\tt(g),\tau_{\alpha}(\tt(g))^{-1}g\tau_{\alpha}(\ss(g)),\ss(g))\\
\end{array}
$$
is an isomorphism under which $\dircon(y,x)$ is sent to $(y,1_{x},x)$ for all $x,y\in U_{\alpha}$.
Hence, $\calG$ is locally exact.
	
($\Leftarrow$) Conversely, let $\tau_\alpha\colon U_\alpha\to \mathcal G^{x_0}_{U_\alpha}$ be a section atlas with constant transition functions $\tau_\alpha(x)^{-1}\, \tau_\beta(x)=g_{\alpha \beta} 1_x$ on the intersections $U_{\alpha \beta}=U_\alpha\cap U_{\beta}$, and set $\calU:=\bigcup_\alpha (U_\alpha\times U_\alpha)$. Then, the map $\dircon: \calU\to \calG$ given by \eqref{eq:exactGamma}, namely $\dircon(y,x)=\tau_\alpha(y)\,\tau_\alpha(x)^{-1}$ for $x,y\in U_\alpha$, is well defined on $\calU$. Indeed, on every intersection $U_{\alpha \beta}$ we have $\tau_\beta(x)=\tau_\alpha(x)\,g_{\alpha\beta}$, therefore 
\begin{align*}
\dircon\vert_{U_\beta\times U_\beta}(y,x) &= \tau_\beta(y)\,\tau_\beta(x)^{-1} \\ 
&= \tau_\alpha(y)\,g_{\alpha \beta}\,g_{\alpha \beta}^{-1}\,\tau_\alpha(x)^{-1}
=  \dircon\vert_{U_\alpha\times U_\alpha}(y,x). 
\end{align*} 
That $\dircon\vert_{U_\alpha\times U_\alpha}$ defines a morphism is easily verified.
\end{proof}

\begin{example}
The trivial Lie groupoid $\calG(M\times G)= M\times G\times M$ is locally exact, because it is clearly flat and any global section $\tau:M\to M\times G$ gives a trivial direct connection $\dircon_\tau(x,y)=(y,\tau(y)\,\tau(x)^{-1},x)$ as in Example \ref{ex:trivial dircon}. 
\end{example}

\begin{remark}\label{rk:fundamentalgroupoid} 
The \textbf{fundamental groupoid} of $M$ is the set $\Pi(M)$ of homotopy classes $[\gamma]$ of continuous paths $\gamma:[0,1]\to M$ with source $s([\gamma])=\gamma(0)$, target $t([\gamma])=\gamma(1)$, partial composition induced by the concatenation of paths, and inversion induced by the inversion of orientation.
This is a gauge groupoid, whose vertex group at a point $x_0$ in $ M$ is the fundamental group $\pi_1(M,x_0)$, see \cite[Examples 1.1.1, 1.3.4]{Mackenzie:2005}. 

By Proposition \ref{prop:flat groupoid}, there exists a flat direct connection $\dircon_0: \Pair(M) \longlocmor \Pi(M)$ which can be constructed as follows. Let $\{U_\alpha\}$ be a cover of $M$ by simply connected open subsets, then $\calU=\bigcup_\alpha (U_\alpha\times U_\alpha)$ defines a diagonal domain. For any $(y,x)\in \calU$, we define $\dircon_0(y,x)$ as the homotopy class of any path in $\calU$ connecting $x$ and $y$. This is well defined as any two paths in $U_\alpha$ are homotopic, and clearly gives a local groupoid morphism. 

In \cite[Prop. 6.1.8]{Mackenzie:2005}, Mackenzie showed an important factorisation result, that we rephrase in terms of direct connections: 
\begin{enumerate}
\item The flat direct connection  $\dircon_0: \Pair(M) \longlocmor \Pi(M)$  is uniquely determined up to germ equivalence (as local maps). 
\item Any flat direct connection $\dircon:\Pair(M)\longlocmor \mathcal{G}$ on a Lie groupoid $\calG$ factorises in a unique way through $\dircon_0$, i.e. there exists a unique groupoid morphism $H:\Pi(M)\to \mathcal G$, such that $\dircon=H \circ \dircon_0$.
\end{enumerate} 
The second statement is the Riemann--Hilbert correspondence. It can be used to prove that a flat infinitesimal connection $\delta:TM\to \Lie(\calG)$ integrates to a flat direct connection $\dircon^\delta$ on $\calG$, which is unique up to germ equivalence, cf. Remark \ref{rem:one-to-one}. 
\end{remark}




\section{Jet prolongation of direct connections}
\label{sec:jetconnections}

Jet prolongations of groupoids were first considered by C. Ehresmann \cite{Ehresmann:1955} and later revisited by I. Kol\'a\v{r} \cite{Kolar:2007}. 
After reviewing the basic concepts regarding jet prolongation of Lie groupoids and Lie algebroids, in this section we study the jet prolongation of direct connections on jet groupoids. 
We start by introducing the \emph{jet primitive bisection} of any direct connection on a jet groupoid (Definition \ref{def:jet primitive bisection}). Using this tool, we show that the jet prolongation of a direct connection is determined by those of the trivial direct connection $\Id$ on the pair groupoid, and describe their associated infinitesimal connection (Theorem \ref{thm:jet prolongation dircon}). 
We then exhibit three jet prolongations of $\Id$ on $\Pair(M)$, based on three Jacobi fields of variation curves through geodesics which differ by the initial value of the field and of its covariant derivative (Theorems \ref{thm:parallel direct connection}, \ref{thm:Jacobi direct connection} and \ref{thm:exponential direct connection}). 
In the three cases we also compute the associated infinitesimal connection and show that the first two are, in a sense, orthogonal to each other, while the third one reproduces the exponential jet connection on the Lie algebroid $J^nTM$ given by I. Kol\'a\v{r} in \cite{Kolar:2007} and by W. M. Mikulski in \cite{Mikulski:2007}. 



\subsection{Jet prolongation of groupoids}
\label{subsec:jet groupoid} 

Following \cite{Ehresmann:1955}, \cite[\S 1]{Kolar:2003}, \cite[\S 4.5]{Meinrenken:2017}, we recall the definition and main properties of jet groupoids, using local bisections of Section \ref{subsec:bisections}. 
The \textbf{$n$-jet prolongation of $\calG$} is the set of $n$-jets of local bisections $\sigma:M\locmap \calG$, that is, the jet space
\begin{align} \label{eq:J^n calG}
J^n\calG &:= J^n\Bis_\loc(\calG) = \{ j^n\sigma\in J^n(M,\calG)\ |\ \ss\circ\sigma=\Id,\ \mbox{$\varphi_\sigma=\tt\circ\sigma$ is a diffeomorphism}\},
\end{align}
together with the structure of a Lie groupoid on $M$ induced by that of $\calG$: 
\begin{enumerate}[label={(\roman*)}]
\item \textbf{source} and \textbf{target} maps $\ss^n,\tt^n:J^n\calG\to M$ given respectively by the surjective submersions $\ss^n(j^n_x\sigma)=x$ and $\tt^n(j^n_x\sigma)=\tt\big(\sigma(x)\big) =\varphi_\sigma(x)$ for any $x\in M$,  

\item the \textbf{multiplication} $j^n_{x'}\sigma' j^n_x\sigma = j^n_x\big(\sigma'\bullet \sigma\big)$ defined if $x'= \varphi_\sigma(x)$, where $\bullet$ is the product of bisections given in Section \ref{subsec:bisections}, 

\item \textbf{unit} $\uu^n(x)\equiv 1^n_x=j^n_x \uu$, where $\uu:M\to \calG$ is the unit bisection of $\calG$,

\item \textbf{inverse} $(j^n_x\sigma)^{-1}= j^n_{\varphi_\sigma^{-1}(x)}\sigma^{-1}$ where $\sigma^{-1}$ is the inverse local bisection as given above. 
\end{enumerate}
From (i) it follows that 
\begin{align} \label{eq:(J^nG)^x_y}
(J^n\calG)^x_y=\{j^n_x\sigma\ |\ \text{$\sigma$ is a local bisection on $\calG$ defined around $x$ such that $\varphi_\sigma(x)=y$}\}.
\end{align} 
Since $\calG$ is fibered over $M$ by the source map, the jet prolongation $J^n\calG$ enjoys the same properties as the jet bundle of a fibre bundle described in Section \ref{subsec:jet bundle}. In particular:

\begin{itemize}
\item $J^n\calG$ is a closed submanifold of the jet space $J^n(M,\calG)$ \cite[\S 12.16]{Kolar-Michor-Slovak:1993},


\item the jet projections $\pi^k_{k-1}:J^k\calG\to J^{k-1}\calG,\ j^k_x\sigma \mapsto j^{k-1}_x\sigma$, for any $1\leq k\leq n$, give a sequence of Lie groupoid morphisms \cite[Example 1.12]{Meinrenken:2017}
\begin{align} \label{eq:jet groupoid tower}
J^n\calG\rightarrow J^{n-1}\calG\rightarrow \cdots \rightarrow J^1\calG\rightarrow J^0\calG=\calG.
\end{align}
\end{itemize}
In particular, the natural projection $\pi^n_0:J^n\calG\to \calG$ is then also a groupoid morphism. In fact, for any $\sigma,\sigma'$ in $  \Bis_{\rm loc}(\calG)$ and for any $x$ in $ M$ such that $\sigma'$ is defined in $x'=\varphi_\sigma(x)$, we have
\begin{align*}
\pi_0^n(j^n_{x'}\sigma'\,j^n_x\sigma) &= \pi_0^n(j^n_x(\sigma'\bullet\sigma)) = (\sigma'\bullet\sigma)(x) \\ 
&= \sigma'(x')\,\sigma(x) = \pi_0^n(j^n_{x'}\sigma')\,\pi_0^n(j^n_x\sigma). 
\end{align*}

A smooth map $\phi: \calG\longrightarrow \calG'$ between two Lie groupoids over $M$ can be prolonged to jets using the jet prolongation of the pushforward $\Bis_{\rm loc}(\phi)$ in \eqref{eq:MapsBisections-local} on local bisections, namely
\begin{align}\label{eq:J^nPhi on calG}
J^n\phi: J^n\calG\longrightarrow J^n\calG',\ j^n_x\sigma \mapsto J^n\phi(j^n_x\sigma) := j^n_x(\phi\circ\sigma). 
\end{align}
If $\phi$ is a groupoid morphism, then $J^n\phi$ is also a groupoid morphism, because $\Bis_{\rm loc}(\phi)$ is a morphism of pseudo-groups (cf. Section \ref{subsec:bisections}).
Thanks to the locality of the notion of jets, the jet prolongation \eqref{eq:J^nPhi on calG} actually extends to any local map $\phi:\calG\locmap \calG'$. 
Finally, the \textbf{$n$-jet prolongation} $J^n$ is a functor on Lie groupoids with smooth local maps, which preserves morphisms. 

\begin{prop}\label{prop:J^nPhi}
Any smooth local map $\phi:\calG\locmap \calG'$ yields a local jet prolongation $J^n\phi: J^n\calG\locmap J^n\calG'$.
If $\phi$ is a groupoid morphism, then $J^n\phi$ is also a groupoid morphism.
\end{prop}

\begin{proof}
Suppose that $\phi$ is defined on an open neighbourhood $\calU\subset \calG$ of the units $\uu(M)\subset \calG$. We only need to define an open neighbourhood $\calU^n\subset J^n\calG$ of the units $\uu^n(M)\subset J^n\calG$ on which the jet prolongation \eqref{eq:J^nPhi on calG} is well defined. Note that, being a local map in the sense of Definition \ref{def:localmap}, $\phi$ preserves the source, the target and the units. Then, since $\calU$ is open in $\calG$, the subset 
\begin{align*}  
\Bis_\loc(\calU) = \{\sigma\in \Bis_\loc(\calG)\ |\ \mbox{$\sigma$ is defined on $U\subset M$ such that $\sigma(U)\subset \calU$}\}
\end{align*}
is open in $\Bis_\loc(\calG)$, by the topological properties in \cite[\S 4.4]{Michor:1980}. Therefore, the set
\begin{align*}  
\calU^n := J^n\Bis_\loc(\calU) = \{j^n_x\sigma\in J^n\calG\ |\ \mbox{$\sigma$ is defined on $U\subset M$ such that $\sigma(U)\subset \calU$}\}
\end{align*}
satisfies the required properties, and moreover the jet projection $\pi^n_0$ restricted to $\calU^n$ has image in $\calU$. Indeed, $\calU^n$ contains the units $\uu^n(M)=\{j^n_x\uu\ |\ x\in M\}$, because $\uu(M)\subset \calU$, and it is open in $J^n\calG$ in the Fr\'echet topology induced by the supremum norm in all derivatives of order no larger than $n$ on compact subsets, which corresponds to compact open topology in \cite[\S 4.3]{Michor:1980}. 
Finally, the fact that $J^n$ preserves groupoid morphisms can be easily verified with a direct computation. 
\end{proof}


\begin{example}[Jet prolongation of gauge groupoids]\label{ex:J^nPair(M)} 
It was proved in \cite[\S 4]{Kolar:2007}  that the $n$-jet prolongation of a gauge groupoid $\calG(P) \rightrightarrows M$ is isomorphic to the gauge groupoid of the $n$-jet principal bundle $W^nP$, namely 
\begin{align}\label{eq:J^nG(P)=G(W^nP)}
J^n\calG(P) \cong \calG(W^nP).  
\end{align}
Furthermore, the jet prolongation of a morphism $\phi:\calG(P)\to \calG(P')$ yields a morphism of gauge groupoids $J^n\phi:\calG(W^n P)\to\calG(W^n P')$.
This gives rise to several examples:
\begin{enumerate}[label={(\alph*)}]
\item For $P=M\times G$ we have $W^n(M\times G) \cong F^nM \times T^n_dG$ (see Lemma \ref{lem:W^n(MxG)} in the Appendix, eq. \eqref{eq:W^n(MxG) formula 2}), hence
\begin{align}\label{eq:JnGMtimesG}
J^n(M\times G\times M)\cong (F^nM\times T^n_d G\times F^nM)/GL^n_d,  
\end{align}
where $j^n_0 h\in GL^n_d$ acts on $(j^n_0\varphi',j^n_0 g,j^n_0 \varphi) \in F^nM\times T^n_d G\times F^nM$ (with $\varphi(0)=x$ and $\varphi'(0)=y$) by simple composition $(j^n_0(\varphi'\circ h),j^n_0 (g\circ h),j^n_0 (\varphi\circ h))$. 

\item
In particular, for $G=\{1\}$ we have $T^n_d(\{1\})=\{1\}$ and $W^n(M\times \{1\})=F^nM$, therefore
\begin{align}\label{eq:J^nPair(M)} 
J^n\Pair(M) = \mathrm{inv}J^n(M,M) \cong \calG(F^nM) = F^nM \times_{\GLnd} F^nM.  
\end{align}
For $n=1$, since $GL^1_d= GL_d(\R)$ and $F^1M=F(TM)$ is a frame bundle, the $1$-jet groupoid is a frame groupoid \cite[Example 1.12]{Meinrenken:2017} 
\begin{align}\label{eq:J^1Pair(M)} 
J^1\Pair(M) \cong \calG(FM) \cong \GL(TM). 
\end{align}

\item
The frame groupoid $\GL(E)$ of a vector bundle $E$ with rank $r$ is the gauge groupoid of the frame bundle $FE$ with fibre $GL_r:=GL_r(\R)$. It follows that
\begin{align}\label{eq:J^nIso(E)} 
J^n\GL(E) \cong \calG(W^nFE) = W^nFE \times_{W^n_d GL_r} W^nFE.  
\end{align}

\item Suppose that a frame bundle $FE$ admits a reduction first to a principal $G$-bundle $P$, and further to the trivial principal bundle $M\times \{1\}$ (which forces $E$ to be trivial, i.e. $E\cong M\times \R^r$). Then we have a sequence of gauge subgroupoids
\begin{align*}
\calG(M\times \{1\})=\Pair(M) \subset \calG(P)\subset \calG(FE)=\GL(E). 
\end{align*}
The jet bundle $J^nE$ has structure group $GL(P^n_{d,r})$, frame bundle $F(J^nE)$ and frame groupoid $\GL(J^nE)$, which all admit reductions to the jet prolongation of the starting objets: 
we have a sequence of subgroups
\begin{align} \label{eq:sequence jet groups}
W^n_d\{1\} = \GLnd
& \subset W^n_dG = \GLnd\ltimes T^n_dG \\ \nonumber
& \hspace{1cm} \subset W^n_dGL_r = \GLnd\ltimes T^n_dGL_r \subset GL(P^n_{d,r}), 
\end{align}
a corresponding sequence of reduced principal bundles 
\begin{align} \label{eq:sequence jet principal bundles}
W^n(M\times \{1\}) = F^nM
& \subset W^nP= F^nM \times_M J^nP \\ \nonumber 
& \hspace{1cm} \subset W^nFE = F^nM \times_M J^nFE \subset F(J^nE),
\end{align}
and an induced sequence of subgroupoids
\begin{align} \label{eq:sequence jet groupoids}
\calG(F^nM) = J^n\Pair(M)
& \subset \calG(W^nP)=J^n\calG {(P)} \\ \nonumber
& \hspace{1cm} \subset \calG(W^nFE)= J^n\GL(E) \subset \calG(F(J^nE))=\GL(J^nE). 
\end{align}
\end{enumerate}
\end{example}


\begin{example}[Frame groupoid of a jet bundle]\label{ex:F(J^nE)}
Let $E\to M$ be a vector bundle.
The frame groupoid $\GL(J^nE)$ of the jet bundle of $E$ is the gauge groupoid of the frame bundle $F(J^nE)$ and it is \emph{not} a jet prolongation, since it contains $J^n\GL(E)=\calG(W^nFE)$ as a proper subgroupoid.
Indeed, by Lemma \ref{lem:W^n_dG subgroup of GL(P^n_{d,r})}, $GL(P^n_{d,r})$ contains $W^n_dGL_r$ as a proper subgroup, hence the frame bundle $F(J^nE)$ contains the jet prolongation $W^nFE$ as a proper principal subbundle. Then, by \eqref{eq:F(J^nE) vs W^nFE} and \eqref{eq:J^nG(P)=G(W^nP)} we have
\begin{align} \nonumber
  \GL(J^nE) & \cong \calG(F(J^nE))
  = F(J^nE) \times_{GL(P^n_{d,r})} F(J^nE)
 \label{eq:Iso(J^nE)} \\  
  & \cong W^nFE \times_{W^n_dGL_r} F(J^nE) 
  \\   \nonumber
& \cong W^nFE \times_{W^n_dGL_r} GL(P^n_{d,r})\times_{W^n_d GL_r} W^nFE,   
\end{align}
while by \eqref{eq:J^nIso(E)} we have 
$J^n\GL(E) \cong W^nFE \times_{W^n_dGL_r} W^nFE$. 
The inclusion $J^n\GL(E)\hookrightarrow \GL(J^nE)$ is given by $[A,B]_\sim \mapsto [A,1^n,B]_\sim$ for any $A,B\in W^nFE$, where $1^n$ is the unit in the group $GL(P^n_{d,r})$, with equivalence relations $(AH,B)\sim (A,B\,H^{-1})$ and $(AH,G,BH') \sim (A,HG(H')^{-1},B)$ for any $G\in GL(P^n_{d,r})$ and $H,H'\in W^n_d GL_r$.
\end{example}

\begin{example}[Action of a jet frame groupoid on the jet bundle]
\label{ex:Action of J^nIso(E) on J^nE}
The inclusion $J^n\GL(E)\hookrightarrow \GL(J^nE)$ of Example \ref{ex:F(J^nE)} can also be understood as the canonical action of a jet $j^n_{x_0}\sigma\in J^n_{x_0}\GL(E)$ on the bundle $J^nE$, where $\sigma:U\to \GL(E)$ is a local bisection defined on an open neighbourhood of $x_0\in M$, that is, as a linear map $j^n_{x_0}\sigma:J^n_{x_0}E\to J^n_{y_0}E$, where $y_0=\varphi_\sigma(x_0)$. 
For any $j^n_{x_0}f\in J^n_{x_0}E$, where $f:U\to E_U$ is a local section around $x_0$, this action is given by  
\begin{align} \label{eq:action of j^n sigma on j^n f} 
  j^n_{x_0}\sigma \cdot j^n_{x_0}f
  & := j^n_{y_0}\Big(y\mapsto \sigma\big(\varphi_\sigma^{-1}(y)\big)\cdot f\big(\varphi_\sigma^{-1}(y)\big)\Big)
  \\ \nonumber
  & = j^n_{x_0}\Big(x\mapsto \sigma(x)\cdot f(x)\Big) \circ j^n_{y_0}\varphi_\sigma^{-1}, 
\end{align}
where $\varphi_\sigma^{-1}$ is the inverse diffeomorphism of $\varphi_\sigma$ and where we denote by $\cdot$ the linear action of $\sigma(x)\in \GL(E)^x_{\varphi_\sigma(x)}$ on $f(x)\in E_x$, as in \eqref{eq:actionIso(E)}, and by $\circ$ the jet composition as in \eqref{eq:jet composition}. 
Indeed, the expression $\tilde{f}(y)= \sigma\big(\varphi_\sigma^{-1}(y)\big)\cdot f\big(\varphi_\sigma^{-1}(y)\big)$ of eq. \eqref{eq:action of j^n sigma on j^n f} gives the section of $E$ on $\varphi_\sigma(U)$ making the following diagram commutative: 
$$
\begin{tikzcd}
E_U \arrow[r, "\sigma\cdot"] & E_{\varphi_\sigma(U)} \\
U \arrow[u, "f"] \arrow[r, "\varphi_\sigma"'] & \varphi_\sigma(U) \arrow[u,"\tilde{f}"'] 
\end{tikzcd}
$$
\end{example}


\subsection{Jet prolongation of Lie algebroids} \label{subsec:jet Lie algebroids}

The $n$-jet bundle of a Lie algebroid $A\to TM$ is the vector bundle $J^nA\to M$ of its underlying vector bundle $A\to M$. It turns out \cite[Proposition 1]{Kumpera:2015} \cite[Example 6.11]{Meinrenken:2017} that this jet bundle is again a Lie algebroid $J^nA\to TM$, with anchor and bracket given by
\begin{align} \label{eq:jet-bracket}
a^n(j^n_xX) =a(X_x) \qquad\mbox{and}\qquad [j^n_xX,j^n_xY]_{J^nA} = j^n_x([X,Y]_A)
\end{align}
for any smooth local sections $X,Y$ of $A$ defined around $x$ in $M$.
If $\phi:A\to A'$ is a smooth map between Lie algebroids which commutes with the anchors, then there exists a jet prolongation
\begin{align} \label{eq:J^nPhi-algebroid}
J^n\phi:J^nA\to J^nA',\ j^n_xX\mapsto J^n\phi(j^n_xX):=j^n_x(\phi\circ X) 
\end{align}
which commutes with the jet prolongation of the anchors. 
If $\phi$ also preserves the bracket, hence it is a morphism of Lie algebroids, then by \eqref{eq:jet-bracket} the map $J^n\phi$ is also a morphism of Lie algebroids. 

\begin{example}[Jet prolongation of Atiyah algebroids]\label{ex:J^nA(P)}
For a principal $G$-bundle $P$, the $n$-jet prolongation of its Atiyah algebroid $A(P)$ is isomorphic to the Atiyah algebroid of its $n$-jet principal bundle $W^nP$ \cite[\S 2 Eq. (8)]{Kolar:2008}
\begin{align} \label{eq:J^nA(P)=A(W^nP)}
J^n A(P) \cong A(W^nP). 
\end{align}
A morphism $\phi:P\to P'$ of principal bundles over $M$ induces a morphism $A(\phi):A(P)\to A(P')$ between their algebroids \cite[Definition 3.2.8]{Mackenzie:2005} and a principal prolongation \eqref{eq:W^nPhi on P} $W^n\phi: W^nP\to W^nP'$. One can check that the jet prolongation \eqref{eq:J^nPhi-algebroid} $J^n A(\phi): J^nA(P)\to A(P')$ coincides with the algebroid map $A(W^n \phi): A(W^nP)\to A(W^nP')$.
\end{example}

\begin{example}[Lie algebroid of jet groupoids]\label{ex:J^nLgauge}
If $\calG$ is a Lie groupoid with Lie algebroid $A(\calG)$, there is an isomorphism of Lie algebroids \cite[Proposition 1 and Theorem 1]{Kumpera:2015}, \cite[Example 9.5 (e)]{Meinrenken:2017}
\begin{align} \label{eq:J^nL(G)=A(J^nG)}
A\big(J^n\calG\big) \cong J^n A(\calG).
\end{align}
A smooth local map $\phi:\calG\locmap \calG'$ of Lie groupoids over $M$ induces a vector bundle map $d^{\tt}\phi\vert_M:A(\calG)\to A(\calG')$ between their Lie algebroids by Proposition \ref{prop:t-differential}, and a local jet prolongation $J^n \phi:J^n\calG\locmap J^n\calG'$ by Proposition \ref{prop:J^nPhi}. One can check that the jet prolongation $J^n(d^{\tt}\phi\vert_M):J^nA(\calG)\to J^nA(\calG')$ coincides with the $\tt$-differential $d^{\tt}(J^n \phi)\vert_M:A(J^n\calG)\to A(J^n\calG')$. 
\end{example}

\begin{examples}[Summary on jet algebroids and jet groupoids]
The previous results relate in many ways Atiyah algebroids, gauge groupoids and their jet prolongations:
\begin{align*}
  J^nA(\calG(P)) \overset{\eqref{eq:J^nL(G)=A(J^nG)}}{\cong}
  A(J^n\calG(P)) \overset{\eqref{eq:J^nG(P)=G(W^nP)}}{\cong}
  A(\calG(W^nP)) \overset{\eqref{eq:L(G(P))=A(P)}}{\cong}
  A(W^nP) \overset{\eqref{eq:J^nA(P)=A(W^nP)}}{\cong}
  J^nA(P).
\end{align*}
These isomorphisms give some key examples: 
\begin{enumerate}[label={(\alph*)}]
\item
\cite[Proposition 1]{Kolar:2009}, \cite[\S 2 Eq. (6)]{Kolar:2008}
For $P=M\times \{1\}$, we have $A(M\times \{1\})=TM$, $W^n(M\times \{1\})=F^nM$ and $\calG(M\times \{1\})=\Pair(M)$, therefore  
\begin{align}\label{eq:J^nTM}
A(J^n\Pair(M))\cong A(F^nM)\cong J^nTM.
\end{align}
The case $n=1$ is particularly important for connections, and gives a sequence of key isomorphisms from $J^1TM$ to the bundle $\Tlin TM$ of Example \ref{exs:Lie algebroids} (e): 
\begin{align}\label{eq:J^1TM}
J^1TM \overset{\eqref{eq:J^nTM}}{\cong} A(FM) \overset{\eqref{eq:L(G(P))=A(P)}}{\cong} A(J^1\Pair(M)) \overset{\eqref{eq:J^1Pair(M)}}{\cong} A(\GL(TM)) = \Tlin TM. 
\end{align}
Using the isomorphism of Example \ref{ex:Der(E)}, we also have $J^1TM \cong \Der(TM)$.

\item
For $P=M\times G$, we have $A(M\times G)=TM\oplus (M\times \mathfrak{g})$, $W^n(M\times G)=F^nM\times T^n_dG$ and $\calG(M\times G)=M\times G\times M$, therefore
\begin{align*}
J^n\big(TM\oplus (M\times \mathfrak{g})\big) \cong A\big(J^n(M\times G\times M)\big)\cong A(F^nM\times T^n_dG). 
\end{align*}
Since $F^nM$ is a principal $GL^n_d$-bundle, this algebroid fits in Example \ref{exs:Lie algebroids} (d), and we have
\begin{align}\label{eq:JnTMplusg}
J^n\big(TM\oplus (M\times \mathfrak{g})\big) \cong A(F^nM) \oplus (M\times \mathrm{Lie}(T^n_dG)) \cong J^nTM \oplus (M\times T^n_d\mathfrak{g}).
\end{align}

\item
If $P=FE$, then $A(FE)=\Tlin E$ is the Lie algebroid of linear vector fields of $E$, which is isomorphic to that of derivative endomorphisms $\Der(E)$, $W^n(FE)\cong F^nM\overset{\pi}{\times} J^nFE$ and $\calG(FE)=\GL(E)$. Therefore
$$
J^n \Tlin E \cong A(J^n\GL(E))\cong A\big(F^nM\times_M J^nFE\big) \cong J^n\Der(E).
$$
\end{enumerate}
\end{examples}

 
\subsection{Jet prolongation of direct connections}
\label{subsec:jet prolongation of direct connections}

For usual connections on jet prolongations of fibre bundles, the notion of \emph{$n$-jet prolongation of a connection} in given in \cite[Section 45.4]{Kolar-Michor-Slovak:1993} and depends on the choice of a linear connection on the jet bundle $J^nTM$. 
In this section we show that an analogue definition can be given for direct connections on groupoids. 

\begin{defn}
Let us call \textbf{$n$-th order direct connection on a Lie groupoid $\calG$} any direct connection $\Sigma:\Pair(M)\locmap J^n\calG$. 
Then $\Sigma$ determines: 
\begin{itemize}
\item 
a direct connection $\Sigma^0$ on $\calG$, called its \textbf{$0$-th order projection},  given by the composite map
\begin{align} \label{eq:projection of n-direct connection}
\Sigma^0 = \pi^n_0\circ \Sigma:\Pair(M)\locmap \calG ,\ 
\end{align}
where $\pi^n_0:J^n\calG\to J^0\calG=\calG$ is the jet projection described in \S \ref{subsec:jet groupoid}, 

\item 
an $n$-th order direct connection $\Delta_\Sigma$ on $\Pair(M)$, called its \textbf{projection to $J^n\Pair(M)$}, given by the composite map 
\begin{align} \label{eq:D direct connection on J^nPair(M)}
\Delta_\Sigma = J^n \aa \circ \Sigma:\Pair(M)\locmap J^n\Pair(M),\ 
\end{align}
where $\aa\colon\calG\to \Pair(M)$ is the anchor of $\calG$ and $J^n \aa\colon J^n\calG\to J^n\Pair(M)$ is its transform under the jet functor $J^n$ defined in \eqref{eq:J^nPhi on calG}. 
Of course its $0$-th order projection on $\Pair(M)$ is always $\Delta^0_\Sigma:=\pi^n_0 \circ \Delta_\Sigma = \Id$. 
\end{itemize}
\end{defn}

These two direct connections are related by the equality 
\begin{align} \label{eq: a Sigma_0 = a^n Delta_Sigma}
    \aa\circ \Sigma^0  = \aa_M^n\circ \Delta_\Sigma = \aa^n \circ \Sigma = \Id_{\Pair(M)}, 
\end{align}
where $\aa_M^n:J^n\Pair(M)\to \Pair(M)$ is the anchor of $J^n\Pair(M)$, as defined in section \ref{subsec:jet groupoid}. In other words, we have a commutative diagram
$$
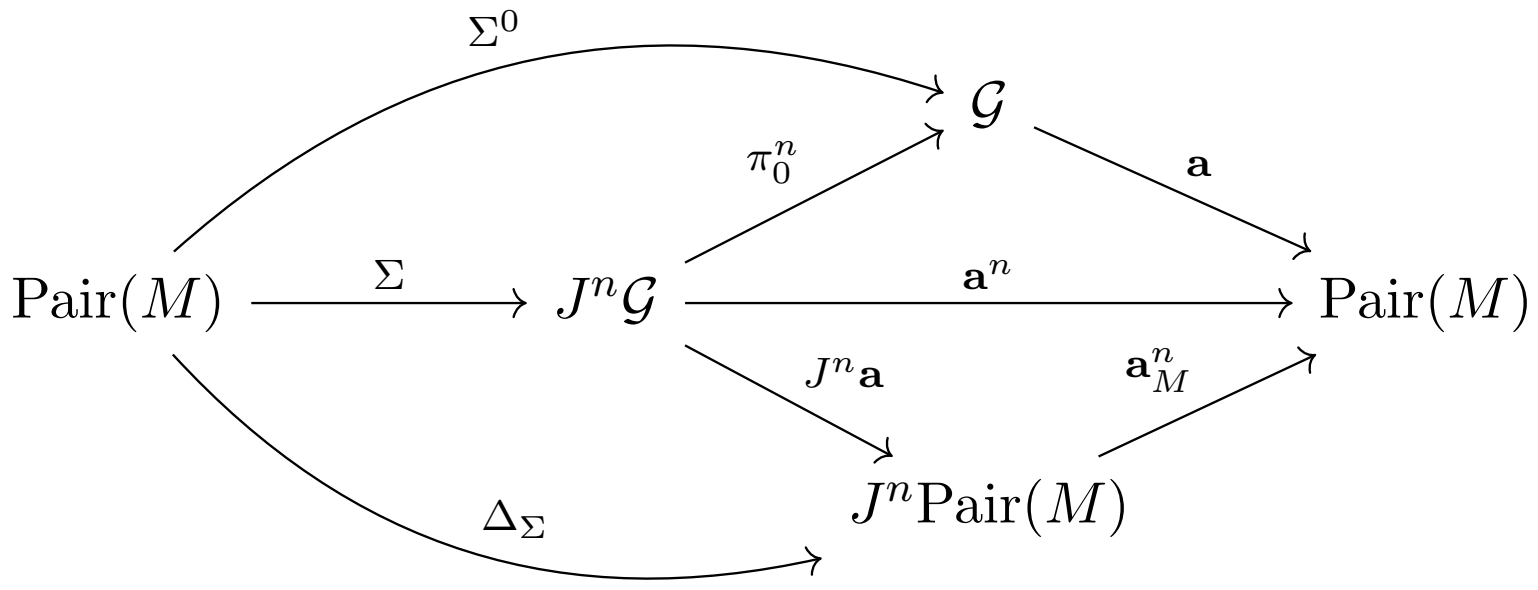
\begin{tikzcd}
&& & \calG \arrow[dr, "\aa"] & \\ 
\Pair(M) \arrow[rr, "\Sigma"] \arrow[urrr, bend left, "\Sigma^0"] \arrow[drrr, bend right, "\Delta_\Sigma"] && J^n\calG \arrow[rr, "\aa^n"] \arrow[ur, "\pi^n_0"] \arrow[dr, "J^n \aa"] & & \Pair(M) \\ 
&& & J^n\Pair(M) \arrow[ur, "\aa^n_M"] & 
\end{tikzcd}
$$

\begin{defn} \label{def:n-jet prolongation of direct connection}
The direct connection $\Sigma$ is called an \textbf{$n$-jet prolongation} if 
\begin{align*}
    \Sigma = J^n \Sigma^0\circ \Delta_\Sigma,
\end{align*}
where $J^n\Sigma^0:J^n\Pair(M)\to J^n\calG$ is the transformed of $\Sigma^0$ under the jet functor $J^n$ given in \eqref{eq:J^nPhi on calG}. 
In other words, in the spirit of \cite[Section 45.4]{Kolar-Michor-Slovak:1993}, given a direct connection $\dircon:\Pair(M)\locmap \calG$ and an $n$-th order direct connection $\Delta:\Pair(M)\locmap J^n\Pair(M)$, we call \textbf{$n$-jet prolongation of $\dircon$ along $\Delta$} the direct connection on $J^n\calG$ given by the composite local map 
\begin{align} \label{eq:n-jet prolongation of direct connection}
\dircon^n = J^n\dircon \circ \Delta. 
\end{align}
We then have $(\dircon^n)^0=\dircon$ and
\begin{align*}
    \Delta_{\dircon^n} &= J^n \aa \circ \dircon^n 
    = J^n \aa \circ J^n\dircon \circ \Delta
    = J^n (\aa \circ \dircon) \circ \Delta
    = \Delta.
\end{align*}
Moreover, the jet prolongation $\dircon^n$ is compatible with the filtration on $J^n \calG$, that is, 
\begin{align}\label{eq:filtrGamman}
  \pi^n_{k} \dircon^n = \dircon^k ,
\end{align}
where $\pi^n_k: J^n \calG\to J^k \calG$ are the canonical jet projections for any $k=0,...,n-1$.
\end{defn}

\begin{remark}\label{rem:holonomic bisection}
The above definition of jet prolongation for direct connections is \emph{relative} to a given higher-order connection $\Delta$ on $\Pair(M)$ with no restrictions. One might ask if there are meaningful restrictions on $\Delta$ which  guarantee nice properties of the jet prolongation of direct connections. 

One nice property is to \emph{produce holonomic bisections} of $J^n\Pair(M)$. 
According to Remark \ref{rem:dircon vs bisection}, if $\Delta$ is defined on $\calU$ then there exists a global bisection $\widehat{\Delta}_\phi:M \to J^n\Pair(M)$ for any diffeomorphism $\phi\in \Diff_\calU(M)$. 
Requiring $\widehat{\Delta}_\phi$ to be holonomic for every $\phi$ makes no sense. 
However the bisection is holonomic for $\phi = \id$, since in this case $\widehat{\Delta}_\id(x)=\Delta(x,x)=1^n_x=j^n_x \id$. 

We can then require that there exists a non-trivial diffeomorphism $\phi$ for which $\widehat{\Delta}_\phi$ is holonomic, that is,
\begin{align*}
    \widehat{\Delta}_\phi(x) &= \Delta(\phi(x),x) = j^n_x\phi.
\end{align*}
Such a diffeomorphism does not necessarily exist, but is there is one, it induces a holonomic bisection of $J^n\calG$ via
\begin{align*}
    \widehat{\Sigma}_\phi(x) &= \Sigma(\phi(x),x) &= j^n_x\big(z\mapsto \dircon(\phi(z),z)\big).
\end{align*}
\end{remark}

Hence jet prolongations of direct connections on a groupoid $\calG$ are completely determined by $n$-th order direct connections $\Delta$ on $\Pair(M)$. 
Let us describe them explicitly. 

\begin{defn} \label{def:jet primitive bisection}
Let $\Sigma$ be an $n$-th order direct connection on $\calG$ defined on a diagonal domain $\calU$. 
We say that $\Sigma$ \textbf{admits a jet primitive bisection} if there is a smooth collection of local bisections $\{\sigma^x_y\}$ of $\calG$, globally given by a smooth map 
\begin{align*}
\sigma_\Sigma: \tilde{\calU} \subset \calU \times M \to \calG,\ (y,x,z)\mapsto \sigma^x_y(z),
\end{align*}
where $\tilde{\calU}$ is an open subset containing $\{(y,x,x)\ |\ (y,x)\in \calU\}$, such that, for every $(y,x)\in \calU$, 
\begin{enumerate}[label={(\roman*)}]
\item the map $\sigma^x_y:\tilde{\calU}_{y,x} = \{z\in M\ |\ (y,x,z)\in \tilde{\calU} \} \to \calG$ is a local bisection, 
\item $\Sigma(y,x)=j^n_x \sigma^x_y$,  
\item the local diffeomorphisms $\varphi^x_y=\tt\circ \sigma^x_y$ satisfy $\varphi^x_y(x)=y$ and $\varphi^x_x=\id$. 
\end{enumerate}
\end{defn}

\begin{lem} \label{lem:jet prolongation Delta^n}
Every smooth direct connection $\Delta:\Pair(M)\locmap J^n\Pair(M)$ defined on $\calU$ admits a smooth jet primitive bisection 
\begin{align*}
\hat{\varphi}: \tilde{\calU}\to \Pair(M),\ (y,x,z) \mapsto \widehat{\varphi^x_y}(z) = (\varphi^x_y(z),z)
\end{align*}
with $\varphi^x_y(x)=y$ and $\varphi^x_x = \id$. 
\end{lem}

\begin{proof}
For any fixed $(y,x)\in \calU$, one can choose an arbitrary representative bisection $\widehat{f}$ for the jet element $\Delta(y,x)\in J^n\Pair(M)^x_y$, of the form $\widehat{f}(z)=(f(z),z)$ for a local diffeomorphism $f$ defined in an open neighbourhood of $x$, which necessarily satisfies the properties
\begin{align*}
   \Delta(y,x)=j^n_x\widehat{f} \qquad\text{and}\qquad f(x)=y.
\end{align*}
Among such local bisections, we need to show that there always exists one which is smooth in the variables $(y,x)$. 
The proof is in three steps: 
\begin{itemize}
\item First consider the local diffeomorphism  $F_{y,x}= \exp_y^{-1}\circ f\circ \exp_x:T_xM \longrightarrow T_yM$, such that 
$$
f=\exp_y \circ\, F \circ \exp_x^{-1}.
$$
Note that $F_{y,x}\neq df_x$, the two maps coincide only if $f$ is an affine map, which would imply that $M$ is flat.

\item Then compute the Taylor polynomial of order $n$ at $0_x$ of the (real) map $F_{y,x}$, and read it as a polynomial map $\phi_{y,x}:T_xM\to T_yM$ with the same jet as $F_{y,x}$, that is, 
$$
j^n_{0_x} \phi_{y,x} = j^n_{0_x} F_{y,x}, 
$$
and in particular an invertible linear part. By the inverse function theorem $\phi_{y,x}$ is invertible. 
\item Finally define a new local bisection $\widehat{\varphi^x_y}$ of $\Pair(M)$ associated to the local diffeomorphism 
$$
\varphi^x_y = \exp_y \circ\, \phi_{y,x} \circ \exp_x^{-1}.  
$$
\end{itemize}
Now, the diffeomorphism $\varphi^x_y$ satisfies the wished properties, namely
\begin{align*}
    \varphi^x_y(x) &= \exp_y \big(\phi_{y,x}(0)\big) = y
\end{align*}
and, using the properties of jet composition, 
\begin{align*}
    j^n_x \varphi^x_y &= (j^n_{0_y} \exp_y) \circ (j^n_{0_x} \phi_{y,x}) \circ (j^n_x \exp_x^{-1}) \\ 
    & = (j^n_{0_y} \exp_y) \circ (j^n_{0_x} F_{y,x}) \circ (j^n_x \exp_x^{-1}) \\ 
    & = j^n_x\big( \exp_y \circ F \circ \exp_x^{-1} \big)
    = j^n_x f. 
\end{align*}
Finally, $\varphi^x_y$ is smooth in $(y,x)$, because the exponential maps are smooth and the polynomial $\phi_{y,x}$ is also smooth, independently of the original function $F_{y,x}$ (and therefore $f$), of which it only retains the value of derivatives at $0_x$.  
\end{proof}

\begin{thm} \label{thm:jet prolongation dircon}
A direct connection $\Sigma$ on $J^n\calG$ is the jet prolongation of a direct connection $\dircon$ on $\calG$ if and only if it admits a smooth jet primitive bisection $\sigma_\Sigma=\{\sigma^x_y\}$ satisfying 
\begin{align} \label{eq:jet condition on primitive bisection}
j^n_x \sigma^x_y = j^n_x(z\mapsto \dircon(\varphi^x_y(z),z)), 
\end{align}
where $\{\varphi^x_y\}$ is a jet primitive bisection of $\Delta_\Sigma$. 
In particular, one can then choose $\sigma^x_y(z) = \dircon(\varphi^x_y(z),z)$.

In this case, the infinitesimal connection $\delta^\Sigma:TM \to J^n A(\calG)$ of $\Sigma$ is given on any $Y\in T_xM$ by
\begin{align*}
\delta^\Sigma(Y) &= j^n_x\left(z\mapsto \delta^\dircon\left(\frac{d}{dt}\varphi^x_{y(t)}(z)\big|_{t=0}\right)\right) 
= J^n \delta^\dircon \circ \delta^{\Delta_\Sigma}(Y),
\end{align*}
where $\delta^\dircon$ is the infinitesimal connection of $\dircon$, $y(t)$ is a small geodesic with $y(0)=x$ and $\dot{y}(0)=Y$, and 
\begin{align} \label{eq:easy delta^n}
\delta^{\Delta_\Sigma}(Y) &= j^n_x\left(z\mapsto \frac{d}{dt}\varphi^x_{y(t)}(z)\big|_{t=0}\right)
\end{align}
is the infinitesimal connection $\delta^{\Delta_\Sigma}:TM\to J^nTM$ of $\Delta_\Sigma$. 
\end{thm}

\begin{proof}
Suppose that $\Sigma = J^n\dircon \circ \Delta$ is a jet prolongation. 
Apply Lemma \ref{lem:jet prolongation Delta^n} to $\Delta$, choose a family of local diffeomorphisms $\{\varphi^x_y\}$ such that $\Delta(y,x)=j^n_x\widehat{\varphi^x_y}$ and define 
\begin{align*}
    \sigma^x_y(z)=\dircon(\varphi^x_y(z),z). 
\end{align*}
Then, by the jet prolongation identity \eqref{eq:n-jet prolongation of direct connection}, these local bisections represent $\Sigma$.

Viceversa, given a smooth map $\sigma:\tilde{\calU}\to \calG$ of local bisections $\sigma^x_y$ satisfying \eqref{eq:jet condition on primitive bisection}, the map 
\begin{align*}
    \Sigma: \calV = \{ (y,x)\in \calU\ |\ \exists z\in M,\ (x,y,z)\in \tilde{\calU} \} \to J^n\calG,
\end{align*}
defined by $\Sigma(y,x) = j^n_x \sigma^x_y$ is automatically a jet prolongation of the direct connection $\dircon$ along the direct connection $\Delta$ given by 
$$
\Delta(y,x) = j^n_x \widehat{\varphi^x_y},
$$
where $\widehat{\varphi^x_y}(z)=(\varphi^x_y(z),z)$ is the local bisection of $\Pair(M)$ associated to the diffeomorphisms  $\varphi^x_y$ provided by eq. \eqref{eq:n-jet prolongation of direct connection}. 

For the infinitesimal connections, let $y(t)$ be a small geodesic starting at $x$ in the direction of $Y\in T_xM$. Then, according to Proposition \ref{prop:dircttoinf}, we have 
\begin{align*} 
\delta^\Sigma(Y) &=\frac{d}{dt}\Sigma(y(t),x)\big|_{t=0}
= \frac{d}{dt} j^n_x\left(z\mapsto \dircon(\varphi^x_{y(t)}(z),z)\right)\big|_{t=0} \\
&= j^n_x\left(z \mapsto \frac{d}{dt} \dircon(\varphi^x_{y(t)}(z),z)\big|_{t=0}\right) 
= j^n_x\left(z \mapsto \delta^\dircon\left(\frac{d}{dt} \varphi^x_{y(t)}(z)\big|_{t=0}\right)\right) ,
\end{align*}
where the equality between the first and the second line is due to the fact that the derivatives in $t$ and in the coordinates of $z$ commute since all the maps are smooth and the variables $t$ and $z$ are mutually independent. 
\end{proof}

After this characterisation, direct connections on $J^n\calG$ can be studied via the collection of their representative local bisections.

\begin{example}[Jet prolongation of direct connections on frame groupoids]
\label{ex:jet prolongation of Gamma on Iso(E)}
Let $\dircon:\Pair(M)\locmap \GL(E)$ be a direct connection on a frame groupoid defined on a diagonal domain $\mathcal{U}\subset \Pair(M)$.  
Then, for any $(y_0,x_0)$ in $\mathcal{U}$, $\dircon(y_0,x_0):E_{x_0}\to E_{y_0},\ e_{x_0}\mapsto \dircon(y_0, x_0)\cdot e_{x_0}$ is a linear isomorphism. 

Let $\Delta^n:\Pair(M)\locmap J^n\Pair(M)$ be a direct connection on the jet pair groupoid defined on the same diagonal domain $\mathcal{U}$, with jet primitive diffeomorphisms $\varphi^{x_0}_{y_0}$.

\begin{enumerate}[label={(\alph*)}]
\item 
The jet prolongation $\dircon^n:\Pair(M)\locmap J^n\GL(E)$ along $\Delta^n$ is best described by linear operators acting on the jet bundle $J^nE$, as in Example \ref{ex:Action of J^nIso(E) on J^nE}. 
For $(y_0,x_0)\in \mathcal{U}$, the linear isomorphism $\dircon^n(y_0,x_0):J^n_{x_0}E\to J^n_{y_0}E$ is given by
\begin{align*}
\dircon^n(y_0,x_0)\cdot j^n_{x_0}f
& = j_{y_0}^n\Big(y\mapsto \dircon(y,\psi_{x_0}^{y_0}(y))\cdot f(\psi_{x_0}^{y_0}(y))\Big) \\ 
& = j_{x_0}^n\Big(x\mapsto \dircon(\varphi^{x_0}_{y_0}(x),x)\cdot f(x)\Big)\circ j^n_{y_0}\psi_{x_0}^{y_0}, 
\end{align*}      
where $f$ is a local section of $E$ around $x_0$ and $\psi_{x_0}^{y_0}$ is the inverse diffeomorphism of $\varphi^{x_0}_{y_0}$.

\item 
If $E=M\times \R^r$ is the trivial bundle on $M=\R^d$, with constant global section a $\mathbf{1}: x\mapsto \mathbf{1}_x=(x,(1,\ldots,1))$, consider the trivial direct connection on $\GL(E)$ given by $\dircon(y,x)\cdot \mathbf{1}_x = \mathbf{1}_y$. Then, for a section $f: M\to E,\ x\mapsto f(x)\, \mathbf{1}_x$, we have $\dircon(y,x)\cdot f(x)\, \mathbf{1}_x = f(x)\,\, \mathbf{1}_y$. It follows that 
\begin{eqnarray*}
\dircon^n(y_0, x_0)\cdot j_{x_0}^nf\, \bf 1_{x_0}
&=&  j_{y_0}^n\left(y\mapsto \dircon(y,\psi_{x_0}^{y_0}(y))\cdot f(\psi_{x_0}^{y_0}(y))\, \mathbf{1}_{y}\right) \\ 
&=&  j_{y_0}^n\left(y\mapsto f(\psi_{x_0}^{y_0}(y))\, \mathbf{1}_{\psi_{x_0}^{y_0}(y)}\right) \\ 
&=&  \left(j_{x_0}^nf \circ j_{y_0}^n\psi_{x_0}^{y_0}\right)\, \mathbf{1}_{y_0} \in J^n_{y_0}E.
\end{eqnarray*} 

\item 
If moreover $M=\R^d$ is equipped with the trivial connection $\nabla^M$, then $\psi_{x_0}^{y_0}(y)= y+x_0-y_0$ and the above formula boils down to:
\begin{align*}
\dircon^n(y_0,x_0)\cdot j_{x_0}^n f\,\textbf{1}_{x_0} = j_{x_0}^n f\,\mathbf{1}_{y_0}.
\end{align*} 
\end{enumerate}
\end{example}

We conclude this section by emphasising that the existence of a direct connection on a jet groupoid forces it to be a gauge groupoid, but it does not force the direct connection to be a jet prolongation. 

\begin{remark}\label{rem:dircon on J^nG}
If a jet groupoid $J^n\calG$ admits a direct connection $\Sigma$, then by Proposition \ref{prop:direct} it is the gauge groupoid $J^n\calG \cong \calG(Q)$ of some principal bundle $Q\to M$. On the other hand, since $\Sigma$ induces a direct connection $\Sigma_0$ on $\calG$, the latter is also the gauge groupoid $\calG \cong \calG(P)$ of some principal bundle $P\to M$. The isomorphism \eqref{eq:J^nG(P)=G(W^nP)} then implies that $Q \cong W^nP$, and therefore $J^n\calG$ is the gauge groupoid of a jet prolonged principal bundle.

However, this does not ensure that the direct connection $\Sigma$ on $J^n\calG$ is necessarily the $n$-jet prolongation of its $0$-th order projection $\Sigma_0$ on $\calG$. 
In fact, the identity \eqref{eq:n-jet prolongation of direct connection} is non-trivial: using the facts that $\Sigma_0=\dircon$ and $\aa^n\circ \Sigma=\Delta_\Sigma=\Delta^n$, the equality 
\begin{align*}
    \Sigma & = J^n \dircon \circ \Delta^n = J^n (\dircon \circ a) \circ \Sigma 
\end{align*}
implies that $J^n (\dircon \circ a)$ fixes the image of $\Sigma$, and this is not always ensured, as shown by next example. 
\end{remark}

\begin{example} \label{ex:dircon not jet prolongation} [Direct connection on a jet groupoid which is not a jet prolongation]. 
Consider the trivial groupoid $\calG$ of Example \ref{ex:trivialGgroupoid}, with $M=\R$ and structure group $(\R,+)$, so that 
\begin{align*}
    \aa(y,g,x)=(y,x), \qquad (z,h,y)\cdot (y,g,x)=(z,h+g,x),\qquad 1_x = (x,0,x). 
\end{align*}
For any $x,y\in \R$, consider the global bisection of $\calG$
\begin{align*}
    \sigma^x_y(z) &= (z+(y-x),(y-x)(z-x),x), \qquad z\in \R
\end{align*}
associated to the translation diffeomorphism $\varphi^x_y(z)=(z+(y-x),z)$, and the global direct connection $\Sigma:\Pair(\R)\to J^1\calG$ given by 
\begin{align*}
    \Sigma(y,x) &= j^1_x \sigma^x_y = j^1_x\big(z\mapsto (z+(y-x),(y-x)(z-x),x)\big).     
\end{align*}
Then, its $0$-th order direct connection on $\calG$ is
\begin{align*}
    \dircon(y,x) &= \pi^1_0 \Sigma(y,x) = \sigma^x_y(x) = (y,0,x),
\end{align*}
and the induced direct connection on $J^1\Pair(\R)$ is 
\begin{align*}
    \Delta_\Sigma(y,x) = j^1_x\big(z\mapsto (z+(y-x),z\big). 
\end{align*}
Therefore, the jet prolongation of $\dircon$ along $\Delta_\Sigma$ is 
\begin{align*}
    \big(J^n\dircon \circ \Delta_\Sigma\big)(y,x) &= j^1_x\big(z\mapsto (z+(y-x),0,z)\big),
\end{align*}
which clearly differs from $\Sigma(y,x)$. 
\end{example}

 
\subsection{Jet prolongation of direct connections on the pair groupoid}
\label{subsec:jet direct connections on Pair(M)}

In this section we construct three canonical direct connections on $J^n\Pair(M)$, starting from a given torsion-free linear connection $\nabla^M$ on $TM$, and the  associated infinitesimal connections on $A(J^n\Pair(M))\cong J^nTM$. 

To define them, we adopt the general strategy outlined in Lemma \ref{lem:jet prolongation Delta^n}, which amounts to finding suitable jet primitive bisections $\widehat{\varphi^x_y}(z)=(\varphi^x_y(z),z)$ of $\Delta^n$ on $J^n\Pair(M)$, given by local diffeomorphisms $\varphi^x_y$ of $M$ such that $\varphi^x_y(x)=y$ and $\varphi^x_x(z)=z$. 
In particular, $\Delta^0$ is necessarily the identity map on $\Pair(M)$, and $\Delta^1:\Pair(M)\to J^1\Pair(M)\cong \GL(TM)$ is given by the differential of $\varphi^x_y$ at $x$, namely 
\begin{align} \label{eq:easy Delta^1}
\Delta^1(y,x)=(d\varphi^x_y)_x:T_xM\overset{\cong}{\to} T_yM. 
\end{align}

Let us fix some notations. 
For any point $x\in M$, denote by $\exp_x:V_{0_x}\subset T_xM\to M$ the exponential map along geodesics. 
As in Section \ref{sec:direct connections}, let 
\begin{align*}
    \calU &= \{ (y,x)\in M\times M\ |\ \text{there exists a convex open subset of $M$ containing both $x$ and $y$} \}
\end{align*}
be the largest convex diagonal domain in $\Pair(M)$. For any $(y,x)\in \calU$, denote by $U\subset M$ a convex open set such that $(y,x)\in U\times U\subset \calU$ (omitting the dependence of $U$ on the points $x$ and $y$).  
For any $(y,x)\in \calU$, denote by $\Par(y,x):T_xM\to T_yM$ the parallel transport of $\nabla^M$ along the unique geodesic from $x$ to $y$. 
In the sequel, we will use \emph{Jacobi vector fields} and related \emph{variation curves through geodesics}, see \cite[Ch.6 and Ch.10]{Lee:1997} or \cite[Ch.5]{doCarmo:1992}\footnote{The references \cite[Ch.6 and Ch.10]{Lee:1997} or \cite[Ch.5]{doCarmo:1992} treat Jacobi fields for Riemannian manifolds and the related Levi-Civita connection $\nabla^M$, but the results hold for any torsion-free linear connection on a manifold.}. For this, let us denote by $R^M$ the curvature tensor of $\nabla^M$, and by $D_t V(t)=\nabla^M_{\dot{y}(t)} V(t)$ the covariant derivative along a curve $y(t)$ of a vector field $V(t)= V(y(t))$ defined along the curve. 
We recall that a \textbf{Jacobi field} is any vector field $J(t)$ along a geodesic $y(t)=\exp_x(tY)$ such that 
\begin{align} \label{eq:Jacobi equation}
    D_t^2 J(t) + R^M_{y(t)}(J(t),\dot{y}(t))\, \dot{y}(t) = 0 \qquad \text{for all $t\in [0,1]$},
\end{align}
which is uniquely determined by its initial conditions $J(0)$ and $D_tJ(0)$ \cite[Definition 2.1]{doCarmo:1992}. 

In the sequel, we denote by $J^A_{B,C}(t)$ the Jacobi field along the geodesic at $x$ in direction $A$, with initial conditions $J(0)=B$ and $D_t J(t)\big|_{t=0}=C$. 


\subsubsection{Parallel jet direct connection}
\label{subsubsec:parallel jet dircon on Pair(M)}

\begin{defn}\label{def:tau^n}
We call \textbf{parallel direct connection} the local map $\Par^n:\Pair(M) \locmap J^n\Pair(M)$ defined on any $(y,x) \in \mathcal{U}$ by 
\begin{align*}
\Par^n(y,x) = j^n_x\big(z\mapsto \widehat{\alpha^x_y}(z)\big), 
\end{align*}
with jet primitive bisections $\widehat{\alpha^x_y}(z)=(\alpha^x_y(z),z)$ given by the local diffeomorphisms 
\begin{align} \label{eq:parallel bisection}
\alpha^x_y:U_x\to M,\ z\mapsto \alpha^x_y(z)=\exp_y\big(\Par(y,x)Z\big), 
\end{align}
where $Z=\exp_x^{-1}(z)\in T_xM$. 
\end{defn}

\begin{thm} \label{thm:parallel direct connection}
The parallel direct connection $\Par^n$ is the higher order extension of the parallel transport $\Par(y,x):T_xM\to T_yM$ along geodesics. 
Its infinitesimal connection $\delta^{\Par^n}:TM\to J^nTM$ is
\begin{align*} 
\delta^{\Par^n}(Y) &= j^n_x\left(z\mapsto J^Z_{Y,0}(1)\right), \qquad Y\in T_xM
\end{align*}
where $J^Z_{Y,0}(s)$ is the Jacobi field along the geodesic $z(s)=\exp_x(sZ)$ from $x$ to $z\in U_x$ with initial conditions $J(0)=Y$ and $D_s J(s)\big|_{s=0}=0$.
\end{thm}

\begin{proof}
For the first assertion we use \eqref{eq:easy Delta^1}. For $z\in U_x$, let $z(s)=\exp_x(sZ)$ be the geodesic from $x$ to $z$ with $Z=\dot{z}(0)\in T_xM$. Then we have 
\begin{align*}
\Par^1(y,x) Z &= (d \alpha^x_y)_x(\dot{z}(0)) 
= \frac{d}{ds}\Big(\alpha^x_y(z(s))\Big)\big|_{s=0} 
= \frac{d}{ds}\Big(\exp_y(\Par(y,x)\, sZ)\Big)\big|_{s=0} 
= \Par(y,x)\,Z. 
\end{align*}
Note also that, by the isomorphisms of groupoids $\calG(FM)\cong \GL(TM)$, this can equivalently be seen as the jet composition with $j^1_x \alpha^x_y$ and gives the parallel transport $\Par^{FM}(y,x):F_xM\to F_yM$ on the frame bundle of $M$, that is, 
\begin{align*}
\Par^{FM} (y,x)(\psi) = j^1_x {\alpha^x_y}\circ j^1_0\tilde{\psi} = j^1_0(\alpha^x_y\circ \tilde{\psi})
\end{align*}
where a frame $\psi:\R^d\overset{\cong}{\to} T_xM$ is seen as the differential at $0$ of a local diffeomorphism $\tilde{\psi}:\R^d\locdif M$ such that $\tilde{\psi}(0)=x$, and therefore it is seen as a jet $j^1_0\tilde{\psi}$. 

For the second assertion we use eq. \eqref{eq:easy delta^n}, and we need to compute the derivative $\frac{d}{dt} \alpha^x_{y(t)}(z)\big|_{t=0}$. To do this, we consider the variation curve of the geodesic $z(s)$ along the geodesic $y(t)=\exp_x(tY)$ from $x$ to $y$ (cf. Fig. \ref{fig:variation alpha}) given by  
\begin{align} \label{eq:alpha} 
\alpha(t,s) = \alpha^x_{y(t)}(z(s)) =\exp_{y(t)}\left(sZ(t)\right), \qquad (t,s)\in [0,1]\times [0,1]
\end{align}
where we choose $Z(t)=\Par(y(t),x)\,Z$. This variation curve has 
\begin{itemize}
\item initial data $\alpha(0,s)=z(s)$, $\alpha(t,0)=y(t)$ and end-point $\alpha(1,1)=\exp_y(\Par(y,x)Z)=\alpha^x_y(z)$, 
\item geodesic main curves $s\mapsto \alpha(t,s)$ in the direction of $\partial_s \alpha(t,s)\big|_{s=0} = Z(t)$,
\item transversal curves $t\mapsto \alpha(t,s)$ in the direction of the Jacobi field $J_\alpha(s)=\partial_t \alpha(t,s)\big|_{t=0}$ along $z(s)$, so that $\alpha(t,s)=\exp_{z(s)}(tJ_\alpha(s))$.
\end{itemize}
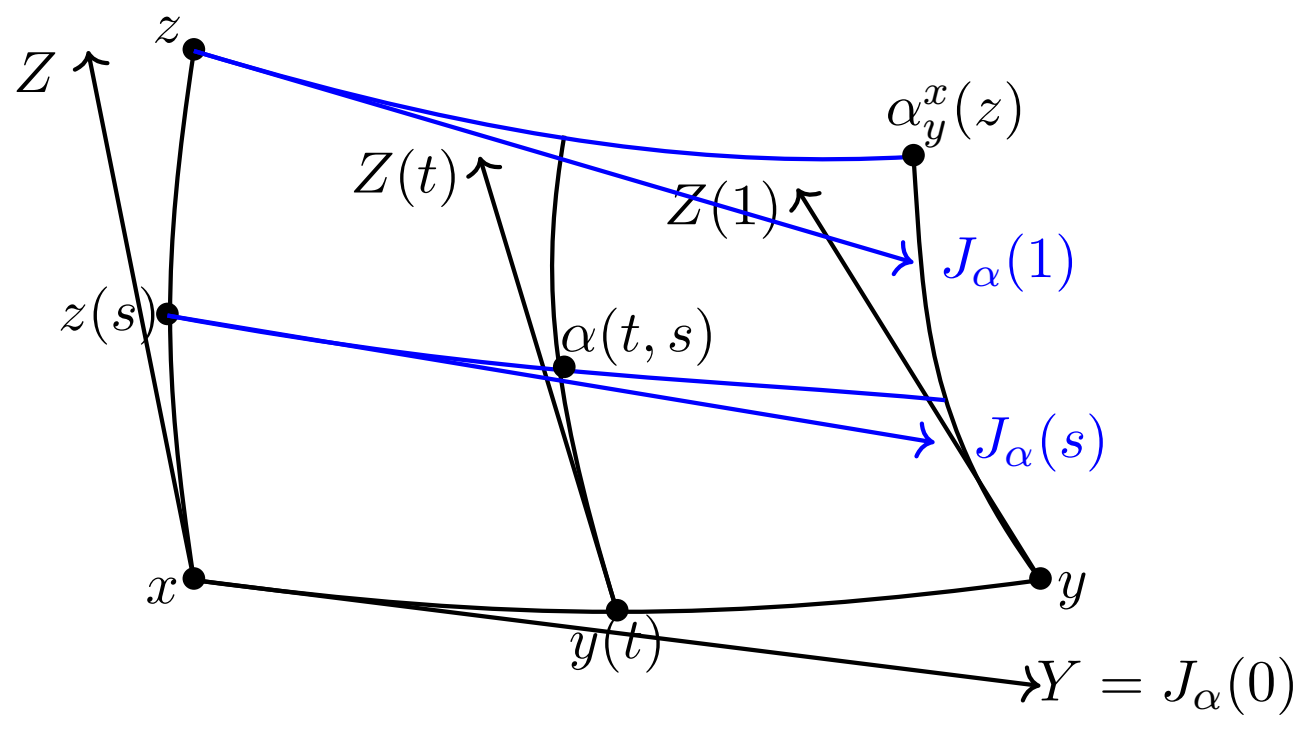
\begin{figure}[h]
\centering
\begin{tikzpicture}[scale=.7]
\draw (0,0) node {$\bullet$}; 
\draw (-0.3,-0.1) node {$x$};
\draw[thick] (0,0) 
      .. controls (3,-0.4) and (5,-0.4) .. (8,0);
\draw[thick] (0,0)
      .. controls (-0.3,2) and (-0.3,3) .. (0,5);
\draw[thick,->] (0,0) -- (-1,5);
\draw (-1.5,4.8) node {$Z$}; 
\draw (4,-0.3) node {$\bullet$};
\draw (4,-0.6) node {$y(t)$};
\draw[thick] (4,-0.3)
      .. controls (3.3,2) and (3.3,3) .. (3.5,4.2);
\draw[thick,->] (4,-0.3) -- (2.7,4); 
\draw (2,3.8) node {$Z(t)$};
\draw (8,0) node {$\bullet$}; 
\draw (8.3,-0.1) node {$y$};
\draw[thick] (8,0)
      .. controls (6.9,1.5) and (6.9,2.5) .. (6.8,4);
\draw[thick,->] (8,0) -- (5.7,3.7); 
\draw (5,3.5) node {$Z(1)$};
\draw[thick,->] (0,0) -- (8,-1); 
\draw (9.2,-1) node {$Y = J_\alpha(0)$};
\draw (-0.25,2.5) node {$\bullet$};
\draw (-0.8,2.5) node {$z(s)$};
\draw[thick,blue] (-0.25,2.5)
      .. controls (3,1.9) and (5,1.9) .. (7.1,1.7);
\draw[thick,blue,->] (-0.25,2.5) -- (7,1.3); 
\draw[blue] (8,1.3) node {$J_\alpha(s)$};
\draw (0,5) node {$\bullet$}; 
\draw (-0.25,5.2) node {$z$};
\draw[thick,blue] (0,5)
      .. controls (3,4.1) and (5,3.9) .. (6.8,4);
\draw[thick,blue,->] (0,5) -- (6.8,3); 
\draw[blue] (7.7,3) node {$J_\alpha(1)$};
\draw (3.5,2) node {$\bullet$};
\draw (4.2,2.3) node {$\alpha(t,s)$};
\draw (6.8,4) node {$\bullet$};
\draw (7.2,4.4) node {$\alpha^x_y(z)$};
\end{tikzpicture}
\caption{variation curve $\alpha(t,s)$ of the geodesic $z(s)$, with vertical main curves in the direction of $Z(t)=\exp_{y(t)}(\Par(y(t),x)Z)$ and horizontal transversal curves in the direction of the Jacobi field $J_\alpha(s)=J^Z_{Y,0}(s)$.}\label{fig:variation alpha}
\end{figure} 
Using eq. \eqref{eq:easy delta^n}, we then have
\begin{align*}
\delta^{\Par^n}(Y) &= j^n_x\left(z\mapsto \frac{d}{dt} \alpha^x_{y(t)}(z)\big|_{t=0}\right) 
= j^n_x\left(z\mapsto \partial_t \alpha(t,1)\big|_{t=0}\right)
= j^n_x\left(z\mapsto J_\alpha(1)\right). 
\end{align*}
The Jacobi field $J_\alpha(s)$ satisfies the initial conditions $J_\alpha(0)=\partial_t \alpha(t,0)\big|_{t=0}=Y$ and
\begin{align*}
D_s J_\alpha(s)\big|_{s=0} &= D_s \left(\partial_t\alpha(t,s)\big|_{t=0}\right)\big|_{s=0} 
= D_t \left(\partial_s\alpha(t,s)\big|_{s=0}\right)\big|_{t=0} 
= D_t Z(t) = 0, 
\end{align*}
where we used \cite[Lemma 6.3]{Lee:1997} and the fact that $Z(t)$ is parallel along $y(t)$. 
Therefore $J_\alpha$ coincides with $J^Z_{Y,0}$. 
\end{proof}


\subsubsection{Jacobi jet direct connection}
\label{subsubsec:Jacobi jet dircon on Pair(M)}

\begin{defn}\label{def:Delta^n}
Similarly, we call \textbf{Jacobi direct connection} the local map $\Jac^n:\Pair(M) \locmap J^n\Pair(M)$ defined on any $(y,x)\in\mathcal{U}$ by 
\begin{align}\label{eq:Delta^n}
\Jac^n(y,x) = j^n_x\big(z\mapsto \widehat{\beta^x_y}(z)\big), 
\end{align}
with jet primitive bisections $\widehat{\beta^x_y}(z) = (\beta^x_y(z), z)$ given by the local diffeomorphisms 
\begin{align*} 
\beta^x_y:U_x\to M,\ z\mapsto \beta^x_y(z)=\exp_z\big(\Par(z,x)Y\big), 
\end{align*}
where $Y=\exp_x^{-1}(y)\in T_xM$. 
\end{defn}

\begin{thm} \label{thm:Jacobi direct connection}
The Jacobi direct connection $\Jac^n$ is the higher order extension of the Jacobi transport map $\Jac^1(y,x):T_xM\to T_yM$ given by
\begin{align*}
    \Jac^1(y,x) Z = J^Y_{Z,0}(1), \qquad Z\in T_xM
\end{align*}
where $J^Y_{Z,0}(t)$ is the Jacobi field along the geodesic $y(t)=\exp_x(tY)$ from $x$ to $y$ with initial conditions $J(0)=Z$ and $D_t J(t)\big|_{t=0}=0$.

Its infinitesimal connection $\delta^{\Jac^n}:TM\to J^nTM$ is the higher order extension of the parallel transport along geodesics, that is, 
\begin{align*} 
\delta^{\Jac^n}(Y) &= j^n_x\left(z\mapsto \Par(z,x)Y\right). 
\end{align*}
\end{thm}

\begin{proof}
For the first assertion, we use again eq. \eqref{eq:easy Delta^1} and get $\Jac^1(y,x) = (d \beta^x_y)_x:T_xM\to T_yM$. For any $Z\in T_xM$ small enough to have $z=\exp_x(Z)\in U_x$ and a geodesic $z(s)=\exp_x(sZ)$ from $x$ to $z$ (if $Z$ is not small enough we rescale it and consider a geodesic on a smaller interval of $s$), we consider the variation curve of $y(t)$ along $z(s)$ (cf. Fig. \ref{fig:variation beta}) given by
\begin{align} \label{eq:beta}
\beta(s,t) = \beta^x_{y(t)}(z(s)) = \exp_{z(s)}\left(t Y(s)\right),\qquad (s,t)\in [0,1]\times [0,1]
\end{align}
where we choose $Y(s)=\Par(z(s),x)Y$. This variation curve has 
\begin{itemize}
\item initial data $\beta(0,t)=y(t)$, $\beta(s,0)=z(s)$ and end-point $\beta(1,1)=\exp_z(\Par(z,x)Y)=\beta^x_y(z)$, 
\item geodesic main curves $t\mapsto \beta(s,t)$ in the direction of $\partial_t \beta(s,t)\big|_{t=0}=Y(s)$, 
\item transversal curves $s\mapsto \beta(s,t)$ in the direction of the Jacobi field $J_\beta(t)=\partial_s \beta(s,t)\big|_{s=0}$ along the geodesic $y(t)$, so that $\beta(s,t)=\exp_{y(t)}(sJ_\beta(t))$. 
\end{itemize}
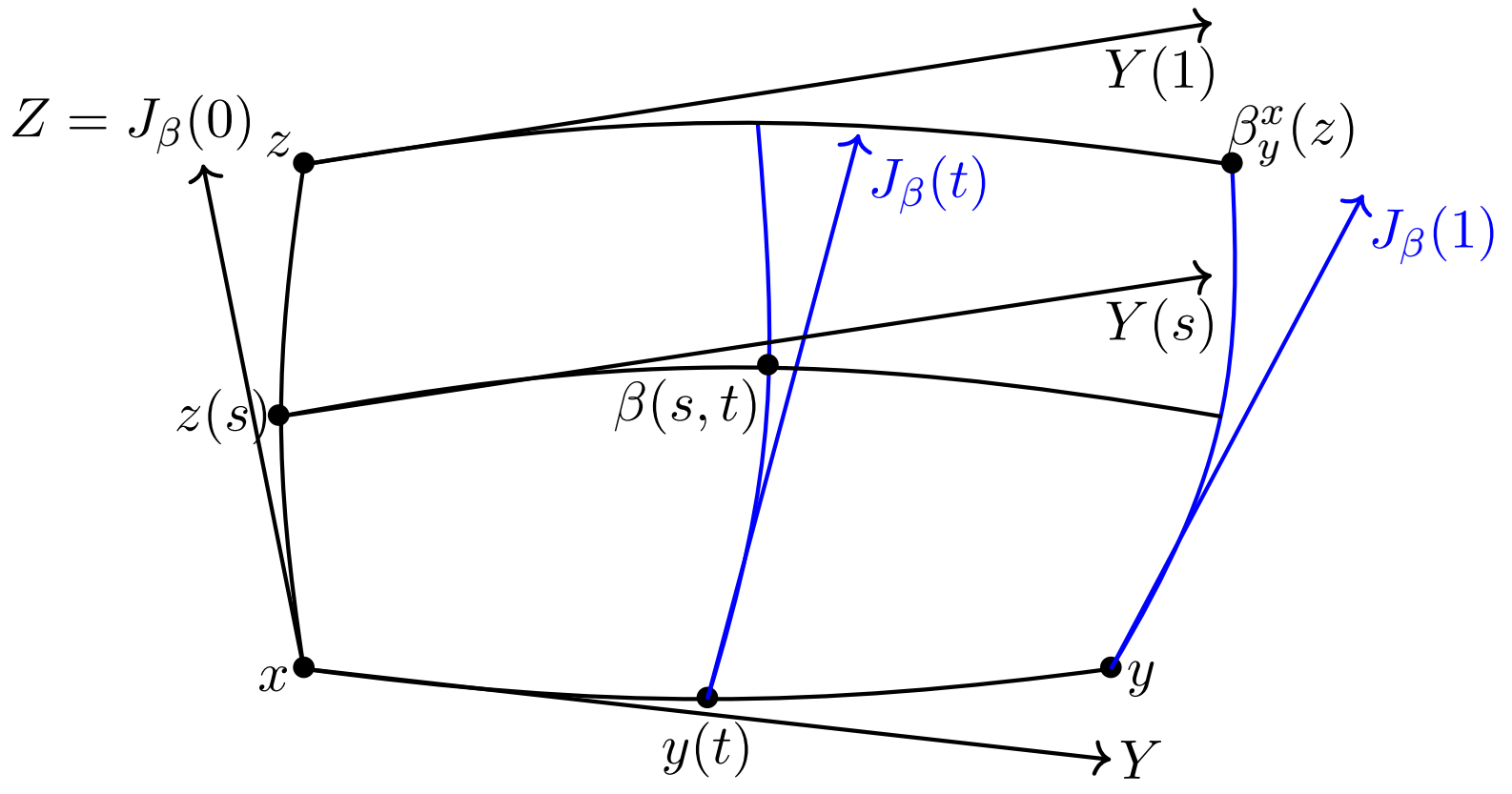
\begin{figure}[h]
\centering
\begin{tikzpicture}[scale=.7]
\draw (0,0) node {$\bullet$}; 
\draw (-0.3,-0.1) node {$x$};
\draw[thick] (0,0) 
      .. controls (3,-0.4) and (5,-0.4) .. (8,0);
\draw[thick] (0,0)
      .. controls (-0.3,2) and (-0.3,3) .. (0,5);
\draw[thick,->] (0,0) -- (-1,5);
\draw (-1.7,5.4) node {$Z=J_\beta(0)$}; 
\draw (4,-0.3) node {$\bullet$};
\draw (4,-0.8) node {$y(t)$};
\draw[thick,blue] (4,-0.3)
      .. controls (4.7,2) and (4.7,3) .. (4.5,5.4);
\draw[thick,blue,->] (4,-0.3) -- (5.5,5.3); 
\draw[blue] (6.2,4.8) node {$J_\beta(t)$};
\draw (8,0) node {$\bullet$}; 
\draw (8.3,-0.1) node {$y$};
\draw[thick,blue] (8,0)
      .. controls (9.2,2) and (9.3,3) .. (9.2,5);
\draw[thick,blue,->] (8,0) -- (10.5,4.7); 
\draw[blue] (11.2,4.3) node {$J_\beta(1)$};
\draw[thick,->] (0,0) -- (8,-0.9); 
\draw (8.3,-0.9) node {$Y$};
\draw (-0.25,2.5) node {$\bullet$};
\draw (-0.8,2.5) node {$z(s)$};
\draw[thick] (-0.25,2.5)
      .. controls (3,3.1) and (5,3.2) .. (9.1,2.5);
\draw[thick,->] (-0.25,2.5) -- (9,3.9); 
\draw (8.5,3.4) node {$Y(s)$};
\draw (0,5) node {$\bullet$}; 
\draw (-0.25,5.2) node {$z$};
\draw[thick] (0,5)
      .. controls (3,5.5) and (5,5.6) .. (9.2,5);
\draw[thick,->] (0,5) -- (9,6.4); 
\draw (8.5,5.9) node {$Y(1)$};
\draw (4.6,3) node {$\bullet$}; 
\draw (3.8,2.6) node {$\beta(s,t)$};
\draw (9.2,5) node {$\bullet$};
\draw (9.8,5.3) node {$\beta^x_y(z)$};
\end{tikzpicture}
\caption{variation curve $\beta(s,t)$ of the geodesic $y(t)$, with horizontal main curves in the direction of $Y(s)=\exp_{z(s)}(\Par(z(s),x)Y)$ and vertical transversal curves in the direction of the Jacobi field $J_\beta(t)=J^Y_{Z,0}(t)$.}\label{fig:variation beta}
\end{figure} 
Then we have 
\begin{align*}
\Jac^1(y,x)(Z) &= (d \beta^x_y)_x(\dot{z}(0)) 
= \frac{d}{ds} \beta^x_y(z(s))\big|_{s=0}
= \partial_s \beta(s,1)\big|_{s=0} 
= J_\beta(1).
\end{align*}
Since $J_\beta(t)$ has initial values $J_\beta(0)=Z$ and $D_t J_\beta(t)\big|_{t=0}=D_s Y(s)\big|_{s=0}= D_s \Par(z(s),x)Y|_{s=0}=0$, it coincides with $J^Y_{Z,0}$. 

Finally, let us show that the map $\Jac^1(y,x): T_xM\to T_yM,\ Z\mapsto J^Y_{Z,0}(1)$ defines a direct connection on $\GL(TM)\cong J^1\Pair(M)$. In fact, it is well defined for $(y,x)\in \mathcal{U}$, it is linear in $Z$ because $J^Y_{Z,0}(t)=J_\beta(t)$ is the solution of an ODE whose initial conditions are linear in $Z$, and finally it is invertible: its inverse is the linear map $\Phi(x,y):T_yM\to T_xM,\ U\mapsto \tilde{J}(1)$ where $\tilde{J}(t):=\partial_s\tilde{\beta}(s,t)$ is the Jacobi field detemined by the inverse variation curve $\tilde{\beta}(s,t)=\beta(s,1-t)\big|_{s=0}$ of the inverse geodesic $x(t)=y(1-t)$ from $y$ to $x$, along the geodesic $u(s)=\exp_y(sU)$, hence $\tilde{J}(0)=U$, and with second initial condition 
$$
D_t\tilde{J}(t)\big|_{t=0} = D_s[-\Par(u(s),z(s)\,\Par(z(s),x)\,Y]|_{s=0}
$$
instead of $D_s[-\Par(u(s),y)\,\Par(y,x)\,Y]|_{s=0}$. 
The map $\Phi(x,y)$ is indeed the inverse of $\Jac^1(y,x)$, because $\tilde{J}(t)=J_\beta(1-t)$: if $U=J_\beta(1)=\tilde{J}(0)$, then $\tilde{J}(1)=J_\beta(0)=Z$.  

For the second assertion, we use again eq. \eqref{eq:easy delta^n} and get 
\begin{align*}
\delta^n(Y) &= j^n_x\left(z\mapsto \frac{d}{dt} \beta^x_{y(t)}(z)\big|_{t=0}\right) 
= j^n_x\left(z\mapsto \partial_t \beta(1,t)\big|_{t=0}\right) 
= j^n_x\left(z\mapsto Y(1)\right) 
= j^n_x\left(z\mapsto \Par(z,x)\,Y\right). 
\end{align*}
\end{proof}

\begin{prop}\label{prop:different bisections}
The local maps $\Par^n$ and $\Jac^n$ define two generally distinct direct connections on $J^n\Pair(M)$. 
They coincide if and only if $\nabla^M$ is flat (that is, $M$ is an affine manifold). 
\end{prop}

\begin{proof}
As shown in the previous proofs, the primitive bisections of the two direct connections are associated to two variation curves $\alpha(s,t)$ and $\beta(t,s)$ with the same initial data but distinct end-points, velocity vectors and Jacobi fields. 
The two variation curves are necessarily distinct, that is, $\alpha(t,s)\neq \beta(s,t)$ except for $s=0$ or $t=0$. In fact, the equality $\exp_{y(t)}\left(sZ(t)\right)=\exp_{z(s)}\left(tY(s)\right)$ means that $Y(s)=J_\alpha(s)$ is the variation field of $\alpha(s,t)$ and $Z(t)=J_\beta(t)$ is the variation field of $\beta$. These equalities in general do not hold, because for generic points $y,z\in U_x$ (hence generic "small" tangent vectors $Y,Z\in T_xM$) the fields $Y(s)=\Par(z(s),x)Y$ and $Z(t)=\Par(y(t),x)Z$ are not Jacobi fields: 
\begin{align*}
    D^2_s Y(s) + R^M_{z(s)}(Y(s),\dot{z}(s))\,\dot{z}(s) = \Par(z(s),x)\  R^M_x(Y,Z)\,Z \neq 0, \\ 
    D^2_t Z(t) + R^M_{y(t)}(Z(t),\dot{y}(t))\,\dot{y}(t) = \Par(y(t),x)\  R^M_x(Z,Y)\,Y \neq 0. 
\end{align*}
In particular, the two end-points $\alpha^x_y(z)=\exp_y(\Par(y,x)Z)$ and $\beta^x_y(z)=\exp_z(\Par(z,x)Y)$ are distinct. 

Finally, the two direct connections $\Par^n$ and $\Jac^n$ coincide if and only if we have $\Par^n(y,x)=\Jac^n(y,x)$ for any points $x,y\in M$ sufficiently close to each other. 
In particular, since $\Par^1(y,x)=\Jac^1(y,x)$, the Jacobi fields are given by $J^Y_{Z,0}(t)=\Par(y(t),x)Z$ for any $t\in [0,1]$. As parallel transport satisfies $D_t(\mathcal P(y(t),x)Z)=0$, it follows that $D_t^2J_Z(t)=0$. The Jacobi equation therefore yields $R_x(Z,Y)Y=0$ for all $Y,Z\in T_xM$. By polarisation and the first Bianchi identity, this implies $R_x=0$, and hence $R=0$.
\end{proof}


\subsubsection{Exponential jet direct connection}
\label{subsubsec:exponential jet dircon on Pair(M)}

\begin{defn}\label{def:Exp^n}
We call \textbf{exponential direct connection} on $J^n\Pair(M)$ the local map $\Exp^n:\Pair(M) \locmap J^n\Pair(M)$ defined on any $(y,x) \in \mathcal{U}$ by 
\begin{align*}
\Exp^n(y,x) = j^n_x\big(z\mapsto \widehat{\gamma^x_y}(z)\big), 
\end{align*}
with jet primitive bisections $\widehat{\gamma^x_y}(z)=(\gamma^x_y(z),z)$ given by the local diffeomorphisms 
\begin{align} \label{eq:parallel_bisection2}
\gamma^x_y:U_x\to M,\ z\mapsto \gamma^x_y(z)=\exp_x\big(Y+Z\big), 
\end{align}
where $Y,Z\in T_xM$ are defined by $y=\exp_x(Y)$ and $z=\exp_x(Z)$. 
\end{defn}
 
Let us recall that W.M. Mikulski \cite[\S 3]{Mikulski:2007} and I. I. Kol\'a\v{r} \cite[\S 5]{Kolar:2009} introduced a canonical \textbf{exponential jet prolongation} of the Lie algebroid connection $\varepsilon: TM\to J^1TM\cong A(\GL(TM))$ associated to a torsion-free linear connection $\nabla^M$ on $TM$: it is the Lie algebroid connection $\varepsilon^n: TM\longrightarrow J^nTM$ defined on a tangent vector $Y\in T_xM$ as the $n$-jet 
\begin{align*}
\varepsilon^n(Y):= j_{0_x}^{n}\left(d\exp_x(\widetilde{Y})\right) = j^n_x\left(z\mapsto (\exp_x)_*(\widetilde{Y})(z)\right),
\end{align*}
where $\widetilde{Y}$ is the vector field on $T_xM$ obtained by translation of $Y$, that is, $\tilde{Y}(Z)=Z+Y\in T_Z(T_xM)$ for any $Z\in T_xM$, and $(\exp_x)_*: T(T_xM) \to TM$ transforms vector fields on $T_xM$ into vector fields on $M$. 
\begin{thm} \label{thm:exponential direct connection}
The exponential direct connection $\Exp^n$ is the higher order extension of the Jacobi transport map $\Exp(y,x):T_xM\to T_yM$ given by 
\begin{align*}
    \Exp(y,x) Z = (d\exp_x)_Y(Z) = J^Y_{0,Z}(1), \qquad Z\in T_xM
\end{align*}
where $J^Y_{0,Z}(t)$ is the Jacobi field along the geodesic $y(t)=\exp_x(tY)$ from $x$ to $y$ with initial conditions $J(0)=0$ and $D_tJ(t)\big|_{t=0}=Z$.  

Its infinitesimal connection $\delta^{\Exp^n}:TM\to J^nTM$ is the jet prolongation 
\begin{align*}
    \delta^{\Exp^n}_x(Y) = j^n_x(z\mapsto \Exp(z,x)Y) = 
    j^n_x\left(z\mapsto (d\exp_x)_Z(Y)\right)= j^n_x\left(z\mapsto J^Z_{0,Y}(1)\right), \qquad Y\in T_xM 
\end{align*}
which coincides with the Lie algebroid exponential jet prolongation $\varepsilon^n$.
\end{thm}

\begin{proof}
For the first assertion, by eq. \eqref{eq:easy Delta^1} we have $\Exp^1(y,x)=(d\gamma^x_y)_x:T_xM\to T_yM$ with $\gamma_y^x$   as in  eq. \eqref{eq:parallel_bisection2}. Consider the variation curve of $y(t)$ by a spray of geodesics starting at $x$ in directions $Y(s)=Y+sZ$ (cf. Fig. \ref{fig:variation gamma}), given by 
\begin{align*} 
\gamma(s,t) = \exp_x\left(tY(s)\right) = \exp_x\left(tY+stZ)\right), \qquad (s,t)\in [0,1]\times [0,1]. 
\end{align*}
This variation curve has 
\begin{itemize}
\item initial data $\gamma(0,t)=y(t)$, $\gamma(s,0)=x$ and end-point $\gamma(1,1)=\exp_z(Y+Z)=\gamma^x_y(z)$, 
\item geodesic main curves $t\mapsto \gamma(s,t)$ in the direction of $\partial_t \gamma(s,t)\big|_{t=0} = Y(s)$,
\item transversal curves $s\mapsto \gamma(s,t)$ in the direction of the Jacobi field $J_\gamma(t)=\partial_s \gamma(s,t)\big|_{s=0}$ along $y(t)$, so that $\gamma(s,t)=\exp_{y(t)}(sJ_\gamma(t))$.
\end{itemize}

\begin{figure}[h!]
\centering
\begin{tikzpicture}[scale=.7]
\draw (0,0) node {$\bullet$}; 
\draw (-0.3,-0.1) node {$x$};
\draw[thick] (0,0) 
      .. controls (3,-0.4) and (5,-0.4) .. (8,0);
\draw[thick] (0,0)
      .. controls (-0.3,2) and (-0.3,3) .. (0,5);
\draw[thick,->] (0,0) -- (-1,5);
\draw (-1.3,5.4) node {$Z$}; 
\draw (4,-0.3) node {$\bullet$};
\draw (4,-0.8) node {$y(t)$};
\draw[thick,blue] (4,-0.3)
      .. controls (3.8,0.7) and (3.5,1.4) .. (3,2);
\draw[thick,blue,->] (4,-0.3) -- (3.5,2.2); 
\draw[blue] (4.3,1.9) node {$J_\gamma(t)$};
\draw (8,0) node {$\bullet$}; 
\draw (8.3,-0.1) node {$y$};
\draw[thick,blue] (8,0)
      .. controls (7.7,1) and (7.6,2) .. (6,4.8);
\draw[thick,blue,->] (8,0) -- (7,5.1); 
\draw[blue] (7.8,4.7) node {$J_\gamma(1)$};
\draw[thick,->] (0,0) -- (8,-1); 
\draw (8.3,-1) node {$Y$};
\draw (-0.25,2.5) node {$\bullet$};
\draw (-0.8,2.5) node {$z(s)$};
\draw[thick] (0,0)
      .. controls (3,0.6) and (5,1.2) .. (7.2,2.4);
\draw[thick,->] (0,0) -- (7.5,1.5); 
\draw (6.8,0.9) node {$Y+sZ$};
\draw (0,5) node {$\bullet$}; 
\draw (-0.25,5.2) node {$z$};
\draw[thick] (0,0)
      .. controls (3,1.7) and (5,3.4) .. (6,4.8);
\draw[thick,->] (0,0) -- (7,4); 
\draw (6.8,3.2) node {$Y+Z$};
\draw (3.6,0.9) node {$\bullet$}; 
\draw (2.7,1) node {$\gamma(s,t)$};
\draw (6,4.8) node {$\bullet$};
\draw (5.2,5) node {$\gamma^x_y(z)$};
\end{tikzpicture}
\caption{variation curve $\gamma(s,t)$ of the geodesic $y(t)$, with radial main curves in the direction of $Y(s)=\exp_x(Y+sZ)$ and vertical transversal curves in the direction of the Jacobi field $J_\gamma(t)=J^Y_{0,Z}(t)$.}
\label{fig:variation gamma}
\end{figure}
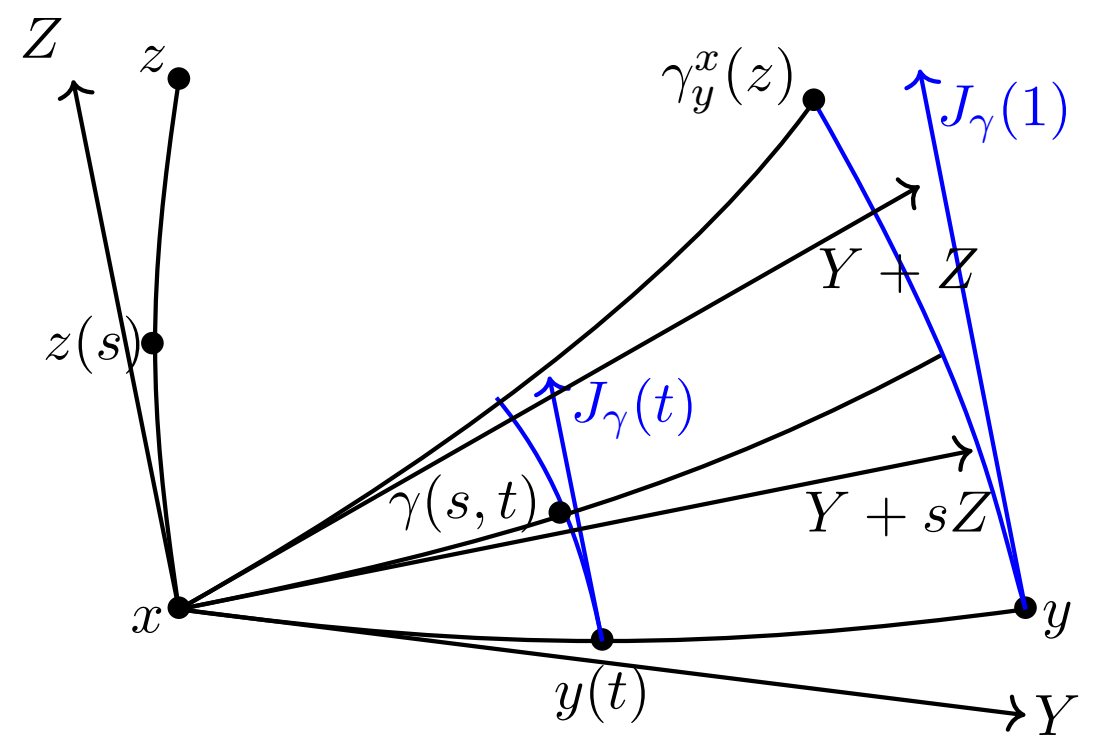 
Since $\gamma^x_y(z(s))=\exp_x(Y+sZ)=\gamma(s,1)$, we have
\begin{align*}
    \Exp^1(y,x)(Z) = (d\gamma^x_y)_x(\dot{z}(0)) 
    = \frac{d}{ds} \gamma^x_y(z(s))\big|_{s=0} 
    = \partial_s \gamma(s,1)\big|_{s=0} = J_\gamma(1).  
\end{align*}
Since $J_\gamma(t)$ has initial values $J_\gamma(0)=0$ and $D_tJ_\gamma(t)\big|_{t=0}=D_s(Y+sZ)\big|_{s=0}=Z$, it coincides with $J^Y_{0,Z}$. 
Finally, we have $J^Y_{0,Z}(1)=\partial_s \gamma(s,1)\big|_{s=0}= (d\exp_x)_Y(Z)$ \cite[Ch.5, Section 2]{doCarmo:1992}. 

We turn to the infinitesimal connection $\delta^{\Exp^n}$. 
Using eq. \eqref{eq:easy delta^n}, we have
\begin{align*}
\delta^{\Exp^n}(Y) &= j^n_x\left(z\mapsto \frac{d}{dt} \gamma^x_{y(t)}(z)\big|_{t=0}\right).
\end{align*}
In order to disentangle $z$ from $y$, we explicitly write the dependence of $Y$ in the variation curve and consider the disentangled variation $u\mapsto uY$ at fixed $z$, replacing the notation $\gamma(s,t)$ by
\begin{align*}
    \gamma_Y(s,t):=\exp_x(t(Y+sZ)) \qquad \gamma_{uY}(1,1) = \exp_x(uY+Z).
\end{align*}
We have
\begin{align*}
\delta^{\Exp^n}(Y) &= j^n_x\left(z\mapsto \partial_u \gamma_{uY}(1,1)\big|_{u=0}\right)
= j^n_x\left(z\mapsto (d\exp_x)_Z(Y)\right).  
\end{align*}
The vector $(d\exp_x)_Z(Y)$ coincides with the Jacobi field of the exponential direct connection $\Exp^1$ described above, but for inverted initial conditions $Z$ and $Y$, which amounts to inverting the roles of the points $y$ and $z$, and we have
\begin{align*}
\delta^{\Exp^n}(Y) &= j^n_x\left(z\mapsto J^Z_{0,Y}(1)\right) = j^n_x\left(z\mapsto \Exp(z,x) Y \right).  
\end{align*}
Since 
\begin{align*}
    (\exp_x)_*(\widetilde{Y})(z) = (d\exp_x)_Z(Y)
\end{align*}
we have the equality of infinitesimal connections $\delta^{\Exp^n} = \varepsilon^n$. 
\end{proof}

Note that the value $\Exp(y,x)^1 Z = \partial_s \gamma(s,1)\big|_{s=0}$ is computed using a \emph{geodesic at $x$ fixed by its endpoint} that depends on $y$ and $z=\exp_x(Z)$, while $\delta^\Exp(Y)$ requires to \emph{vary that geodesic as $y$ moves while $z$ is fixed}. This variation is not computed by the derivative of $\gamma(s,t)$ with respect to $t$ with fixed $s=1$, because the parameter $t$ which makes $y$ move along $y(t)=\exp_x(tY)$ also moves  $Z$ to $tZ$. 

\begin{remark}
For a general manifold $M$ and a generic linear connection $\nabla^M$, the exponential jet direct connection $\Exp^n$ is clearly different from the parallel one $\Par^n$ and from the Jacobi one $J^n$. As in Proposition \ref{prop:different bisections}, they coincide when $\nabla^M$ is flat.
\end{remark}


\subsubsection{Interpolation of   jet direct connections}
\label{subsubsec:interpolation}

A comparison of the three preceding examples reveals a significant underlying relationship\footnote{The underlying geometric meaning can be explained in terms of \emph{smooth local loops} \cite{Sabinin:1999}.}.
We first collect the three families of local diffeomorphisms, the associated direct connections and the induced infinitesimal connections. 
For any $x$ in $M$ and for any $y,z$ sufficiently close to $x$, we set $Y=\exp_x^{-1}(y)$ and $Z=\exp_x^{-1}(z)$. 

We first observe that by Theorem \ref{thm:parallel direct connection}, we have $\delta^{\Par^n}(Y) = j^n(z\mapsto J^Z_{Y,0}(1))$, and exchanging the role of the two points $y$ and $z$, by Theorem \ref{thm:Jacobi direct connection} we have $J^Z_{Y,0}(1)=\Jac^1(z,x)Y$. 
Then, from Theorems \ref{thm:parallel direct connection}, \ref{thm:Jacobi direct connection}, \ref{thm:exponential direct connection}, we can set up the following comparison table:
\bigskip 
\begin{center}
\begin{tabular}{c c c}
\toprule
Diffeomorphisms & Direct connection & Infinitesimal connection \\[3pt]
$\varphi^x_y$ & $\Delta^n(y,x)=j^n_x(z\mapsto \widehat{\varphi^x_y})$ & $\delta^n(Y)=\frac{d}{dt}\big|_0\Delta^n(\exp_x(tY),x)$\\[3pt] 
\midrule
$\alpha^x_y(z)=\exp_y\bigl(\Par(y,x)Z\bigr)$  & $\Par^n(y,x)$ & $j_x^n\bigl(z\mapsto \Jac(z,x)Y\bigr)$\\[3pt]
$\beta^x_y(z)=\exp_z\bigl(\Par(z,x)Y\bigr)$ & $\Jac^n(y,x)$ & $j_x^n\bigl(z\mapsto \Par(z,x)Y\bigr)$\\[3pt]
$\gamma^x_y(z)=\exp_x\bigl(Y+Z\bigr)$ & $\Exp^n(y,x)$ & $j_x^n\bigl(z\mapsto \Exp(z,x)Y\bigr)$\\
\bottomrule
\end{tabular}
\end{center}
\bigskip 
This table shows that the parallel and the Jacobi direct connections are interchanged in the infinitesimal connection, whilst the exponential connection is self-dual (with swapped role of the two points $y$ and $z$).
\bigskip

We now exploit the fact that $J^1\Pair(M) \cong GL(TM)$ is the frame groupoid of $TM$, and hance the fact that the $1$-st order direct connection $\Delta^1(y,x)$ is a linear isomorphism $T_xM \to T_yM$. 
In the three cases, the value $\Delta^1(y,x) Z$ of any $Z$ in $T_xM$ is that of a Jacobi field $J^A_{B,C}$ along the geodesic in the direction $A$ with specific initial conditions $J(0)=B$ and $D_t J(0)=C$. Again, set $Y=\exp_x^{-1}(y)$ and $Z=\exp_x^{-1}(z)$. We then have a second comparison table:

\begin{center}
\begin{tabular}{c c c c}
\toprule
Direct connection & Geodesic velocity & $0$-th order initial cond. & $1$-st order initial cond.\\[3pt]
$\Delta(y,x)Z=J^A_{B,C}(1)$ & $A$ & $B$ & $C$ \\[3pt] 
\midrule
$\Jac(y,x) Z$ & Y & Z & 0 \\[3pt]
$\Exp(y,x) Z$ & Y & 0 & Z \\
\bottomrule
\end{tabular}
\end{center}
\bigskip 
This table shows that $\Jac$ and $\Exp$, with Jacobi fields along the same geodesic with velocity $Y$, mutually exchange the two initial conditions $(Z,0)$ and $(0,Z)$ of the Jacobi field. This exchange admits a continuous interpolation. 

\begin{remark}\label{rem:interpolation}
For any $\theta \in [0,1]$, set $z(\theta)=\exp_x(\theta Z)$. Then, the collection of local diffeomorphisms 
\begin{align*}
    \varphi^x_{\theta,y}(z) &= \exp_{z(\theta)}\big(\Par(z(\theta),x)[Y+(1-\theta)Z]\big), \quad \theta \in [0,1]
\end{align*}
satisfy $\varphi^x_{\theta,y}(x)=y$ and $\varphi^x_{\theta,x}(z)=z$ and therefore define a direct connection on $J^n\Pair(M)$ by
\begin{align*}
    \Delta_\theta^n(y,x) = j^n_x \widehat{\varphi^x_{\theta,y}}. 
\end{align*}
The value of the $1$-st order connection $\Delta_\theta^1(y,x)Z$ is that of the Jacobi field $J^Y_{\theta Z,(1-\theta)Z}(1)$ along the geodesic $\exp_x(tY)$, with initial conditions $J(0)=\theta Z$ and $D_tJ(0)=(1-\theta)Z$. 
Hence the direct connection
\begin{align*}
    \Delta_\theta^1(y,x) = \theta \Jac(y,x) + (1-\theta) \Exp(y,x), 
\end{align*}
is a continuous deformation of $\Exp(y,x)$ into $\Jac(y,x)$. 
Its infinitesimal connection is
\begin{align*}
    \delta_\theta^n(Z) = j^n_x\big(z\mapsto \theta \Exp(z,z_\theta)\Par(z_\theta,x)Z \big).
\end{align*}
where $z_\theta=\exp_x(\theta Z)$. 
\end{remark}

 
\subsection{Flat jet groupoids} \label{subsec:flat jet groupoids}

According to Proposition \ref{prop:flat groupoid}, a jet groupoid $J^n\calG$ on a connected manifold $M$ is flat if and only if it admits a flat direct connection. 
In this section we show that, when $\calG$ is flat, in general  $J^n\calG$ is not necessarily flat, and we give a sufficient condition (on $M$) to have flat jet groupoids.

\begin{example}[Non flat jet groupoid]
Over the $2$-sphere $S^2$, the pair groupoid $\Pair(S^2)$ is flat (as any pair groupoid), because the trivial connection $\Id: \Pair(S^2)\to \Pair(S^2)$ is flat. However its first jet prolongation $J^1\Pair(S^2)$ cannot be flat. 
In fact, since $J^1\Pair(S^2) \cong \GL(TS^2)$, a flat direct connection on $J^1\Pair(S^2)$ would induce a flat linear connection on $TS^2$, which does not exist. 
\end{example}

\begin{thm}\label{thm:flat jet groupoid}
Assume that the manifold $M$ is affine, i.e. $TM$ is equipped with a flat torsion-free linear connection. Then:  
\begin{enumerate}
\item the groupoid $J^n\Pair(M)$ is flat, 
\item the jet prolongation $J^n\calG$ of any flat groupoid $\calG$ on $M$ is flat.
\end{enumerate}
\end{thm}

\begin{proof}
1. If $TM$ is endowed with a flat and torsion-free linear connection $\nabla^M$, with corresponding parallel transport $\Par$, then the parallel direct connection $\Par^n$ on $J^n\Pair(M)$ given in Definition \ref{def:tau^n} is flat, namely it is a groupoid morphism.
Indeed, consider three nearby points $x,y,z \in M$ and 
\begin{align*}
    \Par^n(y,x) &= j^n_x\big(u\mapsto (\alpha^x_y(u),u)\big) &\text{with}\qquad \alpha^x_y(u) = \exp_y(\Par(y,x) U), \qquad u=\exp_x(U) \\ 
    \Par^n(z,y) &= j^n_y\big(v\mapsto (a^y_z(v),v)\big) &\text{with}\qquad a^y_z(v) = \exp_z(\Par(z,y) V), \qquad v=\exp_y(V). 
\end{align*}
Since $\nabla^M$ is flat, we have $\Par(z,y) \Par (y,x) = \Par(z,x)$ and therefore 
\begin{align*}
    \alpha^y_z(\alpha^x_y(u)) &= \exp_z(\Par(z,y) \Par(y,x) U) = \exp_z(\Par(z,x)U) = a^x_z(u). 
\end{align*}
Then, the multiplication in $J^n\Pair(M)$ gives 
\begin{align*}
    \Par^n(z,y) \Par^n(y,x) = j^n_x\big(u\mapsto (\alpha^y_z(\alpha^x_y(u)),u)\big) 
    = j^n_x\big(u\mapsto (a^x_z(u),u)\big) = \Par^n(z,x).
\end{align*}
This proves that the groupoid $J^n\Pair(M)$ is flat.

2. Now let $\calG$ be a flat groupoid, equipped with a flat direct connection $\dircon$. Then its jet prolongation along $\Par^n$ is the composition $\dircon^n = J^n\dircon \circ \Par^n$. Since $J^n\dircon$ is also a groupoid morphism, by Proposition \ref{prop:J^nPhi}, the direct connection $\dircon^n$ is flat. Therefore $J^n\calG$ is a flat groupoid. 
\end{proof}


\section{Direct connections on jet bundles}
\label{sec:direct connections on jet bundles}

In this section we provide an example of a direct connection on the frame groupoid of a jet bundle that is not the jet prolongation of a direct connection on the bundle (Theorem \ref{thm:Taylor direct connection}). The idea comes from the theory of regularity structures in stochastic PDEs, introduced by M. Hairer in \cite{Hairer:2014b} for functions on $\R^n$ and further developed for functions on a manifold by A. Dahlqvist, J. Diehl and B. Driver in \cite{DahlqvistDiehl-Driver:2019}, and for sections of a vector bundle by M. Hairer and H. Singh in \cite{Hairer-Singh:2023}. 
To define this higher direct connection, we start with a direct connection $\dircon$ on the frame groupoid of a vector bundle $E$ and give a local approximation formula for sections of $E$ with displacement along the fibres ruled by $\dircon$. With respect to the usual Taylor expansion of sections, both the principal part and the remainder we introduce contain new extra terms in all covariant derivatives, defined in terms of the derivatives of the curvature of $\dircon$ along geodesics, that we call \emph{transitivity defects} (Theorem \ref{thm:Taylor with nabla}). These terms vanish systematically if $\dircon$ is the parallel transport of a usual connection, so the non-local correction terms to the Taylor expansion are a new feature of general direct connections. 
We conclude with the sketch of a \emph{local} model for $E$-valued jet regularity structures based on the generalised Taylor expansion induced by a general direct connection on $E$ (Remark \ref{rem:dircon on J^nG}). 

 
\subsection{Taylor expansion with respect to direct connections} 
\label{subsec:Taylor expansion} 

In this section we introduce a covariant Taylor expansion of local sections of a vector bundle $E$ induced by a generic direct connection $\dircon$ on $E$. 
Expressed as a polynomial in higher covariant derivatives, the principal part of this expansion comes with extra coefficients containing derivatives of the curvature of $\dircon$ at the base point. 
These new terms also arise in the remainder, at every lower order of derivatives. 

As before, consider a base manifold $M$ equipped with a torsion-free linear connection $\nabla^M$ on $TM$, and a vector bundle $E\to M$ of finite rank equipped with a direct connection $\dircon:\Pair(M)\locmap \GL(E)$ on its frame groupoid.  
Suppose that $\dircon$ is defined on a convex diagonal domain $\calU$. 

\begin{defn}\label{def:Taylor}
Let $x_0$ be a given point in $M$, $U_{x_0}$ a convex neighbourhood of $x_0$ such that $U_{x_0}\times U_{x_0}\subset \calU$ and $\psi$ in $\Gamma_{U_{x_0}}(M,E)$ a smooth local section defined on $U_{x_0}$. 
We call 
\begin{itemize}
\item 
\textbf{$n$-th order Taylor polynomial of $\psi$ at $x_0$} the local section of $E$ defined on any $x\in U_{x_0}$ by 
\begin{align}\label{eq:Taylor}
\Tay_{x_0}^n \psi(x) 
&:= \dircon(x,x_0) \sum_{k=0}^n \frac{1}{k!} \frac{d^k}{dt^k} \dircon(x(t),x_0)^{-1} \psi(x(t))\big|_{t=0},
\end{align} 

\item 
\textbf{$n$-th order remainder of $\psi$ at $x_0$} the local section of $E$ defined on any $x\in U_x$ by 
\begin{align}\label{eq:Remainder}
\Rem_{x_0}^n \psi(x) 
&:= \dircon(x,x_0) \frac{1}{n!} \int_0^1 (1-t)^n \frac{d^{n+1}}{dt^{n+1}} \dircon(x(t),x_0)^{-1} \psi(x(t))\, dt,
\end{align} 
\end{itemize}
where $x(t)$ is the geodesic from $x_0$ to $x$, defined for $t\in [0,1]$. 
\end{defn}

\begin{defn} \label{def:function F}
Fix $x_0\in M$. It is convenient to consider the open neighborhood of the null vector in $T_{x_0}M$ given by 
\begin{align*}
    V_{x_0} &= \{X\in T_{x_0}M\ |\ (x,x_0)\in \calU,\ x=\exp_{x_0}(X)\}, 
\end{align*}
and to introduce the vector valued function $F: V_{x_0}\to E_{x_0}$ 
that \textbf{recentres the value $\psi(x)$ at $x_0$} using $\dircon$, that is 
\begin{align*} 
F(X)=\dircon(x,x_0)^{-1}\psi(x), \qquad X\in V_{x_0},\ x =\exp_{x_0}(X).
\end{align*}
Then, setting $x(t) = \exp_{x_0}(tX)$ for any $t\in [0,1]$, we have 
\begin{align*}
    F(tX) &= \dircon(x(t),x_0)^{-1}\psi(x(t)), \qquad t\in [0,1] \\ 
    F(0) &= \psi(x_0), 
\end{align*}
and therefore, for $x=\exp_{x_0}(X)$, we have 
\begin{align*}
    \Tay_{x_0}^n \psi(x) &= \dircon(x,x_0) \sum_{k=0}^n \frac{1}{k!} (d^k F)_0(X^{\otimes k}) \\ 
    \Rem_{x_0}^n \psi(x) &= \dircon(x,x_0) \frac{1}{n!} \int_0^1 (1-t)^n (d^{n+1} F)_{tX}(X^{\otimes n+1}). 
\end{align*}
\end{defn}

\begin{lem}[Covariant Taylor expansion] \label{lem:Taylor for F}
With the above notations, for any $n\geq 0$ we have 
\begin{align*}
\psi(x) &= \Tay_{x_0}^n \psi(x) + \Rem_{x_0}^n \psi(x)
\end{align*}
for any $x\in U_{x_0}$, and $\Rem_{x_0}^n \psi(x)$ converges to zero as $x$ tends to $x_0$. 
\end{lem}

\begin{proof}
Since the function $F:V_{x_0}\to E_{x_0}$ given in Definition \ref{def:function F} is vector valued, we can apply the usual Taylor theorem with integral remainder, see for instance \cite[Theorem 2.4.15]{Abraham-Marsden-Ratiu:1988}, which says that 
\begin{align*}
F(X) &= \sum_{k=0}^n \frac{1}{k!}\,(d^k F)_0(X^{\otimes k})
    + \frac{1}{n!} \int_0^1 (1-t)^n\,(d^{n+1}F)_{tX}\big(X^{\otimes n+1}\big)\,dt. 
\end{align*} 
Applying $\dircon(x,x_0)$ on both sides yields the Taylor expansion of $\psi$ ot $x_0$. 
\end{proof}

Let us compare the above approximation of $\psi$ produced by $\dircon$ to the usual approximation in terms of the higher covariant derivatives of the linear connection on $E$ induced by $\dircon$. 
We show that the difference between the two approximations is encoded in the derivatives of the curvature $K$ of $\dircon$, given in Definition \ref{def:curvature}, along geodesics.

\begin{defn} \label{def:transitivity defects}
Fix $x_0\in M$. We call \textbf{transitivity defects} of $\dircon$ in the direction $X\in V_{x_0}$, the derivatives of the curvature operator $K_{x_0}$ along the geodesic $x(t)=\exp_{x_0}(tX)$, namely the operators (implicitly depending on $X$)
\begin{align} \nonumber
    \TDef_1(t) &:= \frac{d}{ds} K_{x_0}(x(s),x(t))\big|_{s=t} 
    \\ \label{eq:Cn(t)}
    \TDef_2(t) &:= \frac{d}{du} \TDef_1(u)\big|_{u=t} = \frac{\partial^2}{\partial u^2} K_{x_0}(x(u),x(u))\big|_{u=t} + \frac{\partial^2}{\partial u \partial s} K_{x_0}(x(s),x(u))\big|_{s=u=t} 
    \\ \nonumber
    \TDef_n(t) &:= (\TDef_1)^{(n-1)}(t),\qquad n\geq 2
\end{align}
\end{defn}

If $\dircon$ is the parallel transport of a linear connection, then the curvature $K_{x_0}$ is the identity along geodesics and the transitivity defects $\TDef_n(t)$ vanish for any $t\in [0,1]$ and for any $n\geq 1$. 
For a general direct connection, it follows from Lemma \ref{lem:d/dt Gamma} that $\TDef_1(0)=0$, but the transitivity defects for $n>1$ do not necessarily vanish even at $t=0$ (above $x_0$). 

\begin{example}
Consider the general direct connection $\dircon(y,x)=e^{\Lambda(y)-\Lambda(x)+(y-x)^2h(y,x)}$ of Example \ref{ex:dircon on RxR} in the case $h(y,x)=c$ is constant, with curvature $K_x(z,y)=e^{-2c(z-y)(y-x)}$ given in Example \ref{ex:curvature on RxR}. 
Taking $x=x_0$, $y=x(u)=x_0+uX$ and $z=x(s)=x_0+sX$, equations \eqref{eq:Cn(t)} give 
\begin{align*}
    \TDef_1(t)= -2ctX^2, \qquad \TDef_2(t)=-2cX^2 \qquad\text{and}\qquad \TDef_n(t)=0 \quad \forall n\geq 3,  
\end{align*}
therefore $\TDef_2(0) \neq 0$. 
\end{example}

To prove our statement, it is convenient to introduce an operator that records the discrepancy between the direct transport $\dircon(x,x_0)$ and the parallel transport $\Par^\nabla(x,x_0)$ of the linear connection $\nabla$ induced by $\dircon$.

\begin{defn} \label{rem:parallel defect} 
Fix $x_0\in M$. Let $\nabla:=\nabla^\dircon$ be the linear connection on $E$ induced by $\dircon$, and let $\Par^\nabla$ be its associated parallel transport on $E$. 
We call \textbf{parallel defect} in the direction $X\in V_{x_0}$ the operator
$$
\PDef(t) := \dircon(x(t),x_0)^{-1} \Par^\nabla(x(t),x_0)\in GL(E_{x_0}), 
$$
defined along the geodesic $x(t)=\exp_{x_0}(tX)$. 
\end{defn}

\begin{lem} \label{lem:parallel defects}
The parallel defect $\PDef(t)$ satisfies the differential equation with initial condition 
\begin{align*}
    \frac{d}{dt}\PDef(t) &= \TDef_1(t)\, \PDef(t) \\ 
    \PDef(0) &=1_{x_0}. 
\end{align*}
Therefore its derivatives are given by the closed formula 
\begin{align*}
    \frac{d^m}{dt^m} \PDef(t) &= \sum_{r=1}^m \underset{k_1+\cdots +k_r=m}{\sum_{k_1,...,k_r\geq 1}} \prod_{i=1}^r \binom{k_i+\cdots +k_r-1}{k_i-1}\, \TDef_{k_1}(t)\cdots \TDef_{k_r}(t)\, \PDef(t), 
\end{align*}
where the sum runs over ordered compositions of $m$. 
\end{lem}

\begin{proof}
The first derivative of $\PDef$ at $t$ is
$$
\frac{d}{ds} \PDef(s)\big|_{s=t} = \frac{d}{ds} \dircon(x(s),x_0)^{-1} \Par^\nabla(x(s),x_0)\big|_{s=t}. 
$$
To compute it, fix any vector $e_0\in E_{x_0}$ and consider the parallel transported section 
$$
\psi(t) = \Par^\nabla(x(t),x_0)e_0 \in E_{x(t)},
$$
which gives 
$$
\PDef(t) e_0 = \dircon(x(t),x_0)^{-1} \psi(t), 
$$
and which is parallel with respect to $\nabla$ along $x(t)$, that is, 
$$
\nabla \psi\big|_{x(t)}(\dot{x}(t)) = 0,\qquad \forall t\in [0,1].
$$
Using the curvature identity 
$$
K_{x_0}(x(s),x(t)) = \dircon(x(s),x_0)^{-1} \dircon(x(s),x(t)) \dircon(x(t),x_0), 
$$
we can write
$$
\PDef(s) e_0 = \dircon(x(s),x_0)^{-1} \psi(s) = K_{x_0}(x(s),x(t)) \dircon(x(t),x_0)^{-1} \dircon(x(s),x(t))^{-1} \psi(s),
$$
and therefore 
\begin{align*}
    \frac{d}{ds} \big(\PDef(s) e_0\big)\big|_{s=t} &= \frac{d}{ds}K_{x_0}(x(s),x(t))\big|_{s=t} \, \dircon(x(t),x_0)^{-1} \psi(t)
    \\
    & \hspace{2cm} + \dircon(x(t),x_0)^{-1}\, \frac{d}{dt}\dircon(x(s),x(t))^{-1} \psi(s)\big|_{s=t} \\ 
    &= \TDef_1(t)\,\PDef(t) e_0 + \dircon(x(t),x_0)^{-1}\, \nabla \psi\big|_t(X) \\ 
    &= \TDef_1(t)\,\PDef(t) e_0, 
\end{align*}
because $\psi$ is parallel. 
The initial condition $\PDef(0)=1_{x_0}$ follows from the fact that both $\dircon$ and $\Par^\nabla$ are the identity operators at $t=0$. 

The closed formula is proved by repeatedly  applying the Leibniz rule using the same trick. This is a standard computation for derivatives of noncommutative polynomials, see for instance \cite[Remark 4.3]{LundervoldMuntheKaas:2011}, based on ordered compositions and noncommutative Bell polynomials. 
\end{proof}

\begin{example}[Derivatives of the parallel defect] \label{ex:parallel defect}
The first 5 derivatives of the parallel defect at time $t$ (omitting $t$ in the transitive defects) are
\begin{align*}
    \PDef'(t) &= \TDef_1\, \PDef(t) \\
    \PDef''(t) &= \Big(\TDef_2 + \TDef_1^2\Big)\, \PDef(t) \\
    \PDef^{(3)}(t) &= \Big(\TDef_3 + (2 \TDef_2 \TDef_1 + \TDef_1 \TDef_2) + \TDef_1^2\Big)\, \PDef(t) \\ 
    \PDef^{(4)}(t) &= \Big(\TDef_4 + (3 \TDef_3 \TDef_1 + 3 \TDef_2^2 + \TDef_1 \TDef_3) + (3 \TDef_2 \TDef_1^2 + 2 \TDef_1 \TDef_2 \TDef_1 + \TDef_1 \TDef_2^2) + \TDef_1^4\Big)\, \PDef(t) \\ 
    \PDef^{(5)}(t) &= \Big(\TDef_5 + (4 \TDef_4 \TDef_1 + 6 \TDef_3 \TDef_2 + 4 \TDef_2 \TDef_3 + \TDef_1 \TDef_4) \\ 
    & \hspace{1cm} 
    + (6 \TDef_3 \TDef_1^2 + 8 \TDef_2^2 \TDef_1 + 4 \TDef_2 \TDef_1 \TDef_2 + 3 \TDef_1 \TDef_3 \TDef_1 + 3 \TDef_1 \TDef_2^2 + \TDef_1^2 \TDef_3) \\ 
    & \hspace{1cm} 
    + ( 4 \TDef_2 \TDef_1^3 + 3 \TDef_1 \TDef_2 \TDef_1^2 + 2 \TDef_1 \TDef_2^2 \TDef_1 + \TDef_1^3 \TDef_2) + \TDef_1^5 \Big)\, \PDef(t) .
\end{align*}
Notice that, as a result of the Leibniz formula, $b^{(k)}(t)$ only involves terms $a_{i_1}^{\alpha_1}\cdots a_{i_j}^{\alpha_j}$ of weight $\sum_{\ell=1}^j \ell\cdot \alpha_\ell= k$.
\end{example}

Let $\nabla:=\nabla^\dircon$ be the linear connection on $E$ induced by $\dircon$, as above. Denote by $\nabla^n:\Gamma(E)\to T^*M^{\otimes n}\otimes E$ the higher covariant derivative induced by $\nabla^M$ and $\nabla$, as in \eqref{eq:higher covariant derivatives}, and by $D_\nabla^n:\Gamma(E)\to S^n(T^*M)\otimes E$ its associated symmetrised covariant derivatives, as in \eqref{eq:symmetrised covariant derivatives}. 

\begin{thm} \label{thm:Taylor with nabla}
Fix $x_0\in M$. 
In terms of the higher covariant derivatives of $\nabla$, the Taylor polynomial of $\psi$ at $x_0$ of any order $n\geq 0$, and its remainder, are given for any $x=\exp_{x_0}(X)\in U_{x_0}$ by 
\begin{align*}
\Tay_{x_0}^n \psi(x) 
&= \dircon(x,x_0) \sum_{k=0}^n \frac{1}{k!}\, \Coef_{n,k}(x_0) \, D_\nabla^k \psi\big|_{x_0}\big(X^{\otimes k}\big),
\end{align*} 
and 
\begin{align*}
\Rem_{x_0}^n \psi(x) 
&= \dircon(x,x_0) \sum_{k=0}^{n+1} \frac{1}{n!} \binom{n+1}{k} \int_0^1 (1-t)^n \PDef^{(n+1-k)}(t)\, \Par^\nabla(x(t),x_0)^{-1}\, D_\nabla^k \psi\big|_{x(t)}\big(\dot{x}(t)^{\otimes k}\big) dt,
\end{align*} 
where the coefficient is the operator on $E_{x_0}$ given by 
\begin{align} \label{eq:Taylor correction}
    \Coef_{n,k} (x_0) &:= 1_{x_0}+\sum_{m=2}^{n-k} \frac{1}{m!}\, \frac{d^m}{dt^m}\PDef(0). 
\end{align}
\end{thm}

Note that the symmetrised covariant derivatives $D_\nabla^k \psi$ are evaluated on the same tangent vector $k$ times, therefore they coincide with the higher covariant derivatives $\nabla^k \psi$. 

\begin{proof}
Consider the section $\phi(t) = \Par^\nabla(x(t),x_0)^{-1} \psi(x(t))$ parallel transported along the geodesic $x(t)$, 
so that 
\begin{align*}
    \dircon(x(t),x_0)^{-1} \psi(x(t)) &= \dircon(x(t),x_0)^{-1} \Par^\nabla(x(t),x_0)\, \phi(t) = \PDef(t)\, \phi(t). 
\end{align*}
By definition of the parallel transport with respect to the linear connection $\nabla$, we have 
\begin{align*}
    \frac{d}{dt} \Big(\Par^\nabla(x(t),x_0)^{-1} \psi(x(t))\Big) = \Par^\nabla(x(t),x_0)^{-1} \nabla_{\dot{x}(t)}\psi\big|_{x(t)}.
\end{align*}
Iterating, we get 
\begin{align} \nonumber
    \phi^{(k)}(t) &= \frac{d^k}{dt^k} \Par^\nabla(x(t),x_0)^{-1} \psi(x(t)) 
    \\ \nonumber
    &= \Par^\nabla(x(t),x_0)^{-1} \underbrace{\nabla_{\dot{x}(t)}(\nabla_{\dot{x}(t)}\cdots (\nabla_{\dot{x}(t)} \psi)\cdots )}_{\text{$k$ times}}\big|_{x(t)} 
    \\ \label{eq:parallel transported section phi} 
    &= \Par^\nabla(x(t),x_0)^{-1} \nabla^k \psi\big|_{x(t)}(\dot{x}(t)^{\otimes k}), 
\end{align}
where the last equality holds because the covariant derivative is performed along a geodesic, cf. Appendix \ref{subsec:jets and covariant derivatives}, eq. \eqref{eq:nabla^n}. 
In particular, at $t=0$ we have 
$$
\phi^{(k)}(0) = \nabla^k \psi\big|_{x_0}(X^{\otimes k}). 
$$

To prove the assertion for the Taylor polynomial, we then apply the Leibniz rule repeatedly at $t=0$, and obtain 
\begin{align*}
    \sum_{j=0}^n \frac{1}{j!} \frac{d^j}{dt^j} \dircon(x(t),x_0)^{-1} \psi(x(t))\big|_{t=0} 
    &= \sum_{k=0}^n \sum_{m=0}^{n-k} \frac{1}{m!} \PDef^{(m)}(0)\, \frac{1}{k!} \phi^{(k)}(0) \\ 
    &= \sum_{k=0}^n \frac{1}{k!} \Coef_{n,k}(x_0) \, \nabla^k \psi\big|_{x_0}(X^{\otimes k}), 
\end{align*}
where 
$$
\Coef_{n,k}(x_0) = \sum_{m=0}^{n-k} \frac{1}{m!}\, \PDef^{(m)}(0). 
$$
This coefficient agrees with \eqref{eq:Taylor correction}, because $\PDef(0)=1_{x_0}$ and $\PDef'(0)=0$ since $\TDef_1(0)=0$. 

For the remainder, the result follows by applying the Leibniz rule up to order $n+1$ at a generic time $t$, and replacing the derivatives $\phi^{(k)}(t)$ by their expression in eq. \eqref{eq:parallel transported section phi}. 
\end{proof}

\begin{remark}[Comparing Taylor expansions at different orders]
\label{rem:Taylor at different orders}
Unlike in the usual Taylor expansion, the Taylor polynomial at order $n$ is not just that at order $n-1$ plus additional terms: the whole approximation is adjusted differently at every lower order, because the coefficients $\Coef_{n,k}(x_0)$ depend on $n$ and not only on $k$. 
However, one still has
\begin{align*}
    \Tay_{x_0}^{n}\psi(x)-\Tay_{x_0}^{n-1}\psi(x) =
\frac{1}{n!} \Gamma(x,x_0)\, (d^n F)_0(X^{\otimes n}), 
\end{align*}
for $F(tX)$ as in Definition \ref{def:function F}. 
So, while the change from order $n-1$ to order $n$ in covariant derivatives involves several lower derivative orders, it 
gives one homogeneous term of degree $n$ in the $\dircon$ trivialisation.
\end{remark}

The explicit expression of the Taylor polynomial and of the remainder in terms of the transitivity defects is then obtained by replacing $\PDef^{(n+1-k)}(t)$ with its closed expression given in Lemma \ref{lem:parallel defects}. The resulting formulas are very long and are not worth writing out. Instead, we give explicit examples. 

\begin{example}[Lower order corrections for the Taylor polynomial] \label{ex:example Taylor}
Using the derivatives of the parallel defect of Example \ref{ex:parallel defect}, evaluated at $t=0$ (with $\TDef_1(0)=0$), and omitting the dependence of $\TDef_k$ on $t=0$, we obtain the Taylor expansion at order $n=5$: 
\begin{align*}
\Tay_{x_0}^5 f(x) &= \dircon(x,x_0) \left[\left( 1_{x_0} + \frac{1}{2!}\, \TDef_2 + \frac{1}{3!}\, \TDef_3 + \frac{1}{4!}\, (\TDef_4 + 3\,\TDef_2^2) + \frac{1}{5!}\, (\TDef_5+ 6\,\TDef_3 \TDef_2 + 4\, \TDef_2 \TDef_3)\right) f(x_0)\right. 
\\ 
&\hspace{2cm}
+ \left( 1_{x_0} + \frac{1}{2!}\, \TDef_2 + \frac{1}{3!}\, \TDef_3  + \frac{1}{4!}\, (\TDef_4 + 3\,\TDef_2^2)\right) \nabla f\big|_{x_0}(X) 
\\
&\hspace{2cm}
+ \frac{1}{2!}\, \left( 1_{x_0} + \frac{1}{2!}\, \TDef_2 + \frac{1}{3!}\, \TDef_3\right) \nabla^2 f\big|_{x_0}(X^{\otimes 2}) 
\\
&\hspace{2cm}
+ \frac{1}{3!}\, \left( 1_{x_0} + \frac{1}{2!}\, \TDef_2\right) \nabla^3 f\big|_{x_0}(X^{\otimes 3}) 
\\
&\hspace{2cm}
+ \frac{1}{4!}\, \nabla^4 f\big|_{x_0}(X^{\otimes 4}) 
\\
&\hspace{2cm}
\left. + \frac{1}{5!}\, \nabla^5 f\big|_{x_0}(X^{\otimes 5}) \right]. 
\end{align*}
\end{example}

\begin{example}[Approximation of sections] \label{ex:non-local approximation}
Let $\psi:M\to E$ be a smooth section, $x_0\in M$ a given point and $U_{x_0}$ a convex neighbourhood. 
The expansion of order $0$ of $\psi$ at $x_0$ induced by a direct connection $\dircon$ on $GL(E)$ is given, at $x\in U_{x_0}$, by
\begin{align*}
    \psi(x) 
    &= \dircon(x,x_0)\, \psi(x_0)
    \\ 
    & \quad + \dircon(x,x_0)\, \int_0^1 \Big[\TDef_1(t)\,\dircon(x(t),x_0)^{-1}\,\psi(x(t))+\dircon(x(t),x_0)^{-1} \nabla \psi\big|_{x(t)}(\dot{x}(t))\Big]\ dt.
\end{align*}
The expansion of order $1$  is given by
\begin{align*}
    \psi(x) 
    &= \dircon(x,x_0)\, \psi(x_0) + \dircon(x,x_0)\, \nabla \psi\big|_{x_0}(X) \\ 
    &\quad + \dircon(x,x_0)\, \int_0^1 (1-t)\,
    \Big[ (\TDef_2(t)+\TDef_1(t)^2)\,\dircon(x(t),x_0)^{-1}\,\psi(x(t)) 
    \\ 
    &\quad \hspace{2cm}
    + 2\,\TDef_1(t)\,\dircon(x(t),x_0)^{-1}\,\nabla \psi\big|_{x(t)}(\dot{x}(t)) + \dircon(x(t),x_0)^{-1} \nabla^2 \psi\big|_{x(t)}\big(\dot{x}(t)^{\otimes 2}\big) \Big]\ dt.
\end{align*}
The expansion of order $2$ reads
\begin{align*}
    \psi(x) 
    &= \dircon(x,x_0) \Big(1_{x_0}+\frac{1}{2}\, \TDef_2(0)\Big) \psi(x_0) + \dircon(x,x_0)\, \nabla \psi\big|_{x_0}(X) + \frac{1}{2}\dircon(x,x_0)\, \nabla^2 \psi\big|_{x_0}(X^{\otimes 2}) \\ 
    &\quad + \dircon(x,x_0)\, \int_0^1 (1-t)^2\,
    \Big[ \frac{1}{2}\,\big(\TDef_3(t)+2\TDef_2(t)\TDef_1(t)+\TDef_1(t)\TDef_2(t)+\TDef_1(t)^3\big)\,\dircon(x(t),x_0)^{-1}\,\psi(x(t)) 
    \\ 
    &\quad \hspace{2cm}
    + \frac{3}{2}\,\big(\TDef_2(t)+\TDef_1(t)^2\big)\,\dircon(x(t),x_0)^{-1}\,\nabla \psi\big|_{x(t)}(\dot{x}(t)) 
    \\ 
    &\quad \hspace{2cm}
    + \frac{3}{2}\,\TDef_1(t)\,\dircon(x(t),x_0)^{-1} \nabla^2 \psi\big|_{x(t)}\big(\dot{x}(t)^{\otimes 2}
    + \frac{1}{2} \dircon(x(t),x_0)^{-1}\,\nabla^3 \psi\big|_{x(t)}\big(\dot{x}(t)^{\otimes 3}\big) \Big]\ dt.
\end{align*}
Let us note that, in the principal part of the second-order expansion, the usual constant term $\dircon(x,x_0)\psi(x_0)$ involves a new factor $1_{x_0}+\frac{1}{2}\, \TDef_2(0)$ that depends on a transitivity defect, in agreement with Example \ref{ex:example Taylor}. 
\end{example}

\begin{remark}
Suppose $\dircon$ is a transitive direct connection, which is equivalent to the parallel transport of a linear connection $\nabla$ on $E$ by Theorem \ref{thm:transitive = parallel transport}. Then all the correction terms $\TDef_n(t)$ vanish, and all coefficients $\Coef_{n,k}(x_0)$ are identically $1_{x_0}$. The Taylor expansion then coincides with that given in \cite[Definition 65]{DahlqvistDiehl-Driver:2019} for the trivial bundle $E=M\times \R$ (where the displacement from $E_{x_0}$ to $E_x$ can be omitted, since it is induced by the global trivialisation of the bundle), and the covariant Taylor expansion of Lemma \ref{lem:Taylor for F} coincides with that in \cite[Theorem 62]{DahlqvistDiehl-Driver:2019}. 
\end{remark}

\begin{prop} \label{prop:remainder}
With the above notation, for any $x_0\in M$, for any compact subset $K\subset U_{x_0}$ of a given convex neighbourhood of $x_0$ and for any $x\in K$, we have
\begin{align*}
    \Rem_{x_0}^n \psi(x)= \calO\big(\|X\|^{\otimes n+1}\big).
\end{align*}
\end{prop}

\begin{proof}
For $x_0$ fixed, let $F:V_{x_0}\to E_{x_0}$ be the function given in Definition \ref{def:function F} which satisfies $\dircon(x(t),x_0)^{-1} \psi(x(t)) = F(tX)$ for any $t\in [0,1]$, and consequently 
\begin{align*}
    \frac{d^m}{dt^m} \dircon(x(t),x_0)^{-1} \psi(x(t)) &= (d^m F)_{tX}(X^{\otimes m}),\qquad \forall m\geq 0. 
\end{align*}
As in the proof of Lemma \ref{lem:Taylor for F}, we then have 
\begin{align*}
    \Rem_{x_0}^n \psi(x) 
    &= \dircon(x,x_0) \frac{1}{n!} \int_0^1 (1-t)^n \, (d^{n+1} F)_{tX}(X^{\otimes n+1})\, dt. 
\end{align*}
On a sufficiently small compact normal neighborhood $K$ of $x_0$, both the operator $\dircon(x,x_0)$ and $(d^{n+1} F)_{tX}(X^{\otimes n+1})$ are bounded, for $t\in [0,1]$. Hence there exists a constant $\TDef_{\dircon,\psi,K,n}>0$ such that 
\begin{align*}
    \|\Rem_{x_0}^n \psi(x)\| \leq \TDef_{\dircon,\psi,K,n} \|X\|^{n+1}. 
\end{align*}
\end{proof}

 
\subsection{Higher-order direct connections and jet regularity structures} 
\label{subsec:Higher-order dircon} 

The jet bundle $J^nE$ contains   higher covariant derivatives of the sections of $E$ and is therefore   of special interest in the context of regularity structures \cite{Hairer:2014b,Hairer-Singh:2023}. Inspired by \cite[Definition 80]{DahlqvistDiehl-Driver:2019}, in this section we give an example of a higher-order direct connection on $J^nE$ which is not a jet prolongation of a direct connection on $E$. 
We then show that it can be used to define a covariant (and local) model for an $E$-valued jet regularity structure. 

\begin{defn}\label{def:Taylor dircon}
Let $\dircon$ be a direct connection on $\GL(E)$ defined on a convex diagonal domain $\calU$, and for any $x_0\in M$ let $U_{x_0}$ be a convex open neighborhood of $x_0$ such that $U_{x_0}\times U_{x_0}\subset \calU$. 
For any $n\geq 0$, we call \textbf{Taylor direct connection induced by $\dircon$} the collection of linear maps $\Sigma^n(y_0, x_0):J^n_{x_0}E\to J^n_{y_0}E$ defined for any $x_0\in M$ and any $y_0\in U_{x_0}$ by
\begin{align*} 
\Sigma^n(y_0,x_0)\ j_{x_0}^n \psi &:= j_{y_0}^n (\Tay^n_{x_0}\psi). 
\end{align*}
We  recall from \eqref{eq:Taylor} that $\Tay^n_{x_0}\psi(x)=\dircon(x,x_0) \sum_{k=0}^n \frac{1}{k!} d^k F_0(X^{\otimes k})$ is a local section of $E$ defined for any $x\in U_{x_0}$, where $F(tX)=\dircon(x(t),x_0)^{-1} \psi(x(t))$. 
\end{defn}

\begin{example} \label{ex:Taylor direct connection}
For $n=0$, we recover the value of $\psi$ transported by $\dircon$ from $x_0$ to $y_0$: 
\begin{align*}
\Sigma^0(y_0,x_0) \psi(x_0) 
&= \dircon(y_0,x_0) \psi(x_0) .
\end{align*}

For $n=1$, we start with the jet $j^1_{x_0} \psi = \big[ x_0,\psi(x_0),\nabla \psi\big|_{x_0}(X) \big]$ for any $X\in T_{x_0}M$, and wish to describe 
\begin{align*}
    \Sigma^1(y_0,x_0) j^1_{x_0} \psi &= j^1_{y_0} (\Tay^1_{x_0}\psi) = \left[ y_0,\Tay^1_{x_0}\psi(y_0), \nabla (\Tay^1_{x_0}\psi)\big|_{y_0}(Y) \right]  
\end{align*}
for any $Y\in T_{y_0}M$. The $0$-th order component of this jet is given by
\begin{align*}
\Tay^1_{x_0}\psi(y_0) &= \dircon(y_0,x_0) \left( \psi(x_0) + \nabla\psi\big|_{x_0}(X) \right), 
\end{align*}
where $X\in T_{x_0}M$ is the tangent vector such that $y_0=\exp_{x_0}(X)\in U_{x_0}$. Call $x(t) = \exp_{x_0}(tX)$ the geodesic from $x_0$ to $y_0$, with $t\in [0,1]$. 
To compute the first order component of the jet, take a generic vector $Y\in T_{y_0}M$ such that the geodesic $y(s)=\exp_{y_0}(sY)\in U_{y_0}\cap U_{x_0}$. Then 
\begin{align*}
    \nabla (\Tay^1_{x_0}\psi)\big|_{y_0}(Y) &= \partial_s \big[ \dircon(y(s),y_0)^{-1} \Tay^1_{x_0}\psi (y(s))\big]\big|_{s=0} \\ 
    &= \partial_s \big[ \dircon(y(s),y_0)^{-1} \dircon(y(s),x_0) \big(\psi(x_0) + \nabla \psi\big|_{x_0}(X[s])\big)\big]\big|_{s=0}.
\end{align*}
This last expression involves on the one hand the composition
\begin{align*}
    \dircon(y(s),y_0)^{-1} \dircon(y(s),x_0) &= \big(\dircon(y(s),x_0)^{-1} \dircon(y(s),y_0) \big)^{-1} 
    = \big(K_{x_0}(y(s),y_0) \dircon(y_0,x_0)^{-1} \big)^{-1} \\ 
    &= \dircon(y_0,x_0) K_{x_0}(y(s),y_0)^{-1}. 
\end{align*}
On the other hand it involves the vector $X[s]$, which is the velocity of the geodesic $x(t,s)=\exp_{x_0}(tX[s])$ from $x_0$ to $y(s)$. This is the variation curve of $x(t)$ through geodesics called a \emph{spray}, depicted in Fig \ref{fig:inverse problem} (cf. Fig. \ref{fig:variation gamma}). 
\begin{figure}[h!]
\centering
\begin{tikzpicture}[scale=.7]
\draw (0,0) node {$\bullet$}; 
\draw (-0.3,-0.1) node {$x_0$};
\draw[thick] (0,0) 
      .. controls (3,-0.4) and (5,-0.4) .. (8,0);
\draw[thick,red,->] (0,0) -- (-1,5);
\draw[red] (-1.3,5.4) node {$\widetilde{Y}=\Exp(y_0,x_0)^{-1}(Y)$}; 
\draw (4,-0.3) node {$\bullet$};
\draw (4,-0.8) node {$x(t)$};
\draw[thick,blue] (4,-0.3)
      .. controls (3.8,0.7) and (3.5,1.4) .. (3,2);
\draw[thick,blue,->] (4,-0.3) -- (3.5,2.2); 
\draw[blue] (4.3,1.9) node {$J(t)$};
\draw (8,0) node {$\bullet$}; 
\draw (8.3,-0.1) node {$y_0$};
\draw[thick,blue] (8,0)
      .. controls (7.7,1) and (7.6,2) .. (6,4.8);
\draw[thick,blue,->] (8,0) -- (7,5.1); 
\draw[blue] (7.8,4.7) node {$J(1)$};
\draw (9,4.7) node {$=\!Y$};
\draw[thick,->] (0,0) -- (8,-1); 
\draw (8.3,-1) node {$X$};
\draw[thick] (0,0)
      .. controls (3,0.6) and (5,1.2) .. (7.2,2.4);
\draw[thick,->] (0,0) -- (7.5,1.5); 
\draw (6.8,0.9) node {$X[s]=X+s\widetilde{Y}$};
\draw[thick] (0,0)
      .. controls (3,1.7) and (5,3.4) .. (6,4.8);
\draw[thick,->] (0,0) -- (7,4); 
\draw (6.8,3.2) node {$X+Z$};
\draw (3.6,0.9) node {$\bullet$}; 
\draw (2.7,1) node {$x(t,s)$};
\draw (7.2,2.4) node {$\bullet$};
\draw (8.2,2.5) node {$y(s)$};
\draw (6,4.8) node {$\bullet$};
\draw (6.3,5.2) node {$y$};
\end{tikzpicture}
\caption{Spray of variation curves $x(s,t)$ for the inverse problem $J(1)=Y$.}
\label{fig:inverse problem}
\end{figure}
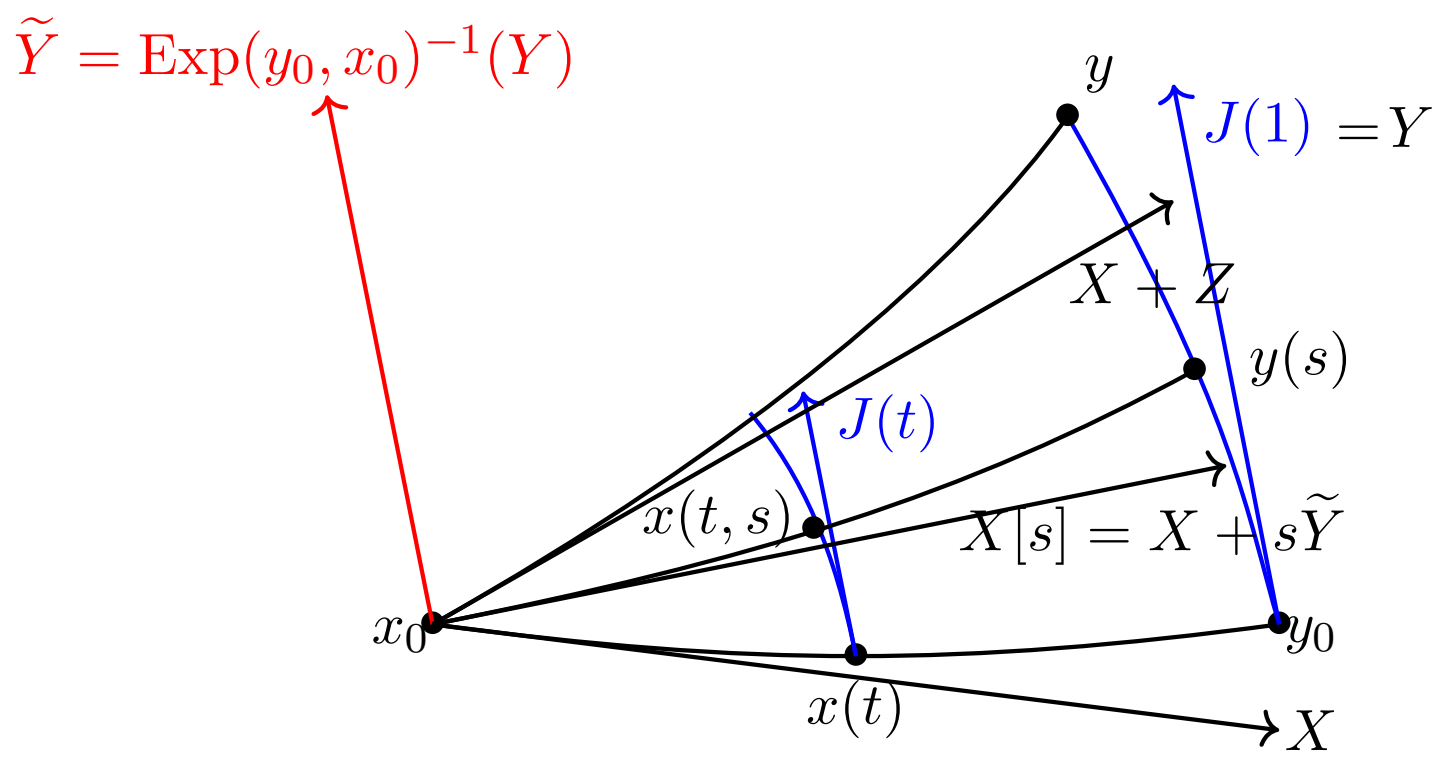 
It determines a Jacobi field $J(t)$ along $x(t)$ with initial conditions $J(0)=0$ and $D_tJ(t)\big|_{t=0}=\widetilde{Y}$ such that the initial velocity of the geodesic $t\to x(t,s)$ is $X[s] = X+s\widetilde{Y}$. In this case, the initial condition $\widetilde{Y}$ is not given, and the Jacobi field is uniquely determined by the final condition $J(1)=Y$. However, Theorem \ref{thm:exponential direct connection} says that $\widetilde{Y}=\Exp^1(y_0,x_0)^{-1}(Y)$, where $\Exp^1(y_0,x_0)$ is the exponential direct connection on $J^1\Pair(M)\cong GL(TM)$ given in Definition \ref{def:Exp^n}. 
We then have
\begin{align*}
    \partial_s(X[s])\big|_{s=0} &= \widetilde{Y} = \Exp^1(y_0,x_0)^{-1}(Y),
\end{align*}
and finally 
\begin{align} \label{eq:nabla Tay^1}
    \nabla (\Tay^1_{x_0}\psi)\big|_{y_0}(Y) &= \dircon(y_0,x_0) \partial_s K_{x_0}(y(s),y_0)^{-1}\big|_{s=0} \left(\psi(x_0) + \nabla\psi\big|_{x_0}(X) \right) \\ \nonumber
    & \hspace{2cm} 
    + \dircon(y_0,x_0) \nabla \psi\big|_{x_0}(\Exp^1(y_0,x_0)^{-1}(Y)).  
\end{align}
\end{example}

\begin{thm} \label{thm:Taylor direct connection}
The collection of maps $\Sigma^n(y_0,x_0)$ defines a direct connection on $\GL(J^nE)$ which, for $n\geq 1$, is not a jet prolongation of a direct connection on $E$. 
\end{thm}

\begin{proof}
For the first assertion, by Proposition \ref{prop:paralleldirconn} it is enough to show that $\Sigma$ is a parallelism. 
The map $\Sigma^n(y_0,x_0)$ depends smoothly on $(y_0,x_0)$, since its construction involves the smooth maps $\dircon$, the local exponential map, their inverses and their derivatives of finite order. It is also clearly a linear operator from $J^n_{x_0}E$ to $J^n_{y_0}E$, because the operator $\Tay_{x_0}^n$ is linear in the derivatives of $\psi$. 
Then we only have to show that $\Sigma^n(x_0,x_0)=\Id_{J^n_{x_0}E}$, that is, 
\begin{align*} 
    \Sigma^n(x_0,x_0)\ j_{x_0}^n \psi &= j_{x_0}^n \psi, 
\end{align*}
or equivalently, by definition of $\Sigma^n$, that
\begin{align} \label{eq:JT=J}
    j_{x_0}^n (\Tay_{x_0}^n\psi) &= j_{x_0}^n \psi. 
\end{align}

This holds trivially for $n=0$. For $n=1$ it follows from Example \ref{ex:Taylor direct connection}. 
In fact, setting $y_0=x_0$ in eq. \eqref{eq:nabla Tay^1}, we see that the first term vanishes because $K_{x_0}(y(s),x_0)=1_{x_0}$ and, since $Y=X$, the second term gives 
\begin{align*}
    \nabla (\Tay^1_{x_0}\psi)\big|_{x_0}(X) &= \nabla \psi\big|_{x_0}(X).  
\end{align*}
Since $\Tay^1_{x_0}\psi(x_0)=\psi(x_0)$, this means that $j^1_{x_0}(\Tay^1_{x_0}\psi) = j^1_{x_0}\psi$. 

More generally, for any $n\geq 0$, fix $x_0$ and let $F:V_{x_0}\to E_{x_0}$ be the function $F(X)=\dircon(\exp_{x_0}(X),x_0)^{-1} \psi(\exp_{x_0}(X))$ as in Definition \ref{def:function F}. Then, by definition of the Taylor section, we have
\begin{align*}
    \Tay_{x_0}^n \psi (x) &= \dircon(x,x_0)\, T_0^n F(X),\qquad x=\exp_{x_0}(X), 
\end{align*}
where 
\begin{align*}
    T_0^n F(X) &:= \sum_{k=0}^n \frac{1}{k!} (d^kF)_0(X^{\otimes k}),
\end{align*}
is the usual Taylor polynomial of $F$, as in the proof of Lemma \ref{lem:Taylor for F}. 
The identity \eqref{eq:JT=J} we want to prove then says that the jet of order $n$ at $x_0$ of the section 
\begin{align*}
    \psi - \Tay_{x_0}^n \psi &= \dircon(x,x_0)\, (F - T_0^n F)
\end{align*}
vanishes. 
Now, the map 
\begin{align*} 
    L_{x_0}^n:J^n_{x_0}E \to J^n_0(V_{x_0},E_{x_0}),\ j_{x_0}^n \psi \mapsto L_{x_0}^n(j_{x_0}^n \psi) := j_0^n F
\end{align*}
is a linear isomorphism of jet spaces, with inverse 
\begin{align*}
    (L_{x_0}^n)^{-1}(j_0^n F) &= j_{x_0}^n\big(x\mapsto \dircon(x,x_0) F(\exp_{x_0}^{-1}(x)\big). 
\end{align*}
We apply it to the section $\psi-\Tay_{x_0}^n\psi$, and get 
\begin{align} \nonumber
    L_{x_0}^n\Big(j_{x_0}^n \psi - j_{x_0}^n \Tay_{x_0}^n \psi\Big) 
    & = L_{x_0}^n j_{x_0}^n \Big(\psi - \Tay_{x_0}^n \psi\Big) 
    \\ \nonumber 
    & = j_0^n \Big(F - T_0^n F \Big) 
    \\ \label{eq:jF-jTF}
    & = j_0^n F - j_0^n(T_0^n F). 
\end{align} 
Since $F$ is a map between vector spaces, the Taylor polynomial of $F$ is naturally defined (once the origin is fixed), and the standard Taylor theorem for $F$ at $0$ ensures that the difference \eqref{eq:jF-jTF} vanishes. 
Since $L_{x_0}^n$ is an isomorphism, this proves eq. \eqref{eq:JT=J}. 

Finally, to show that $\Sigma^n$ is not the jet prolongation of a direct connection on $E$, for $n\geq 1$, it suffices to show that it does not take value in the subgroupoid $J^n\GL(E)$ of $GL(J^nE)$. The reason is that $\Sigma^n j^n_{x_0}\psi$ is not filtered in the order of jets of $\psi$, that is, the invertible matrix associated to $\Sigma^n$ in the linear basis of $J^n_{x_0}E$ spanning jets $j^n_{x_0}\psi$ by the covariant derivatives $\nabla^k \psi\big|_{x_0}$ is not a triangular matrix. This can be seen in $\Sigma^1(y_0,x_0)$, since the $0$-th order component of the jet $j^1_{y_0}(\Tay^1_{x_0}\psi)$ is 
\begin{align*}
\Tay^1_{x_0}\psi(y_0) &= \dircon(y_0,x_0) \psi(x_0) + \dircon(y_0,x_0)\, \nabla \psi\big|_{x_0}(X) 
\end{align*}
which already contains $\nabla \psi\big|_{x_0}(X)$, that is, the full data of $j^1_{x_0}\psi$. 
\end{proof}

Given a vector bundle $E\to M$, we recall from \cite[Definitions 4.13 and 4.20]{Hairer-Singh:2023} that a \textbf{model} for an $E$-valued $n$-jet regularity structure is a pair $(\Pi^n,\Sigma^n)$ where $\Pi^n$ is a continuous linear map from $J^nE$ to the subspace $\Gamma(E)\subset \calD'(E)$ of smooth sections of $E$, and $\Sigma^n$ is a continuous section of vector bundle maps $J^nE\to J^nE$ over diffeomorphisms on $M$.

Inspired by \cite[Definition 14]{DahlqvistDiehl-Driver:2019}, a model can be called \textbf{local} if $\Pi^n:J^nE\to \Gamma_{\loc}(E)$ produces local sections and the transport map $\Sigma^n$ is defined only among nearby fibres, as for local parallelisms in \cite[Theorem 89]{DahlqvistDiehl-Driver:2019}.
The map $\Sigma^n$ can then be interpreted as a direct connection on $\GL(J^nE)$, and we have the following result.

\begin{defn}
Let $\dircon$ be a direct connection on the frame groupoid $\GL(E)$. 
For any $n\geq 0$, and for any $x_0\in M$, let $\Pi^n_{x_0}:J^n_{x_0}E \to \dircon(U_{x_0},E\big|_{U_{x_0}})$ be the map which computes the Taylor expansion of sections of $E$ induced by $\dircon$, as in Definition \ref{def:Taylor}, namely the map given by 
\begin{align*}
    \Pi^n_{x_0}(j^n_{x_0}\psi)(x) &:= \Tay_{x_0}^n \psi(x), \qquad x\in U_{x_0},  
\end{align*}
and, for any $y_0\in U_{x_0}$, let $\Sigma^n(y_0,x_0):J^n_{x_0}E \to J^n_{y_0}E$ be the Taylor direct connection of Definition \ref{def:Taylor dircon}, namely the map given by
\begin{align*} 
\Sigma^n(y_0,x_0)\ j_{x_0}^n \psi &:= j_{y_0}^n (\Tay^n_{x_0}\psi). 
\end{align*}
\end{defn}

\begin{prop} \label{prop:regularity structure}
For any $n\geq 0$, the pair of maps $(\Pi^n,\Sigma^n)$ induced by $\dircon$ gives a local model for $E$-valued $n$-jet regularity structures of transport precision $n+1$. 
\end{prop}

\begin{proof}
This follows from \cite[Proposition 4.21]{Hairer-Singh:2023}, because the generalised Taylor expansion $\Tay^n$ is a (local) admissible realisation, in the sense of \cite[Definition 3.21]{Hairer-Singh:2023}. In fact, Theorem \ref{thm:Taylor with nabla} ensures that the value of $\Tay_{x_0}^n\psi$ only depends on the jet of $\psi$ at $x_0$ (beside the transitivity defects of $\dircon$) and therefore the map $\Pi^n_{x_0}$ is well defined on jets of sections of $E$.  
Moreover it is linear on the fibres of $J^nE$, it is smooth in the base points and in the jet value, and it follows from the identity $j^n_{x_0}(\Tay^n_{x_0} \psi)=j^n_{x_0}\psi$, which holds by Theorem \ref{thm:Taylor direct connection}, that it is a right inverse to taking the jet. 

Finally, the pair $(\Tay^n,\Sigma^n)$ has transport precision $n+1$, in the sense of \cite[Definition 14]{DahlqvistDiehl-Driver:2019}. 
The proof is exactly as in \cite[Theorem 91]{DahlqvistDiehl-Driver:2019}: if we call
\begin{align*}
    \phi &:= \Tay_{x_0}^n \psi = \Pi_{x_0}^n(j_{x_0}^n \psi),  
\end{align*}
so that $\Sigma^n(y_0,x_0) j_{x_0}^n \psi = j_{y_0}^n \phi$, we then have 
\begin{align*} 
    \Pi_{y_0}^n \Sigma^n(y_0,x_0) j_{x_0}^n \psi &= \Pi_{y_0}^n \Big(j^n_{y_0}\phi\Big) =  \Tay_{y_0}^n \phi, 
\end{align*}
and therefore 
\begin{align*}
    \Pi_{x_0}^n (j_{x_0}^n\psi) - \Pi_{y_0}\Sigma^n(y_0,x_0)j_{x_0}^n\psi &= \phi - \Tay_{y_0}^n \phi = \Rem_{y_0}^n \phi. 
\end{align*}
Proposition \ref{prop:remainder} then gives the estimate of order $n+1$. 
In summary, the transitivity defects affect the explicit coefficients and the constants in the estimates, but they do not cause loss of order in transport precision. 
\end{proof}

We finally note that, if wished, a local model can be turned into a global model by extending the sections globally to $M$ using bump functions with compact support \cite{Lee:2012}.

\appendix

\section{Appendix: Jet bundles and jet prolongations}
\label{sec:jets}

In this appendix we review the basic definitions and results on jet bundles and jet prolongation of principal bundles.
Our main references are the standard books by P.\,W. Michor \cite{Michor:1980}, by D.\,J. Saunders \cite{Saunders:1989}, and by I. Kolář, P., Michor, and J. Slovák \cite{Kolar-Michor-Slovak:1993}. 


\subsection{Jets and Taylor polynomials of smooth local functions}\label{subsec:jet bundles}

Let $M$ and $M'$ be two manifolds of dimension respectively $d$ and $d'$.
Given a local map $f:M\locmap M'$ defined around a point $x\in M$, the \textbf{$n$-jet of $f$ in $x$}, denoted by $j^n_xf$, is the equivalence class of local maps from $M$ to $M'$ having the same \textbf{contact of order $n$} of $f$ in $x$, that is, the same value and the same derivatives in $x$ up to order $n$.
Jets of functions can be defined for any integer order $n\geq 0$, and $J^0(M,M')=M\times M'$. 
Denote by $J^n_x(M,M')_y$ the space of $n$-jets of local maps $f:M\locmap M'$ defined around $x$ and such that $y=f(x)$, and set $J_x^n(M,M')  :=\bigcup_{y\in M'}J_x^n(M,M')_y$. The \textbf{jet space} 
\begin{align} \label{eq:J^n(M,M')}
J^n(M,M'):=\bigcup_{(x, y)\in M\times M'}J_x^n(M,M')_y
\end{align}
has two natural projections to $M$ (the \textbf{source} $\pi^n$) and to $M'$ (the \textbf{target} $\pi^n_0$). 

As equivalence classes, jets are by definition independent of local coordinates in $M$ and in $M'$, but a representative of a jet involves derivatives which do depend on the choice of local coordinates. 
Given a choice of coordinates on $M$ and on $M'$, the identification of a local map $f:M\locmap M'$ fixing $f(x)=y$ to its local coordinates expression $\tilde{f}:\R^d\locmap \R^{d'}$, fixing $\tilde{f}(0)=0$, yields an isomorphism \cite[\S 12.6]{Kolar-Michor-Slovak:1993} between $J^n_x(M,M')_y$ and the real vector space   
\begin{align} \label{eq:L^n_{d,d'}}
L^n_{d,d'} &  := J^n_0(\R^d,\R^{d'})_0
\cong \bigoplus_{k=1}^n S^k((\R^d)^*)\otimes \R^{d'}  \cong \overline{\R_n[X_1,...,X_d]}\otimes \R^{d'}
\end{align} 
of dimension $\left[\binom{d+n}{d}-1\right]d'$ which prolongs the matrix space $L^1_{d,d'}\cong M_{d,d'}(\R)$ for $n=1$ and contains the matrix coefficents
$a^j_\alpha=\frac{1}{|\alpha|!}\partial^\alpha\tilde{f}^j(0)\in \R$.
The space $S^k((\R^d)^*)$ is that of symmetric $k$-covectors on $\R^d$, cf. \cite[\S 12.10]{Kolar-Michor-Slovak:1993} and \cite[Theorems 5.1.7 and 6.2.9]{Saunders:1989}.
The collection $(a^j_\alpha)_{j,\alpha}$, with $j=1,...,d'$ and $|\alpha|\geq 1$, together with $f(x)=(a^j_0)_j$, gives the Taylor polynomial of $f$ at $x$ in the chosen coordinates. 
Therefore $J^n(\R^d,\R^{d'}) \cong \R^d\times L^n_{d,d'}\times \R^{d'}$.
It follows that the jet space $J^n(M,M')$ is a smooth manifold \cite[Theorem 1.10]{Michor:1980}.


\subsection{Higher jet groups and higher frame bundles}\label{subsec:jet composition}

Jet composition is the key operation to express the invariance of jets under change of local coordinates. 
For two local maps $f:M\locmap M'$ around $x$ and $f':M'\locmap M''$ around $f(x)\in M'$, the \textbf{jet composition} of $j^n_xf$ and $j^n_{f(x)}f'$ is the $n$-jet  \cite[\S 12.3]{Kolar-Michor-Slovak:1993}
\begin{align} \label{eq:jet composition}
j^n_{f(x)}f' \circ j^n_xf = j^n_x(f'\circ f). 
\end{align}
The jet composition is associative, it preserves the inversion, i.e. if $h:M\locdif M'$ is invertible then $(j^n_xh)^{-1}=j^n_{H(x)}h^{-1}$, and for $M'=M$ it has a two-sided unit $1_x:=j^n_x\mathrm{Id_M}$. 

Consequently, we can read the jets of local diffeomorphisms $h:M\locdif M'$ as elements in the \textbf{space of invertible jets from $M$ to $M'$}
\begin{align} \label{eq:jet diffeo}
\mathrm{inv}J^n(M,M') & := \mbox{\parbox[t]{10cm}{set of invertible $n$-jets $j^n_xh\in J^n(M,M')$ with respect to jet composition.}}  
\end{align}
In particular, this gives the \textbf{$n$-jet group in dimension $d$} \cite[\S 12.6]{Kolar-Michor-Slovak:1993} 
\begin{align} \label{eq:GL^n_d}
\GLnd & := \mathrm{inv} L^n_{d,d} = \mathrm{inv} J^n_0(\R^d,\R^d)_0 
\\ \nonumber
& \cong \mbox{\parbox[t]{13cm}{group of $n$-jets at $0$ of local diffeomorphisms $h:\R^d\locdif \R^d$ preserving $0$, with jet composition.}} 
\end{align}
This group extends the trivial group $\GLzerod\cong \{1_{GL_d(\R)}\}$ and the general linear group $\GLund \cong GL_d(\R)$, where a jet $j^1_0h\in \GLund$ is identified with the differential $dh_0\in GL_d(\R)$. 

The jet of (the inverse of) a local chart $\varphi:\R^d\locmap M$, which represents a choice of local coordinates of $M$ at $x=\varphi(0)$, belongs to the \textbf{$n$-frame bundle of $M$} \cite[\S 12.12]{Kolar-Michor-Slovak:1993}
\begin{align} \label{eq:F^nM}
F^nM &:= \mathrm{inv}J^n_0(\R^d,M) \overset{\pi^n_0}{\longrightarrow} M
\\ \nonumber
& \cong \mbox{\parbox[t]{12cm}{set of $n$-jets at $0$ of local diffeomorphisms $\varphi:\R^d\locdif M$.}} 
\end{align}
This is a principal $\GLnd$-bundle with right action given by the jet composition $j^n_0\varphi \circ j^n_0h = j^n_0(\varphi\circ h)$, where $h$ represents a change of local coordinates on $M$ at $x$ and $\varphi\circ h$ the new local chart.
The $n$-frame bundle is trivial for a flat manifold $M=\R^d$, on which we can decompose $\varphi = \tau_x\circ h$ into a diffeomorphism $h$ preserving $0$ and the translation $\tau_x(y)=x+y$ bringing $0$ to $x$ with jet $j^n_0\tau_x=(x,\Id)$: this gives the isomorphism $F^n\R^d\cong \R^d\times \GLnd$ which identifies $j^n_0\varphi\equiv (x,j^n_0h)$. The $n$-frame bundle is also trivial for $n=0$, since $F^0M \cong M\times\{1\}$, but in general it is not. For instance, $F^1M \cong FM:=F(TM)$ is the usual frame bundle, since a jet $j^1_0\varphi\in F^1_xM$ is equivalent to the pair $(\varphi(0),d\varphi_0)\in F_xM$ with $\varphi(0)=x$, where the differential $d\varphi_0:\R^d\overset{\cong}{\to} T_xM$ gives a linear frame of $TM$ at $x$.

The simultaneous choice of a local chart $\varphi$ on $M$ and $\varphi'$ on $M'$ is then represented on jets by the two $n$-frame bundles $F^nM$ and $F^nM'$, and the change of local coordinates $h$ and $h'$ is represented by the two $n$-jet groups $GL^n_d$ and $GL^n_{d'}$. With the notations of the diagram
\begin{align} \label{fig:JetMM'}
\JetSpace
\end{align}
the group $\GLnd \times GL^n_{d'}$ acts on $L^n_{d,d'}$ from the left by $(j^n_0h,j^nh')\cdot j^n_0\tilde{f}= j^n_0(h'\circ \tilde{f}\circ h^{-1})$, and the jet space $J^n(M,M')$ can be seen as a fibre bundle over $M\times M'$ associated to $F^nM\times F^nM'$ with fibre $L^n_{d,d'}$ \cite[\S 12.12]{Kolar-Michor-Slovak:1993}, that is
\begin{align} \label{eq:J^n(M,M')=F^nMxF^nM' x_{GL^n_dxGL^_d'} L^n_dd'}
J^n(M,M') & \cong (F^n M \times F^n M') \times_{\GLnd \times GL^n_{d'}} L^n_{d,d'} & j^n_xf & 
\mapsto [j^n_0\varphi,j^n_0\varphi',j^n_0\tilde{f}]_{\varphi(0)=x}. 
\end{align}
The isomorphism holds because the local charts $\varphi$ and $\varphi'$ can be chosen arbirarly (with $\varphi(0)=x$) to define $\tilde{f}=(\varphi')^{-1}\circ f\circ \varphi$ from $f$ and viceversa.


\subsection{Higher tangent bundle}\label{subsec:jet tangent}

The local map $\bar{f}:\R^d\locmap M'$ of diagram \eqref{fig:JetMM'} reads $f:M\locmap M'$ in local coordinates on $M$. Its jet belongs to a very useful jet space, the \textbf{$n$-tangent bundle of $M'$ in dimension $d$}, also called the bundle of \textbf{$n$-velocities on $M'$ in dimension $d$} \cite[\S 12.8]{Kolar-Michor-Slovak:1993}, 
\begin{align} \label{eq:T^n_d M'}
T^n_dM' & := J^n_0(\R^d,M') \overset{\pi^n_0}{\longrightarrow} M'
\\ \nonumber
& \cong \mbox{\parbox[t]{12cm}{set of $n$-jets at $0$ of local smooth maps $\bar{f}:\R^d\locmap M'$.}} 
\end{align}
This is the fibre bundle\footnote{The higher tangent bundle $T^n_d M'$ is not a vector bundle for $n>1$, unless $M'$ is a vector space \cite[\S 12.8]{Kolar-Michor-Slovak:1993}. To have a vector bundle, one can consider the dual of the higher cotangent bundle $T^{n*}_d M'=J^n(M',\R^d)_0$, which is a vector bundle over $M'$ \cite[\S 12.14]{Kolar-Michor-Slovak:1993}.} with fibre $L^n_{d,d'}$ associated to the principal bundle $F^nM'$ \cite[\S 12.12]{Kolar-Michor-Slovak:1993}, that is
\begin{align} \label{eq:T^n_d=F^n x_{GL^n_d} L^n_d}
T^n_dM'& \cong F^nM' \times_{GL^n_{d'}} L^n_{d,d'} & j^n_0\bar{f} & 
\mapsto [j^n_0\varphi',j^n_0\tilde{f}], 
\end{align}
where the equivalence class reflects the equality $j^n_0\bar{f} = j^n_0\varphi' \circ j^n_0\tilde{f} = j^n_0(\varphi' \circ h') \circ j^n_0((h')^{-1}\circ \tilde{f})$ for every $j^n_0 h'\in GL^n_{d'}$. 

Since $T^n_dM'$ admits a left $\GLnd$-action given by $j^n_0h \cdot j^n_0\bar{f} = j^n_0(\bar{f}\circ h^{-1})$, it gives rise to a second presentation of the jet bundle $J^n(M,M')$ as the fibre bundle associated to $F^nM$ over $M$ with fibre $T^n_dM'$ (from which it inherits the bundle projection to $M'$): 
\begin{align} \label{eq:J^n=F^n x_{GL^n_d} T^n_d}
J^n(M,M') & \cong F^n M \times_{\GLnd} T^n_d M' & 
j^n_xf &\mapsto [j^n_0\varphi,j^n_0\bar{f}]_{\varphi(0)=x}, 
\end{align}
where the equivalence class now gives 
$[j^n_0 \varphi,j^n_0\bar{f}] = [j^n_0 (\varphi \circ h),j^n_0(\bar{f}\circ h)]$ for every $j^n_0 h\in GL^n_d$. 
The isomorphisms \eqref{eq:T^n_d=F^n x_{GL^n_d} L^n_d} and \eqref{eq:J^n=F^n x_{GL^n_d} T^n_d} together imply \eqref{eq:J^n(M,M')=F^nMxF^nM' x_{GL^n_dxGL^_d'} L^n_dd'}.


\subsection{Jet prolongation of fibre bundles and vector bundles}\label{subsec:jet bundle}

For a fibre bundle $\pi:E\to M$, denote by $\Gamma^\infty_\loc(M,E)$ the set of its smooth local sections $\psi:M\locmap E$, with domain denoted by $U_\psi\subset M$. The \textbf{$n$-jet bundle} of $E$ is the set \begin{align} \label{eq:J^nE}
J^nE & := J^n\Gamma^\infty_\loc(M,E) = \{ j^n_x\psi\in J^n(M,E)\ |\ \pi\circ \psi=\Id_{U_\psi} \}.  
\end{align}
The space $J^nE$ is a closed submanifold of the jet space $J^n(M,E)$ \cite[\S 12.16]{Kolar-Michor-Slovak:1993} and it is a bundle over $M$ \cite[\S 6.2]{Saunders:1989}. 
Moreover, the target jet projection $\pi^n_0$ yields a fibre bundle also on $E$.  
Forgetting the highest derivatives, one gets a filtered tower of affine bundles above $E$
\begin{align} \label{eq:jet bundle tower}
J^nE\overset{\pi_{n-1}^n}{ \rightarrow} J^{n-1}E\overset{\pi_{n-2}^{n-1}}{ \rightarrow} \cdots \overset{\pi_{1}^2}{ \rightarrow} J^1E\overset{\pi_{0}^1}{ \rightarrow} J^0E=E. 
\end{align}
 A bundle map $\phi: E\to E'$ between fibre bundles over $M$ induces a bundle map between their $n$-jet prolongations 
\begin{align}\label{eq:J^nPhi on E}
J^n\phi: J^nE\longrightarrow J^nE',\ j^n_x\psi \mapsto J^n\phi(j^n_x\psi) := j^n_x(\phi\circ \psi).
\end{align}

If $\pi:E\to M$ is a vector bundle of rank $r$, then the source jet projection $\pi^n:J^nE\to M$ is also a vector bundle \cite[\S 12.17]{Kolar-Michor-Slovak:1993}, while the target jet projection $\pi^n_0:J^nE\to E$ is an affine bundle. 
The fibre of $J^nE$ over $M$ is modelled on the vector space 
\begin{align}  
\label{eq:P^n_{d,r}}
P^n_{d,r} := T^n_d\R^r & = J^n_0(\R^d,\R^r) \cong \R^r\oplus L^n_{d,r}
\\ \nonumber  
& \cong \bigoplus_{k=0}^n S^k((\R^d)^*)\otimes \R^r  \cong \R_n[X_1,...,X_d]\otimes \R^r
\end{align} 
of dimension $\binom{d+n}{d}r$, that is, the space of $\R^r$-valued polynomials in $d$ variables, which represent the Taylor polynomials of sections of $E$ in given local coordinates.
With the notations of fig. \eqref{fig:JetMM'}, an element of $P^n_{d,r}$ can be written as the jet $j^n_0\bar{f}=\bar{f}(0)+j^n_0\tilde{f}$ where $j^n_0\tilde{f}\in L^n_{d,r}$. 
In summary, $J^nE$ is the vector bundle associated to its frame bundle $F(J^nE)$ on $M$ with fibre $P^n_{d,r}$, that is, $J^nE \cong F(J^nE) \times_{GL(P^n_{d,r})} P^n_{d,r}$. 

  
\subsection{Jet prolongation of structure groups}\label{subsec:jet group}

The structure group of the jet bundle $J^nE$ is the subgroup of $GL(P^n_{d,r})$ which hosts its transition functions, in particular the jets of transitions functions of $E$. 
If $E$ has trivializing atlas $\{U_\alpha\}$ with associated local sections $\{f_\alpha\}$ and structure group $G\subset GL_r(\R)$, the transition functions are the $G$-equivariant maps $\hat{g}_{\alpha\beta}:{U_\beta}_{|_{U_{\alpha\beta}}}\!\!\times G\to {U_\alpha}_{|_{U_{\alpha\beta}}}\!\!\times G$, defined on the intersection $U_{\alpha\beta}=U_\alpha\cap U_\beta$ of overlapping charts, which relate the local sections by $f_\alpha|_{U_{\alpha \beta}} = \hat{g}_{\alpha\beta} f_\beta|_{U_{\alpha \beta}}$.  
Reading the points $x\in U_\alpha$ in local coordinates $\vec{x}$ via the chart $\varphi=\varphi_\alpha^{-1}:\R^r\locmap U_\alpha$, and those of $U_\beta$ in local coordinates $\vec{x}'=h(\vec{x})$ via the chart $\varphi_\beta^{-1}= \varphi\circ h:\R^d\locmap U_\beta$, where $h=\varphi_\alpha\circ\varphi_\beta^{-1}:\R^r\locdif \R^r$ is the change of coordinates, leads to define the \textbf{$n$-jet prolongation of $G$ in dimension $d$} as the jet space \cite[\S 15.2]{Kolar-Michor-Slovak:1993} 
\begin{align} \label{eq:W^n_dG}
W^n_dG &\cong \mbox{\parbox[t]{13cm}{set of $n$-jets at $(0,1)\in \R^d\times G$ of $G$-equivariant local diffeomorphisms}}
\\ \nonumber
& {\parbox{13cm}{$\hat{h}:\R^d\times G\locdif \R^d\times G$ defined around $(0,1)$ and preserving $0\in \R^d$.}} 
\end{align}
Such maps $\hat{h}$ are necessarily of the form $\hat{h}(\vec{x},a)=(h(\vec{x}),g(\vec{x})a)=(h,\Id)\circ \hat{g}(\vec{x},a)$ for a local map $g:\R^d\locmap G$ defined around $0$ and $\hat{g}:\R^d\times G\locdif \R^d\times G$ given by $\hat{g}(\vec{x},a)=(\vec{x},g(\vec{x}) a)$.
Then the jet $j^n_{(0,1)}\hat{h}=(j^n_0h,1)\circ j^n_{(0,1)}\hat{g}$ is equivalent to the pair of jets $(j^n_0h,j^n_0g)\in GL^n_d\times T^n_dG$ (cf. \eqref{eq:GL^n_d} and \eqref{eq:T^n_d M'}). 
The $n$-tangent space $T^n_dG$ is a group with the pointwise product
\begin{align} \label{eq:jet group}
j^n_0g \cdot j^n_0g' = j^n_0(gg'),  
\end{align}
and $GL^n_d$ acts on $T^n_dG$ on the right by composition. The jet composition of two jets $j^n_{(0,1)}\hat{h}$ and $j^n_{(0,1)}\hat{h}'$ becomes the semidirect group law on their representative pairs
\begin{align} \label{eq:semidirect product}
\big(j^n_0h,j^n_0g\big)\cdot \big(j^n_0h',j^n_0g'\big)
& = \Big(j^n_0(h\circ h'),j^n_0((g\circ h') g')\Big). 
\end{align}
In other words, there is a group isomorphism \cite[\S 15.2]{Kolar-Michor-Slovak:1993} 
\begin{align} \label{eq:W^n_dG=GL_d^n x T^n_dG}
& W^n_dG \cong GL_d^n \ltimes T^n_dG, 
&& j^n_{(0,1)}\hat{h}=(j^n_0h,1)\circ j^n_{(0,1)}\hat{g} \equiv (j^n_0h,j^n_0g). 
\end{align}
The semidirect product law ensures also that $W^n_dG$ acts on $P^n_{d,r}$, with left action given for $j^n_{(0,1)}\hat{h} = (j^n_0h,j^n_0g)$ and $j^n_0\bar{f}=\bar{f}(0)+j^n_0\tilde{f}$ by 
\begin{align} \label{eq:W^n_dG acts on P^n_{d,r}}
j^n_{(0,1)}\hat{h} \cdot j^n_0\bar{f} &:= g(0)\bar{f}(0) + j^n_0(g\circ h^{-1})\cdot j^n_0(\tilde{f}\circ h^{-1}). 
\end{align}

Finally note that the left action of $GL^n_d$ on $T^n_dG$ used in eq. \eqref{eq:J^n=F^n x_{GL^n_d} T^n_d} extends to an action on $W^n_d G \cong GL_d^n \ltimes T^n_dG$ by 
\begin{align} \label{eq:GL^n_d acts of W^n_d G}
    j^n_0 a \cdot (j^n_0 h,j^n_0 g) = \big(j^n_0(h\circ a^{-1}),j^n_0 (g\circ a^{-1})\big).
\end{align}

\begin{lem} \label{lem:W^n_dG subgroup of GL(P^n_{d,r})}
For any subgroup $G
\subseteq GL_r(\R)$, there is a proper inclusion of groups
\begin{align} \label{eq:W^n_d GL_r(R)}
& W^n_d G = \GLnd \ltimes T^n_dG  \hookrightarrow GL(P^n_{d,r}). 
\end{align}
\end{lem}

\begin{proof}
Expressing jets in the standard basis of polynomials representing $P^n_{d,r}$, as after eq. \eqref{eq:L^n_{d,d'}}, shows that the action \eqref{eq:W^n_dG acts on P^n_{d,r}} is a matrix product and the action is linear. This defines the inclusion \eqref{eq:W^n_d GL_r(R)}. 
If we order this basis by increasing order of derivatives, we mimic the affine filtration \eqref{eq:jet bundle tower} of jets: the matrices associated to the jets $j^n_0(g\circ h^{-1})$ and $j^n_0(\tilde{f}\circ h^{-1})$ are then upper triangular. This corresponds to the results in \cite[\S 13, \S 21]{Kolar-Michor-Slovak:1993} on the filtration of jet subgroups and of their induced actions.
In other words, this means that the image of $W^n_dG$ by the inclusion \eqref{eq:W^n_d GL_r(R)} is a proper subgroup of $GL(P^n_{d,r})$. 
\end{proof}


\subsection{Jet prolongation of principal bundles}\label{subsec:jet principal bundle}

The jet bundle $J^nE$ on $M$ gives rise to the frame bundle $F(J^nE)$, which is  a principal  $GL(P^n_{d,r})$-bundle on $M$. On $J^nE$, we can consider to reduce the structure group to the subgroup $W^nGL_r(\R)$, and more generally, if $E=P\times_G \R^k$ is the vector bundle associated to a principal $G$-bundle $P$, we can consider a reduction of the structure group of $J^nE$ to $W^n_dG$.
The corresponding reduction of the principal bundle $F(J^nE)$ is obtained by considering the jet prolongation of (the inverse of) the local trivialisations $\theta_\alpha:U_\alpha\times G\locmap P$, which are local $G$-equivariant bundle automorphisms giving rise, on overlapping charts, to the transition functions $\hat{g}_{\alpha \beta}$. 
Reading these maps in local coordinates on $M$ leads to define the \textbf{principal $n$-jet prolongation} of $P$ as the bundle \cite[\S 15.3]{Kolar-Michor-Slovak:1993}
\begin{align} \label{eq:W^nP}
W^nP &= \mbox{\parbox[t]{12cm}{set of $n$-jets $j^n_{(0,1)}\theta$ of local bundle automorphisms $\theta:\R^d\times G\locmap P$ \\ above base maps $\varphi:\R^d\locmap M$,}}
\end{align}
which is naturally a bundle over $P$ with the target projection $\pi^n_0$, and becomes a bundle over $M$ with projection $\pi\circ \pi^n_0$.
The jet group $W^n_dG$ acts on $W^nP$ on the right by jet composition 
\begin{align} \label{eq:W^n_dG action on W^nP}
j^n_{(0,1)}\theta \cdot j^n_{(0,1)}\hat{h} = j^n_{(0,1)}(\theta\circ \hat{h}) = j^n_{(0,1)}(\theta\circ (h,\Id)\circ \hat{g}) , 
\end{align}
and $W^nP$ is a principal $W^n_dG$-bundle over $M$ \cite[\S 15.3]{Kolar-Michor-Slovak:1993}. 
A bundle map $\phi:P\to P'$ between two principal bundles over $M$ has a jet prolongation 
\begin{align}\label{eq:W^nPhi on P}
W^n\phi: W^nP\to W^nP',\ j^n_{(0,1)}\theta \mapsto W^n\phi(j^n_{(0,1)}\theta)=j^n_{(0,1)}(\phi\circ\theta).
\end{align}

Local trivialisations of $P$ are in bijection with local sections of $P$. This bijection survives on jet prolongations up to a choice of local coordinates: since the choice of the local chart $\varphi$ is totally arbitrary (because of the quotient by the action of $GL^n_d$), the principal prolongation $W^nP$ can be expressed as the pullback bundle of the two bundles $F^nM$ and $J^nP$ over $M$, that is \cite[\S 15.4]{Kolar-Michor-Slovak:1993} \cite[Eq. (33)]{Kolar:2009}
\begin{align} \label{eq:W^nP=F^nM x_M J^nP}
& W^nP \cong F^nM\times_M J^nP, && j^n_{(0,1)}\theta \mapsto (j^n_0\varphi,j^n_xp) 
\end{align}
where $\theta$ is the local trivialisation determined by the local section $p$, modulo a local chart $\varphi$ on $M$ such that $\varphi(0)=x$, i.e. $\theta(\vec{x},a) = p(\varphi(\vec{x}))a$ for any $\vec{x}\in \R^d$ close to $\vec{0}$ and for any $a\in G$. 
 
\begin{examples}[Jet prolongation of principal bundles] \label{ex:W^P}
Some easy examples:
\begin{enumerate}[label={(\alph*)}]
\item When $M=\R^d$, every local chart at a point $x$ is of the form $\varphi = \tau_x\circ h$, where $\tau_x(y)=x+y$ is the translation sending $0$ to $x$ and $h$ is a change of coordinates respecting $0$ \cite[\S 15.3]{Kolar-Michor-Slovak:1993}. This trivializes the frame bundle $F\R^d$. Since $j^n_0\tau_x=x$, the expression $\theta(\vec{x},a)=(\varphi(\vec{x}),g(\vec{x})a)=(x+h(\vec{x}),g(\vec{x})a)$ gives 
\begin{align*} 
W^n(\R^d\times \{1\}) & \cong F^n\R^d \cong \R^d\times GL^n_d, && j^n_{(0,1)}\theta \mapsto j^n_0\varphi \mapsto (x,j^n_0h) \\  
W^n(\R^d\times G) & \cong F^n\R^d\times_{\R^d} J^n(\R^d\times G) \\ 
& \cong \R^d\times W^n_dG && j^n_{(0,1)}\theta \mapsto (x,j^n_0h,j^n_0g) \\
W^n(F\R^d) & = W^n(\R^d\times GL_d(\R)) \cong \R^d\times W^n_dGL_d(\R). 
\end{align*}

\item For general $M$, the $n$-frame bundle $F^nM$ is not necessarily the trivial $GL^n_d$-bundle, the expression $\theta(\vec{x},a)=(\varphi(\vec{x}),g(\vec{x})a)$ gives
\begin{align*} 
W^n(M\times \{1\}) &\cong F^nM \times_M J^n(M\times \{1\}) \cong F^nM, && j^n_{(0,1)}\theta \mapsto j^n_0\varphi 
\\ 
W^n(M\times G) & \cong F^nM\times_{GL^n_d} W^n_dG \cong F^nM\times T^n_dG && 
\text{(see Lemma \ref{lem:W^n(MxG)})} 
\\ 
W^n(FM) & \cong F^nM\times_M J^n (FM) \cong F^nM\times T^n_dGL_d(\R). 
\end{align*}
The last identity follows from the fact that $J^n (FM)$ is a principal $T^n_dGL_d(\R)$-bundle on $M$, cf. \cite[Theorem 1]{Brajercik-Demko-Krupka:2014}. 
\end{enumerate}
\end{examples}

\begin{lem} \label{lem:W^n(MxG)}
For the trivial principal $G$-bundle $P=M\times G$, we have
\begin{align} \label{eq:W^n(MxG) formula 1}
W^n(M\times G) & \cong F^nM\times_{GL^n_d} W^n_dG 
\\ \label{eq:W^n(MxG) formula 2}
& \cong F^nM\times T^n_dG . 
\end{align}
\end{lem}

\begin{proof}
Let us first prove eq. \eqref{eq:W^n(MxG) formula 1}. By eq. \eqref{eq:W^nP=F^nM x_M J^nP} we have $W^n(M\times G) \cong F^nM \times_M J^n(M\times G)$. Since $J^n(M\times G)\cong J^n(M,G)$, by \eqref{eq:J^n=F^n x_{GL^n_d} T^n_d} this is isomorphic to $F^nM\times_{GL^n_d} T^n_dG$. On the other side, we have $W^n_dG \cong GL^n_d \ltimes T^n_dG$ by \eqref{eq:W^n_dG=GL_d^n x T^n_dG}. So, define a map 
$$
\Phi: F^nM \times_M (F^nM\times_{GL^n_d} T^n_dG) \longrightarrow F^nM\times_{GL^n_d} (GL^n_d \ltimes T^n_dG)
$$
by setting 
\begin{align*}
    \Phi\big(j^n_0 \varphi',[j^n_0\varphi,j^n_0 g]\big) := \big[j^n_0\varphi,(j^n_0 h,j^n_0 g)\big], 
\end{align*}
where $j^n_0 h\in GL^n_d$ is the unique element such that $j^n_0\varphi = j^n_0\varphi' \circ j^n_0 h$ and where the equivalence class on the right-hand side gives 
$$
\big[j^n_0\varphi,(j^n_0 h,j^n_0 g)\big] = \big[j^n_0(\varphi\circ a),(j^n_0 (h\circ a),j^n_0 (g\circ a))\big] \qquad \text{for any $j^n_0a \in GL^n_d$}. 
$$
The map $\Phi$ is well defined, because for a different representative $[j^n_0\varphi,j^n_0 g] = [j^n_0(\varphi \circ a),j^n_0 (g\circ a)]$, with $j^n_0 a \in GL^n_d$, we have $\varphi \circ a = \varphi'\circ h'$ with $h'=h\circ a$, and therefore 
\begin{align*}
\Phi\big((j^n_0 \varphi',[j^n_0(\varphi \circ a),j^n_0 (g\circ a)])\big) = [j^n_0(\varphi \circ a),(j^n_0 (h\circ a),j^n_0 (g\circ a))] = \big[j^n_0\varphi,(j^n_0 h,j^n_0 g)\big]. 
\end{align*}
The inverse of $\Phi$ is given by $\Phi^{-1}\big[j^n_0 \varphi,(j^n_0 h,j^n_0 g)\big] := \big(j^n_0(\varphi \circ h^{-1}),[j^n_0 \varphi,j^n_0 g]\big)$. 
The map $\Phi^{-1}$ is also well defined, because for a different representative we have 
\begin{align*}
    \Phi^{-1}\big[j^n_0 (\varphi\circ a),(j^n_0 (h\circ a),j^n_0 (\circ a)g)\big] 
    &= \big(j^n_0((\varphi\circ a) \circ (h\circ a)^{-1}),[j^n_0 (\varphi\circ a),j^n_0 (g\circ a)]\big) \\ 
    &= \big(j^n_0(\varphi \circ h^{-1}),[j^n_0 \varphi,j^n_0 g]\big). 
\end{align*}

To prove eq. \eqref{eq:W^n(MxG) formula 2}, define a map $\Psi: F^nM\times_{GL^n_d} (GL^n_d \ltimes T^n_dG) \longrightarrow F^nM \times T^n_dG$ by setting 
\begin{align*}
    \Psi\big[j^n_0 \varphi,(j^n_0 h,j^n_0 g)\big] := \big(j^n_0 (\varphi\circ h^{-1}),j^n_0(g\circ h^{-1})\big). 
\end{align*}
The map $\Psi$ is well defined because a different representative for the class $\big[j^n_0 \varphi,(j^n_0 h,j^n_0 g)\big]$ does not change the result. 
Its inverse $\Psi^{-1}$ is given by $\Psi^{-1}\big(j^n_0 \varphi,j^n_0 g\big) = \big[j^n_0 \varphi,(j^n_0 \id,j^n_0 g)\big]$.
\end{proof}



\subsection{Associated jet bundles and reductions}\label{subsec:associated jet bundle} 

Let $E=P\times_G V$ be a fibre bundle with fibre $V$ associated to a principal $G$-bundle $P\to M$. Then $J^nE$ is the fibre bundle with fibre $T^n_dV$ associated to the principal $W^n_dG$-bundle $W^nP$ \cite[\S 15.5]{Kolar-Michor-Slovak:1993}, that is, 
\begin{align} \label{eq:J^nE=W^nP x_{W^n_dG} T^n_dV}
J^nE \cong W^nP \times_{W^n_dG} T^n_dV. 
\end{align}

\begin{example} \label{W^nFE vs F(J^nE)}
Let $E$ be a vector bundle with typical fibre $\R^r$ and structure group $GL_r(\R)$. We can consider three interesting jet bundles by reducing the structure group of $J^nE$.  
\begin{enumerate}[label={(\alph*)}]
\item 
Since the jet bundle $J^nE$ is a vector bundle with typical fibre $P^n_{d,r}=T^n_d(\R^r)$, its frame bundle $F(J^nE)$ is a principal $GL(P^n_{d,r})$-bundle and we have
$$
J^nE \cong F(J^nE)\times_{GL(P^n_{d,r})} P^n_{d,r}.  
$$
\item
The frame bundle $FE$ is a principal $GL_n(\R)$-bundle such that $E\cong FE\times_{GL_n(\R)} \R^r$. Then, by \eqref{eq:J^nE=W^nP x_{W^n_dG} T^n_dV}, we have $J^nE \cong W^nFE \times_{W^n_dGL_r(\R)} P^n_{d,r}$. This means that the structure group $GL(P^n_{d,r})$ of $J^nE$ can be reduced to its subgroup $W^n_dGL_r(\R)$, therefore we have
\begin{align} \label{eq:F(J^nE) vs W^nFE} 
F(J^nE) \cong W^nFE \times_{W^n_dGL_r(\R)} GL(P^n_{d,r}). 
\end{align}
\item 
If the structure group of $E$ can be reduced to $G\subset GL_r(\R)$, so that there is a principal $G$-bundle $P$ such that $E\cong P\times_G \R^r$, then \eqref{eq:J^nE=W^nP x_{W^n_dG} T^n_dV} says that $J^nE \cong W^nP \times_{W^n_dG} P^n_{d,r}$. In other words we can further reduce the structure group of $J^nE$ to $W^n_dG\subset W^n_dGL_r(\R)$, and therefore we have 
\begin{align*} 
  W^nFE & \cong W^nP \times_{W^n_dG} W^n_dGL_r(\R),
  \\ 
  F(J^nE) & \cong W^nP \times_{W^n_dG} GL(P^n_{d,r}). 
\end{align*}
\end{enumerate}
Note that by Lemma \ref{lem:W^n_dG subgroup of GL(P^n_{d,r})} there is a proper inclusion of groups $W^n_dG \subsetneq W^n_d GL_r(\R) \subsetneq GL(P^n_{d,r})$, and therefore a proper inclusion of principal bundles $W^nP \subsetneq W^nFE \subsetneq F(J^nE)$. 
\end{example}


\subsection{Jets and covariant derivatives}\label{subsec:jets and covariant derivatives}

Let $M$ be equipped with a torsion-free linear connection $\nabla^M$ on $TM$, and let $E\to M$ be a vector bundle of finite rank equipped with a linear connection $\nabla$.
The \textbf{higher covariant derivatives} 
\begin{align}\label{eq:higher covariant derivatives}
    \nabla^n:\Gamma(E)\to \Gamma(T^*M^{\otimes n}\otimes E)
\end{align}
are defined as follows, cf. \cite[Ch.III, Section 2]{Kobayashi-Nomizu:1963} or \cite[eq. (2.1)]{JafarpourLewis:2014}. 
Let $\psi$ be a section of $E$. For $n$=0 set $\nabla^0\psi=\psi$. For $n=1$ let $\nabla^1\psi = \nabla \psi$ be the given covariant derivative. For $n>1$, first recall that the linear connection $\nabla^M$ on $TM$ induces a linear connection $\nabla^*$ on $T^*M$ such that the Lie derivative $\mathcal{L}_Y$ satisfies the Leibniz rule
\begin{align*}
    \mathcal{L}_Y\langle \alpha,X\rangle = \langle \nabla^*_Y\alpha,X \rangle + \langle \alpha,\nabla_Y X \rangle, 
\end{align*}
where $\alpha \in \Omega^1(M)$ and $X,Y\in \mathfrak{X}(M)$. Then extend the tensor product $\nabla^*\otimes \nabla$ to a linear connection $\nabla^{(n+1)}$ on $T^*M^{\otimes n}\otimes E$ satisfying another Leibniz rule
\begin{align*}
\nabla^{(n+1)}(\alpha_1\otimes \cdots \otimes \alpha_n\otimes \psi)
&= \sum_{j=1}^n  \alpha_1\otimes \cdots\otimes \nabla^*\alpha_j\otimes \cdots  \otimes \alpha_n\otimes \psi\ +\ \alpha_1\otimes \cdots \otimes \alpha_n\otimes \nabla \psi, 
\end{align*} 
for $\alpha_1,...,\alpha_n\in \Omega^1(M)$. Finally, compose these connections into
\begin{align*}
\nabla^{n} \psi:=\nabla^{(n)}\nabla^{(n-1)}\cdots \nabla^{(2)}\nabla \psi.
\end{align*}
This section can be contracted with $n$ vector fields and evaluated at a point $x$, or evaluated at $x$ and contracted with $n$ tangent vectors at $x$, since we have 
\begin{align*}
\nabla^{n} \psi\big|_x(X_x,Y_x,...,Z_x) &= \nabla^{n} \psi(X,Y,...,Z)\big|_x. 
\end{align*}
Note that, when $\nabla^n \psi$ is evaluated along a geodesic $x(t)=\exp_{x_0}(tX)$ and contracted $n$ times with the velocity vector $\dot{x}(t)$, applying successively \cite[Ch.III, Proposition 2.12]{Kobayashi-Nomizu:1963} or \cite[Exercice 4.5, p.54]{Lee:1997} (or directly \cite[Definition 56]{DahlqvistDiehl-Driver:2019}), we have
\begin{align} \nonumber
\nabla^2 \psi\big|_{x(t)} (\dot{x}(t),\dot{x}(t)) &= \nabla_{\dot{x}(t)} \big(\nabla_{\dot{x}(t)} \psi\big)\big|_{x(t)} - \nabla \psi\big|_{x(t)} \big(\nabla^M_{\dot{x}(t)} \dot{x}(t)\big) 
= \nabla_{\dot{x}(t)} \big(\nabla_{\dot{x}(t)} \psi\big)\big|_{x(t)} 
\\ \label{eq:nabla^n}
\nabla^n \psi\big|_{x(t)} \big(\dot{x}(t)^{\otimes n}\big) &= 
\nabla_{\dot{x}(t)}\big(\nabla_{\dot{x}(t)}\cdots (\nabla_{\dot{x}(t)} \psi)\cdots \big)\big|_{x(t)} \qquad n\geq 2.
\end{align}

The \textbf{symmetrised covariant derivatives} $D_\nabla^n:\Gamma(E) \to \Gamma(S^n(T^*M)\otimes E)$ are then defined as 
\begin{align} \label{eq:symmetrised covariant derivatives}
 	D_\nabla^n \psi:=[\Sym^n \otimes \Id_E]\left(\nabla^{n} \psi\right), 
\end{align}
where $\Sym^n: T^*M^{\otimes n} \to S^n(T^*M)=(T^*M^{\otimes n})^{S_n}$ is the usual symmetrizing map given by 
$$
\Sym^n(\alpha_1\otimes \cdots \otimes \alpha_n):= \frac{1}{n!}\,\sum_{\sigma \in \Sigma_n} \alpha_{\sigma^{-1}(1)}\otimes \cdots \otimes \alpha_{\sigma^{-1}(n)}.
$$
The operators $D_\nabla^n$ provide a splitting of the jet bundle exact sequence \cite[Ch. IV, Theorem 7]{Palais:1965}
\begin{align*}
    0 \to S^n(T^*M)\otimes E \longrightarrow J^nE \overset{\pi^n_{n-1}}{\longrightarrow} J^{n-1}E \to 0.
\end{align*}
As a consequence, see \cite[Ch. IV, Corollary, p.90]{Palais:1965}, the given torsion-free linear connection $\nabla^M$ on $TM$, and the linear connection $\nabla$ on $E$, provide a canonical isomorphism of vector bundles \begin{align}\label{eq:JE=ST*MxE}
    J^nE \cong \bigoplus_{k=0}^n S^k(T^*M)\otimes E
\end{align}
under which jets are represented by the symmetrised covariant derivatives as 
\begin{align}\label{eq:jets=symmetrised covariant derivatives}
    j_{x_0}\psi \mapsto [x_0,\psi(x_0),D_\nabla \psi\big|_{x_0},...,D_\nabla^n \psi\big|_{x_0}]. 
\end{align}
This means that jets are represented by the symmetrised covariant derivatives $D_\nabla^k\psi\big|_{x_0}(X_1,...,X_k)$ evaluated in all (small) directions $X_i\in T_{x_0}M$, for any $k=0,...,n$. 


Eq. \eqref{eq:jets=symmetrised covariant derivatives} gives an exact representation of jets. However one can embed jets in all tensors via $\iota:J^nE \hookrightarrow \bigoplus_{k=0}^n (T^*M)^{\otimes k} \otimes E$ and represent them directly with the higher covariant derivatives $\nabla^k\psi$, that is, 
\begin{align}\label{eq:jets=covariant derivatives + redundances}
    j_{x_0}\psi \mapsto [x_0,\psi(x_0),\nabla \psi\big|_{x_0},...,\nabla^n \psi\big|_{x_0}],
\end{align}
with the agreement that the representing covariant derivatives are not independent of one another, but satisfy some Ricci commutation identities involving curvatures terms of $\nabla$, cf. \cite[Lemma 4]{Jentsch:2015} and \cite[eq. (2.4)]{JafarpourLewis:2014}.




\addcontentsline{toc}{section}{References}
\bibliography{references}
\bibliographystyle{alpha}


\bigskip
\bigskip
\bigskip

\noindent Sara Azzali\\
Dipartimento di Matematica\\
Universit\`a degli Studi di Bari Aldo Moro\\
Via Edoardo Orabona 4\\
70125 Bari, Italy\\
\texttt{sara.azzali@uniba.it}
\bigskip

\noindent Alessandra Frabetti\\
Institut Camille Jordan\\
Universit\'e Claude Bernard Lyon 1\\
43 Boulevard du 11 novembre 1918\\
F-69622 Villeurbanne Cedex, France\\
\texttt{frabetti@math.univ-lyon1.fr}
\bigskip

\noindent Youness Bouta\"ib\\
Department of Mathematical Sciences\\
Institute for Financial and Actuarial Mathematics\\
University of Liverpool\\
Liverpool, L69 7ZL, United Kingdom\\
\texttt{youness.boutaib@liverpool.ac.uk}
\bigskip

\noindent Sylvie Paycha\\
Institut f\"ur Mathematik\\
Campus Golm, Haus 9\\
Karl-Liebknecht-Str. 24-25 D-14476 Potsdam\\
\texttt{paycha@math.uni-potsdam.de}

\end{document}